\documentclass[1 [leqno,11pt]{amsart}
\usepackage{amssymb, amsmath,latexsym,amsfonts,amsbsy, amsthm,mathtools,graphicx,CJKutf8,CJKnumb,CJKulem,color}
\usepackage{float}
\usepackage{hyperref}
\usepackage{verbatim}
\renewcommand{\theequation}{\thesection.\arabic{equation}}
\numberwithin{equation}{section}

\allowdisplaybreaks
\let\f=\frac

\newcommand{\beq}{\begin{equation}}
\newcommand{\eeq}{\end{equation}}
\newcommand{\ben}{\begin{eqnarray}}
\newcommand{\een}{\end{eqnarray}}
\newcommand{\beno}{\begin{eqnarray*}}
\newcommand{\eeno}{\end{eqnarray*}}

\newtheorem{theorem}{Theorem}[section]

\newtheorem{lemma}[theorem]{Lemma}
\newtheorem{proposition}[theorem]{Proposition}

\newtheorem{remark}[theorem]{Remark}
\newtheorem{Theorem}{Theorem}[section]

\newtheorem{Corollary}[Theorem]{Corollary}
\newtheorem{Remark}[Theorem]{Remark}
\begin{document}
\begin{CJK*}{UTF8}{gkai}
\title[Taylor-Couette flow]{Nonlinear asymptotic stability and transition\\ threshold for 2D Taylor-Couette flows\\ in Sobolev spaces}

\author{Xinliang AN}
\address{Department of Mathematics, National University of Singapore, Singapore}
\email{matax@nus.edu.sg}

\author{Taoran HE}
\address{Department of Mathematics, National University of Singapore, Singapore}
\email{taoran\underline{~}he@u.nus.edu}

\author{Te LI}
\address{Department of Mathematics, National University of Singapore, Singapore}
\email{matlit@nus.edu.sg}

\date{\today}

\maketitle

\begin{abstract}
In this paper, we investigate the stability of the 2-dimensional (2D) Taylor-Couette (TC) flow for the incompressible Navier-Stokes equations. The explicit form of velocity for 2D TC flow is given by $u=(Ar+\frac{B}{r})(-\sin \theta, \cos \theta)^T$ with $(r, \theta)\in [1, R]\times \mathbb{S}^1$ being an annulus and $A, B$ being constants. Here, $A, B$ encode the rotational effect and $R$ is the ratio of the outer and inner radii of the annular region. Our focus is the long-term behavior of solutions around the steady 2D TC flow. While the laminar solution is known to be a global attractor for 2D channel flows and plane flows, it is unclear whether this is still true for rotating flows with curved geometries. In this article, we prove that the 2D Taylor-Couette flow is asymptotically stable, even at high Reynolds number ($Re\sim \nu^{-1}$), with a sharp exponential decay rate of $\exp(-\nu^{\frac13}|B|^{\frac23}R^{-2}t)$ as long as the initial perturbation is less than or equal to $\nu^\frac12 |B|^{\frac12}R^{-2}$ in Sobolev space.  The powers of $\nu$ and $B$ in this decay estimate are optimal. It is derived using the method of resolvent estimates and is commonly recognized as the enhanced dissipative effect. Compared to the Couette flow, the enhanced dissipation of the rotating Taylor-Couette flow not only depends on the Reynolds number but also reflects the rotational aspect via the rotational coefficient $B$. The larger the $|B|$, the faster the long-time dissipation takes effect. We also conduct space-time estimates describing inviscid-damping mechanism in our proof. To obtain these inviscid-damping estimates, we find and construct a new set of explicit orthonormal basis of the weighted eigenfunctions for the Laplace operators corresponding to the circular flows. These provide new insights into the mathematical understanding of the 2D Taylor-Couette flows.

\end{abstract}

\tableofcontents
\section{\textbf{Introduction}}

Reynolds's famous experiment \cite{Rey} inspired the study of hydrodynamic stability at high Reynolds number. In this regime, the laminer flows could become unstable and transition to turbulence  \cite{DW,SH,TTRD,Ya}.

With the low Reynolds number, Serrin \cite{Se} demonstrated that all equilibria of the forced Navier-Stokes equation on bounded domains are linearly stable. At high Reynolds number, even in the absence of boundaries, the viscosity can significantly complicate the linear problem.

In this paper, we study the 2D (two-dimensional) incompressible Navier-Stokes (NS) equations:
\begin{align}
\label{full nonlinear equation}\left\{
\begin{aligned}
&\partial_tv-\nu\Delta v+v\cdot \nabla v+\nabla p=0,\quad\text{div } v=0,\\
&v(0,x)=v_0(x),
\end{aligned}
\right.
\end{align}
where $v=(v_1,v_2)\in\mathbb{R}^2$ is the fluid velocity, $x=(x_1,x_2)\in\Omega$ represents the space variables, $\Omega\in\mathbb{R}^2$ is an annular region and $t\geq0$ represents the time variable. The unknowns in the equation are the velocity field $v(t,x)=(v^1(t,x),v^2(t,x))$ and the pressure $p(t,x)\in\mathbb{R}$. The constant $\nu>0$ is known as the kinematic viscosity, which is very small in our paper. And the Reynolds number $Re$ is proportional to the inverse of $\nu$.

\begin{comment}
Turbulence occurs when strong shear exists in between two or more adjacent layers of fluid flowing with different velocities.  This kinetic energy is then gets transported to small-scale or perturbations, that meander away from the fluid flow, and have a tendancy to rotate, or produce vorticity (defined as the curl of the velocity vector). If taking the 2D curl to the above Navier-Stokes equations, then \eqref{full nonlinear equation} is transferred to
\begin{align}
\label{full nonlinear equation ns vor}
&\partial_t\omega-\nu\Delta\omega+v\cdot \nabla \omega=0.
\end{align}
Here $\omega=\text{curl }v=\partial_1v^2-\partial_2v^1$ is the vorticity field. For 2D case, it is a scalar. And \eqref{full nonlinear equation ns vor} is called the vorticity equation formulation. In this paper, we will use \eqref{full nonlinear equation ns vor} to understand the asymptotics of solutions to \eqref{full nonlinear equation}.
\end{comment}
In the 2D case, the vorticity field 
\begin{align}
\label{scalar}\omega=\text{curl }v=\partial_1v^2-\partial_2v^1    
\end{align}
is a scalar. By taking the 2D curl of the Navier-Stokes equations, equation \eqref{scalar} can be transformed into its vorticity formulation:
\begin{align}
\label{full nonlinear equation ns vor}
&\partial_t\omega-\nu\Delta\omega+v\cdot \nabla \omega=0.
\end{align}
In this paper, we employ \eqref{full nonlinear equation ns vor} to study the asymptotics of solutions to \eqref{full nonlinear equation} around a Taylor-Couette (TC) flow.

\subsection{Taylor-Couette flow}
\begin{comment}
TC flows were initially used to measure the viscosity of fluids, and to study the stability of fluids inside cylinders.  The setup essentially consists of two coaxial cylinders of different radii, and the fluid is filled between the two.  One of these cylinders is then rotated and the setup is studied using various tools and means of measurement.  

Although demonstrably simple and trivial at the first glance, the flow proves to be extremely hard to explain reasonably, and has been the subject of active experimental, theoretical as well as numerical study for a very long time.  This flow is still actively studied and a lot of unanswered questions till remain.  The reason why the TC flow is a good subject of study in circular geometries comes from the fact that it is well studied and has a lot of literature \cite{CI,F,K,Po}, but rigorous mathematical proofs are still insufficient, even for the two-dimensional TC flow. 

Taylor-Couette flow describes a steady circular flow of viscous fluid bounded between two rotating infinitely long coaxial cylinders and has wide applications ranging from desalination to viscometric analysis. 
\end{comment}
TC flow describes a steady circular solution of the viscous fluid bounded between two rotating infinitely long coaxial cylinders. {The understanding of solutions' asymptotics around a TC flow has board applications}, including desalination, magnetohydrodynamics and viscometric analysis. Despite being a simple type of rotating solutions, perturbation of the TC flow has proven to be a challenging subject and has been extensively studied experimentally, theoretically, and numerically for a long time. Many questions related to TC flow remain unanswered, making it active research field in fluid mechanics \cite{CI,F,K,Po}. However, rigorous mathematical proofs in this field are still insufficient, even for the 2D case.

In the following, we adopt the convention $r = |x|$ and denote the radial vorticity $\omega(x) = \omega(r)$ and stream function $\phi(x) = \phi(r)$. In 2D, the vorticity and the velocity field then take the form
\begin{align}
\label{radical vor}
\left\{
\begin{aligned}
&\omega(x_1,x_2)=\omega(r)=\Delta\phi=\phi''(r)+\f{1}{r}\phi'(r),\\
&v(x_1,x_2)=\left(
  \begin{array}{ccc}
   -\partial_{x_2}\phi \\
   \partial_{x_1}\phi \\
  \end{array}
\right)=\left(
  \begin{array}{ccc}
   -x_2 \\
   x_1 \\
  \end{array}
\right)\f{\phi'(r)}{r}.
\end{aligned}
\right.
\end{align}

We first derive the explicit form of the steady TC flow. First noting that if $\omega=const$, then $\omega$ is a solution to the 2D NS equation \eqref{full nonlinear equation ns vor}. For this case, the stream function $\phi(r)$ can be determined from \eqref{radical vor}, which yields
\begin{align}
\label{Taylor-Couette flow stream function}\phi''(r)+\f{1}{r}\phi'(r)=const.
\end{align}
By employing the polar coordinates,  one can relabel  $v(x_1,x_2)$ and $\omega(r)$ as $U(r,\theta)$ and $\Omega(r)$, respectively. This allows us to solve equation \eqref{Taylor-Couette flow stream function} and obtain the expressions for $U$ and $\Omega$ as follows:
\begin{align}
\label{Taylor-Couette flow}
&U(r,\theta)=\left(
  \begin{array}{ccc}
   U^1 \\
   U^2\\
  \end{array}
\right)=\left(
  \begin{array}{ccc}
   -\sin\theta \\
   \cos\theta\\
  \end{array}
\right)(Ar+\f{B}{r}),\quad
\Omega(r)=2A.
\end{align}
Here $A,B$ are constants and spatial variables $(r, \theta)$ belong to a domain $\mathcal{D}=[1,R]\times \mathbb{S}^1$. The function $U(r,\theta)$ in \eqref{Taylor-Couette flow} is hence a steady state for the 2D incompressible NS equation \eqref{full nonlinear equation}, which is commonly called the  \textit{Taylor-Couette flow}.

\subsection{Hydrodynamic stability at high Reynolds number}
\begin{comment}
    A very interesting thing to note about the Navier-Stokes equation is that it is dependent on initial conditions, and small perturbations in the initial conditions leads to huge changes in the flow, due to the nature of the nonlinear term of the equation.  This phenomenon is more accurately expressed as:
 \begin{itemize}
     \item  on one hand, theoretical analysis shows that some laminar flows are linearly stable for any Reynolds number such as the planar Couette flow \cite{DW};
     \item  on the other hand, the experiments show that they could be unstable and transition to turbulence for small but finite perturbations at high Reynolds number for the planar Couette flow \cite{DHB,TA};
 \end{itemize}
  which is referred to as \textit{subcritical transition}. The regime of subcritical transition is only true for larger values of Reynolds numbers. This leads to the central problem in fluid mechanics: hydrodynamic stability at high Reynolds number.
  \end{comment}
Reynold’s well-known experiment \cite{Rey} revealed that small perturbations can cause significant change of the flow due to the nonlinear nature of the Navier-Stokes equations.  In experiments \cite{DHB,TA} for the TC flow, a small perturbation could lead to instability and transition to turbulence at high Reynolds numbers. This phenomenon is called the \textit{subcritical transition}, which is a central topic in fluid mechanics.

Numerous efforts have been made to comprehend the subcritical transition mechanism \cite{Cha}. Kelvin initially suggested that as the Reynolds number $Re\to\infty$, the basin of attraction for the laminar flow diminishes, allowing the flow to become nonlinearly unstable for small but finite perturbations \cite{Kel}. In \cite{TTRD}, Trefethen, Trefethen, Reddy and Driscoll proposed a way to determine the threshold amplitude by considering a perturbation of size $Re^{-\beta}$ with $\beta\geq0$ and $Re\to +\infty$.  In other words, they aimed to identify the lowest possible value of $\beta$, which could result in a transition to turbulence with a perturbation of size $O(Re^{-\beta})$. Bedrossian, Germain, and Masmoudi \cite{BGM-bams} presented a rigorous mathematical formulation of this approach using the fact that $Re^{-1}\sim\nu$.

\smallskip

{\it
Given a norm $\|\cdot\|_X$, we now hope to determine a nonnegative number $\beta=\beta(X)$ so that
\beno
&&\|u_0\|_X\le \nu^{\beta}\Longrightarrow  {stability},\\
&&\|u_0\|_X\gg \nu^{\beta}\Longrightarrow  {instability}.
\eeno
}
The exponent $\beta$ is referred to as the \textit{transition threshold} in the applied literature. 
\begin{comment}
Recently, for the 2D incompressible NS equation, there have been interesting and important developments and results on the regime of subcritical transition, such as Couette flow \cite{BMV,BWV,BGM-bams,MZ-2,MZ-1} (in the absence of boundaries) and
\cite{CLWZ-2D-C,BH} (in the presence of boundaries), Lamb-Oseen vortices \cite{Ga0,GW,LWZ-O,Ga,GV}, periodic Kolmogorov flow \cite{WZZ,WZ}.
\end{comment}
Significant advancements and findings related to the subcritical transition regime for the 2D incompressible NS equations have emerged recently. These developments include researches on Couette flow without boundaries \cite{BMV,BWV,BGM-bams,MZ-2,MZ-1} and with boundaries \cite{CLWZ-2D-C,BH}, on Lamb-Oseen vortices \cite{Ga0,GW,LWZ-O,Ga,GV} and on periodic Kolmogorov flow \cite{WZZ,WZ}.

\subsection{Main theorems}

The aim of this paper is to explore both linear and nonlinear stability mechanism of the 2D Taylor-Couette flow using rigorous analytical tools. Specifically, we prove nonlinear asymptotic stability of the TC flows for the 2D NS equations. Furthermore, we identify and trace a  critical parameter reflecting the rotating properties of this system.

Our focus in this paper is on the stability of the TC flow, hence we study the below perturbation equation.  Setting $w = \omega - \Omega$, $u = v - U$, with 
$\varphi$ being the stream function satisfying $\Delta\varphi=w$ and $u=(-\partial_2\varphi,\partial_1\varphi)$, we concert \eqref{full nonlinear equation ns vor} to polar coordinates:
\begin{align}
\label{pertubation of NS vor polar coordinates}
\left\{
\begin{aligned}
&\partial_tw-\nu(\partial_r^2+\f{1}{r}\partial_r+\f{1}{r^2}\partial_{\theta}^2)w+(A+\f{B}{r^2})\partial_{\theta}w+\f{1}{r}(\partial_r\varphi\partial_{\theta}w-\partial_{\theta}\varphi\partial_rw)=0,\\
&(\partial_r^2+\f{1}{r}\partial_r+\f{1}{r^2}\partial_{\theta}^2)\varphi=w,\quad w|_{r=1,R}=0\quad \text{with }(r,\theta)\in[1,R]\times\mathbb{S}^1 \text{ and } t\geq0.
\end{aligned}
\right.
\end{align}
\begin{Remark}
      Note that constants $A$ and $B$ in \eqref{Taylor-Couette flow}, \eqref{pertubation of NS vor polar coordinates} and constant $R$ in \eqref{pertubation of NS vor polar coordinates} will serve as parameters in our later arguments. In particular, $B$ will correspond to the rotating effect.
\end{Remark}

Our main results are summarized below.

\begin{theorem}

\label{main statement}
The fully nonlinear two-dimensional incompressible Navier-Stokes equations  exhibit asymptotic stability around the Taylor-Couette flow under perturbations of size $\nu^\f12 |B|^{\f12}R^{-2}$ in the $H^1$ space, with the vorticity controllable at any time  by the initial data. 
\end{theorem}
The more precise mathematical statements are given in the following two theorems. The first main result can be directly derived from Proposition \ref{linear exp decay} and Lemma \ref{zero frequency part}, and describes the asymptotic stability of the linearized TC flow.
\begin{theorem}\label{main theorem-2}
There exist constants $C,c>0$ being independent of $\nu,A,B,R$, such that the solution $w$ for the linear system
\begin{align}
\label{linear-pertur}&\partial_tw-\nu(\partial_r^2+\f{1}{r}\partial_r+\f{1}{r^2}\partial_{\theta}^2)w+(A+\f{B}{r^2})\partial_{\theta}w=0,\quad w|_{t=0}=w(0),
\end{align}
exists globally in time and for any $t\geq0$. Furthermore, the following stability estimates hold
	\begin{align*}
& \|w(t)-\bar{w}(t)\|_{L^2_{\theta}L^2_r}\leq Ce^{-c\nu ^{\f13}|B|^{\f23}R^{-2}t}\|w(0)\|_{L^2_{\theta}L^2_r}, \quad \|\bar{w}(t)\|_{L^2_{\theta}L^2_r}\leq C\|w(0)\|_{L^2_{\theta}L^2_r},
\end{align*}
where $\bar{w}(t)=\f{1}{2\pi}\int_0^{2\pi}w(t,r,\theta)d\theta$.
\end{theorem}
The second main result describes the asymptotic stability of the nonlinear Navier-Stokes equations around the TC flow. Reynolds experiments tell us that, with small viscosity coefficients, even tiny initial perturbations could cause the flow to be chaotic. Therefore, in order to establish the nonlinear asymptotic stability of Navier-Stokes equations, the initial perturbations must be restricted to a certain range with respect to the viscosity coefficient $\nu$. The specific formulation of the result is as follows.
\begin{theorem}\label{main-transition threshold}Given $0<\log R\lesssim\nu^{-\f13}|B|^{\f13}$. There exist constants $\nu_0$ and $c_0, C,c'>0$ independent of $\nu,A,B,R$  such that for any $0<\nu\leq \nu_0$, if the initial data satisfies 
\begin{align*}
E(0)=R\|w(0)\|_{L^2_{\theta}H^1_r}+R^{-2}(\log R)^{-\f32}\|r^2w(0)\|_{L^2_{\theta}L^2_r}+R^{3}\|\f{w(0)}{r^3}\|_{L^2_{\theta}L^2_r}\le c_0\nu^\f12 |B|^{\f12}R^{-2},
\end{align*}
then the solution $w(t,r,\theta)$ to the system \eqref{pertubation of NS vor polar coordinates} is global in time. Moreover, the following stability estimates hold
\begin{align*}
   & \|w(t)-\bar{w}(t)\|_{L^2_{\theta}L^2_r}\le Cc_0e^{-c'\nu^{\f13}|B|^{\f23}R^{-2}t}E(0),
\quad\|\bar{w}(t)\|_{L^2_{\theta}L^2_r}\le Cc_0E(0),\quad \textrm{for any}\ t\geq0,
\end{align*}
where $\bar{w}(t)=\f{1}{2\pi}\int_0^{2\pi}w(t,r,\theta)d\theta$.
\end{theorem}

By setting the radius $R$ to a constant value, such as $R=2$, we can immediately derive the following corollary from Theorem \ref{main-transition threshold}.

\begin{theorem}\label{main corollary}There exist constants $\nu_0$ and $c_0, C,c'>0$ independent of $\nu,A,B$ so that for any $0<\nu\leq \nu_0$ and $|B|\geq\nu_0$, if 
\begin{align*}
E(0)=\|w(0)\|_{L^2_{\theta}H^1_r}\le c_0\nu^\f12 |B|^{\f12},
\end{align*}
then the system \eqref{pertubation of NS vor polar coordinates} admits a global-in-time solution $w(t,r,\theta)$ safisfying the following stability estimates
\begin{align*}
   & \|w(t)-\bar{w}(t)\|_{L^2_{\theta}L^2_r}\le Cc_0e^{-c'\nu^{\f13}|B|^{\f23}t}E(0),
\quad\|\bar{w}(t)\|_{L^2_{\theta}L^2_r}\le Cc_0E(0),\quad \textrm{for any}\ t\geq0,
\end{align*}
where $\bar{w}(t)=\f{1}{2\pi}\int_0^{2\pi}w(t,r,\theta)d\theta$.
\end{theorem}

\subsection{Difficulties, new ingredients and the sketch of the proof}

\subsubsection{Rotational effect and enhanced dissipation}  

Compared to Couette flow, the Taylor-Couette flow involves additional coefficients $A$ and $B$ that account for rotational effects. Both Theorem \ref{main theorem-2} and Theorem \ref{main-transition threshold} demonstrate that for both linearized and nonlinear equations, only the rotation coefficient $B$ affects the stability of the system, and the energy of vorticity is irrelevant to the coefficient $A$.

The impact of the rotation coefficient $B$ on stability is mainly manifested in  the so-called \textit{enhanced dissipation} effect.  The decay rate of the heat equation is $e^{-\nu  R^{-2} t}$, but in the presence of $B$, the dissipation rate becomes faster and can be described by $e^{-c\nu ^{\f13}|B|^{\f23}R^{-2}t}$ (Notice that we are interested in the regime when $0<\nu\ll 1$). This implies that while the heat equation exhibits decay only after the time scale of $\nu^{-1}$, the antisymmetric part $ik\f{B}{r^2}$ in the linearized equation results in the energy decay after a shorter time scale of $\nu ^{-\f13}|B|^{-\f23}R^{2}$. This phenomenon is known as the enhanced dissipation effect. A shear or a diffusion averaging mechanism to trigger it were investigated in \cite{L,RY}. Here we prove this effect in a different setting with rotations.

Theorem \ref{main theorem-2} and Theorem \ref{main-transition threshold} both indicate that the enhanced dissipation effect becomes stronger as $|B|$ increases, which is quantified by a faster asymptotic decay rate $e^{-c\nu ^{\f13}|B|^{\f23}R^{-2}t}$. This enhanced dissipation has a timescale of $t\sim\nu ^{-\f13}|B|^{-\f23}R^{2}$ around the Taylor-Couette flow, with $t\sim\nu ^{-\f13}$ being consistent with Couette flow \cite{BGM-bams,BH,BMV,BWV,CLWZ-2D-C,MZ-1,MZ-2}.

The enhanced dissipation decay rate $e^{-c\nu ^{\f13}|B|^{\f23}R^{-2}t}$ for the TC flow in this paper is optimal. In Section \ref{Section Resolvent estimate}, we derived the sharp resolvent estimates, which are also known as the optimal pseudospectral bound. It is important to note that although the enhanced dissipation decaying rates for the TC flow ($\nu ^{\f13}|B|^{\f23}R^{-2}$) and Couette flow ($\nu ^{\f13}$) have the same power of $\nu$, it does not suggest that the corresponding enhanced dissipation can be directly derived through a scaling argument from results of Couette flow. This is because the (normalized) antisymmetric parts of linearized operators for these two flows are still different. Furthermore, with the aid of dimensional analysis or scaling transformations, it can be shown that the sum of the power of $\nu$ and $B$ is always $1$. Thus, the factor $|B|^{\f23}$ is also  optimal.

Additionally, inspired by the approach used to handle Oseen vortices \cite{GW,LWZ-O,Ga}, in \cite{AHL-1} the authors employed self-similar variables and derived an enhanced dissipation decay rate that is independent of the outer radius $R$.

\subsubsection{Comparison with Couette flow from an operator perspective.}  

Previous studies have investigated plane Couette flow with $(x,y)\in\mathbb{T}\times\mathbb{R}$ ( see \cite{BMV,BWV,BGM-bams,MZ-2,MZ-1}) and Couette flow in a finite channel with $(x,y)\in\mathbb{T}\times[-1,1]$ (see \cite{CLWZ-2D-C,BH}). In these works, the spatial variable $x$ corresponding to the non-shear direction was defined on a torus. However, for the TC flow in this paper, the spatial domain is an annulus region $(r,\theta)\in[1,R]\times\mathbb{S}^1$.

As both $x$ in Couette flow and $\theta$ in TC flow are defined on a torus or $\mathbb{S}^1$, it is natural to apply Fourier transform on the $x$ or $\theta$ direction. The corresponding linearized equations around the Couette flow and TC flow are given as below.
\begin{align*}
\text{Couette flow: }&\partial_t\hat{w}_k-\nu(\partial_y^2-k^2)\hat{w}_k+iky\hat{w}_k=0,\quad \quad \quad\text{  $y\in\mathbb{R}$ or $y\in[-1,1]$},\\
\textrm{TC flow: } &\partial_t\hat{w}_k-\nu(\partial_{r}^2+\f{1}{r}\partial_{r}-\f{k^2}{r^2})\hat{w}_k+(A+\f{B}{r^2})ik\hat{w}_k=0, \quad \text{ $r\in[1,R]$}.
\end{align*}
The linearized operators for the Couette flow and TC flow take the different forms:
\begin{align*}
-\nu(\partial_y^2-k^2)+iky\quad \text{for Couette}  \quad  \text{and}\quad -\nu(\partial_{r}^2+\f{1}{r}\partial_{r}-\f{k^2}{r^2})+ik(A+\f{B}{r^2}) \quad \text{for TC}.
\end{align*}
These two operators exhibit different enhanced dissipation rates of $\nu^{\f13}$ and $\nu^{\f13}|B|^{\f23}R^{-2}$, respectively.  

It can be observed that both the symmetric and antisymmetric parts of the linearized operator for the TC flow have different structures from those of the Couette flow. In the following discussion, we will demonstrate how these differences affect our results and the corresponding proofs for the TC flow.

In previous studies of the planar Couette flow \cite{BMV,BWV,BGM-bams,MZ-2,MZ-1}, where the shear variable $y$ is defined over the entire space $\mathbb{R}$, mathematicians employed the Fourier transform in $y$ to obtain ordinary differential equations with respect to time $t$. However, in the case of TC flow, the radial variable $r$ is defined in the bounded domain $[1, R]$ with boundaries, which makes it inconvenient to conduct the Fourier transform in $r$ directly. Instead, we adopt the method of resolvent estimates to derive enhanced dissipation, with detailed results and proofs provided in Section \ref{Section Resolvent estimate}. Note that this method has also been employed by Chen, Li, Wei and Zhang to study the Couette flow within a finite channel in \cite{CLWZ-2D-C}.

\subsubsection{Resolvent estimates}
In Section \ref{Section Resolvent estimate}, our key results are the resolvent estimates presented in Proposition \ref{resolvent estimate-main}. Define
\begin{equation*}
    L_k:=-\nu(\partial_r^2-\f{k^2-\f14}{r^2})+ikB\f{1}{r^2}
\end{equation*}
and denote
$$F:=(L_k-ikB\lambda)w=-\nu(\partial_r^2-\f{k^2-\f14}{r^2})w+ikB(\f{1}{r^2}-\lambda)w.$$ 
Our initial goal is the following inequality:
\begin{align}
\label{L2 to L2}\nu^{\f13} |kB|^{\f23}\|\f{w}{r}\|_{L^2}\leq C \|rF\|_{L^2},   
\end{align}
which is a resolvent estimate from $L^2$ space to $L^2$ space, as given in Lemma \ref{resolvent estimate-1}. This estimate corresponds to the estimate obtained for the Couette flow in \cite{CLWZ-2D-C}:
\begin{align*}
\nu^{\f13} |k|^{\f23}\|w\|_{L^2}\leq C \|F\|_{L^2}.
\end{align*}
However, due to the specific structure of the TC flow, the corresponding resolvent estimate is weighted in $r$. 
If we take $c'$ as $\f{1}{2C}$,  the inequality \eqref{L2 to L2} directly yields the following estimate:
\begin{align*}
\nu^{\f13} |kB|^{\f23}\|\f{w}{r}\|_{L^2}\leq C \|rF-c'\nu^{\f13} |kB|^{\f23}\f{w}{r}\|_{L^2}, 
\end{align*}
which takes the form of a resolvent estimate for $rF-c'\nu^{\f13} |kB|^{\f23}\f{w}{r}$. Notice that $rF-c'\nu^{\f13} |kB|^{\f23}\f{w}{r}$ can be written as
\begin{align*}
rF-c'\nu^{\f13} |kB|^{\f23}\f{w}{r}=(rL_kr-c'\nu^{\f13} |kB|^{\f23})\f{w}{r}.
\end{align*}
Here, $rL_kr-c'\nu^{\frac{1}{3}} |kB|^{\frac{2}{3}}$ represents the translation of $rL_kr$ to the left by $c'\nu^{\frac{1}{3}} |kB|^{\frac{2}{3}}$, which leads to the enhanced dissipation rate for the following linear evolution equation
\begin{equation*}
   [ \partial_t- (rL_kr-c'\nu^{\f13} |kB|^{\f23}) ]\f{w}{r}=0.
\end{equation*}
This allows us to obtain the dissipative factor $e^{-c\nu^{\f13} |kB|^{\f23}R^{-2}t}$ in space-time estimates for the fully nonlinear system, as shown in detail in Proposition \ref{spacetime estimate proposition-nonlinear}.

Next, we derive the resolvent estimate for the stream function $\varphi$:
\begin{align*}
\nu^{\f16}|kB|^{\f56}|k|^{\f12}\left( \|\varphi'\|_{L^2}+|k|\|\f{\varphi}{r}\|_{L^2} \right) 
\leq CR^2\left((\f{\nu}{|kB|})^{\f16} (\log R)^{\f12} +1 \right) \|rF-c'\nu^{\f13} |kB|^{\f23}\f{w}{r}\|_{L^2}.
\end{align*}
The proof of this estimate relies on the resolvent estimate for $\|w\|_{L^1}$, and its details are given in Lemma \ref{resolvent estimate-2}.

Finally, in Lemma \ref{resolvent estimate-3} and Lemma \ref{resolvent estimate-4} we establish the below estimates for the $H^1$ norms of $w$ and $\varphi$, which are controlled by the $H^{-1}$ norms of the resolvent equations, respectively:
\begin{align*}
\nu\|w\|_{H^1_r}+\nu^{\f23}|kB|^{\f13}\|\f{w}{r}\|_{L^2}\leq& C\|F-c'\nu^{\f13} |kB|^{\f23} R^{-2}w\|_{H^{-1}_r},\\
\nu^{\f{1}{2}}|kB|^{\f12}\|\varphi'\|_{L^2}+\nu^{\f{1}{2}}|k||kB|^{\f12}\|\f{\varphi}{r}\|_{L^2}\leq& C R^2 \|F-c'\nu^{\f13}|kB|^{\f23}R^{-2}w\|_{H^{-1}_r}.
\end{align*}
The resolvent estimates presented above are crucial for our analysis, and are summarized in Proposition \ref{resolvent estimate-main}. Notably, all of the estimates in Proposition \ref{resolvent estimate-main} have been carefully derived in preparation for the subsequent discussion in Section \ref{Enhanced dissipation and invisid damping}.

Based on the above conclusions, we summarize the following table which illustrates the relation between coefficients of the terms on the left of  resolvent estimates and the norm of $F-c'\nu^{\f13}|kB|^{\f23}R^{-2}w$ on the right.

\begin{table}[htbp]
	\centering
	\caption{The coefficients on the left of resolvent estimates}
	\label{tab:1}  
	\begin{tabular}{cccc ccc}
		\hline\hline\noalign{\smallskip}	
	$\|w'\|_{L^2}$ & $\|\f{w}{r}\|_{L^2}$ & $\|\varphi'\|_{L^2}$& $|k|\|\f{\varphi}{r}\|_{L^2}$ &Norm of $F-c'\nu^{\f13}|kB|^{\f23}R^{-2}w$\\
		\noalign{\smallskip}\hline\noalign{\smallskip}
$\nu^{\f23}|kB|^{\f13}$ &$\nu^{\f13}|kB|^{\f23}$  &$\nu^{\f16}|kB|^{\f56}$ &$\nu^{\f16}|kB|^{\f56}$ & $L^2$   \\
$\nu$ & $\nu^{\f23}|kB|^{\f13}$ & $\nu^{\f12}|kB|^{\f12}$ & $\nu^{\f12}|kB|^{\f12}$ & $H^{-1}_r$ \\
		\noalign{\smallskip}\hline
	\end{tabular}
\end{table}
The above powers of $\nu$ in the coefficients are optimal and consistent with the 2D Couette flow in \cite{CLWZ-2D-C}. However, in this paper the coefficients are also dependent on $B$ and $R$, and the powers of $B$ are also sharp. As shown in the table, the coefficients on the left will be multiplied by $(\frac{\nu}{|kB|})^{\frac{1}{3}}$ if the regularity of $F-c'\nu^{\frac{1}{3}}|kB|^{\frac{2}{3}}R^{-2}w$ drops from $L^2$ to $H^{-1}_{r}$.

\subsubsection{Pseudospectrum and enhanced dissipation} In the section \ref{Enhanced dissipation and invisid damping}, the pseudospectral bound of $L$ is defined as
\begin{align*}
\Psi(L)=\inf\{\|(L-i\lambda)f\|:f\in D(L),\lambda\in\mathbb{R},\|f\|=1\}.
\end{align*}
Recall that in Proposition \ref{resolvent estimate-main} we establish the following resolvent estimates
\begin{align*}
&\|r(L_k-i\lambda)w\|_{L^2}\geq C(\nu k^2)^{\f13}|B|^{\f23}\|\f{w}{r}\|_{L^2}.
\end{align*}
Since $r\in [1, R]$, we have
\begin{align*}
&\|(L_k-i\lambda)w\|_{L^2}\geq C(\nu k^2)^{\f13}|B|^{\f23}R^{-2}\|w\|_{L^2}.
\end{align*} 
This provides a lower bound for the pseudospectrum of $L_k$:
\begin{align*}
\Psi(L_k(L^2\rightarrow L^2))\geq C(\nu k^2)^{\f13}|B|^{\f23}R^{-2}.
\end{align*}
Applying the Gearhart-Pr$\ddot{u}$ss type lemma established by Wei in \cite{Wei} and by Helffer-Sj\"{o}strand in \cite{HS}, we can obtain the pointwise estimate of the semigroup bound in Proposition \ref{linear exp decay}
\begin{align*}
	\|w_k^l(t)\|_{L^2}\leq Ce^{-c(\nu k^2)^{\f13}|B|^{\f23}R^{-2}t}\|w_k(0)\|_{L^2}, \quad \textrm{for any} \ t\geq0,
	\end{align*}
where $w_k^l(t)$ satisfies the homogeneous linear equation
\begin{align*}
\partial_tw_k^l+L_kw_k^l=0 \quad \textrm{with} \ w_k^l(0)=w_k(0).
\end{align*}
It is worth noting that here, the exponential decay factor $e^{-c(\nu k^2)^{\f13}|B|^{\f23}R^{-2}t}$ reflects the enhanced dissipation effect.

\subsubsection{Integrated invisid damping} The second part of section \ref{Enhanced dissipation and invisid damping} presents estimations for the stream function $\varphi^l_k$ governed by the initial vorticity. Here  $\varphi_k^l$ satisfies the linear elliptic equation
\begin{equation}\label{def linear phi k l} 
	    (\partial_r^2-\f{k^2-\f14}{r^2}){\varphi}_k^l={w}_k^l,\quad {\varphi}_k^l|_{r=1,R}=0,\quad r\in[1,R].
	\end{equation} 
 The precise estimate is given in Lemma \ref{nonzero frequency linear part for varphi}, which states that, for any $k\in \mathbb{Z}$ and $|k|\geq1$, it holds 
\begin{align*}
     &k^2|B|R^{-4}(\|e^{c'(\nu k^2)^{\f13}|B|^{\f23}R^{-2}t}\partial_r{\varphi}^{l}_k\|_{L^2L^2}^2+|k|^2\|e^{c'(\nu k^2)^{\f13}|B|^{\f23}R^{-2}t}\f{{\varphi}^{l}_k}{r}\|_{L^2L^2}^2) \\
     \lesssim& (\log R)^{-2}R^{-4}\|r^2 w_k(0)\|_{L^2}^2+R^6 \|\f{w_k(0)}{r^3}\|_{L^2}^2+(\f{\nu}{|kB|})^{\f23}R^2\|\partial_r w_k(0)\|_{L^2}^2 \\&+\left((\f{\nu}{|kB|})^{\f13} \log R +1 \right)\left( R^{-2}\|rw_k(0)\|_{L^2}^2 + (\f{\nu}{|kB|})^{\f43}k^4 R^2\|\f{w_k(0)}{r}\|_{L^2}^2\right).
\end{align*}
Note that $\nu$ is much smaller than $1$, and $B$ is a fixed constant, the leading terms on the right of above inequality are 
\begin{equation*}
    (\log R)^{-2}R^{-4}\|r^2 w_k(0)\|_{L^2}^2+R^6 \|\f{w_k(0)}{r^3}\|_{L^2}^2+R^{-2}\|rw_k(0)\|_{L^2}^2.
\end{equation*}
Thus, both the coefficients in front of $\varphi_k^l$ terms on the left and the main initial-vorticity terms on the right are independent of the viscosity coefficient $\nu$. This phenomenon is commonly refer to as the effect of inviscid damping. Therefore, Lemma \ref{nonzero frequency linear part for varphi} presents a space-time version of linear inviscid damping around 2D TC flow. 

To prove this result, in Lemma \ref{nonzero frequency linear part for varphi l1 part} we find and formulate a new set of explicit orthonormal basis of (weighted) eigenfunctions corresponding to the Laplace operator for circular flows, namely $-\partial_r^2+(k^2-\f14)/r^2$. This basis can be represented in the closed form
\begin{equation*}
    \psi_l(r)=(\f{2}{\log R})^{\f12}r^{1/2}\sin (\f{l\pi}{\log R} \log r), \quad \textrm{for} \ l\in \mathbb{N}_+
\end{equation*}
and satisfies the equation
\begin{equation*}
    (\partial_r^2-\frac{k^2-\frac{1}{4}}{r^2})\psi_l=-\frac{(\f{l\pi}{\log R} )^2+k^2}{r^2}\psi_l, \quad \psi_l|_{r=1, R}=0.
\end{equation*}
Note that there is the $1/r^2$ on the right and $-(\f{l\pi}{\log R} )^2-k^2$ is not the canonical eigenvalue. As far as we know, our explicit constructions here are new. With this basis, we move to prove a crucial estimate
\begin{equation}\label{|kB|(partial r phil0+phil0) est}
|kB||k|(\log R)^2\left( \|\partial_r \tilde{\varphi}^{l_0}_k\|_{L^2_t L^2_r}^2+k^2\|\f{\tilde{\varphi}^{l_0}_k}{r}\|^2_{L^2_t L^2_r}\right)\lesssim \|r^{2}w_k(0)\|_{L^2}^2,
\end{equation}
where $k\in \mathbb{Z}$, $B$ and $R$ are constants with $R>1$, $w_k(0) \in L^2([1, R])$, and $\tilde{\varphi}^{l_0}_k$ is defined as the solution to the elliptic equations below:
\begin{equation*}
(\partial_r^2-\f{k^2-\f14}{r^2})\tilde{\varphi}^{l_0}_k=e^{-ikB\f{t}{r^2}} w_k(0),\quad \tilde{\varphi}^{l_0}_k|_{r=1,R}=0,\quad r\in[1,R],\quad t\geq0.
\end{equation*}
Notice that the Laplace operator for the 2D Couette flow is given by $-\partial_y^2+k^2$ within the domain $y\in[-1,1]$. To find corresponding orthonormal eigenfunctions, it is natural to introduce the set $\{\sin\left(\pi j(y+1)/2\right) \}_{j=1}^{\infty}$, as described in \cite{CLWZ-2D-C}. However, for the circular Laplace operator $\partial_r^2-\frac{k^2-\frac{1}{4}}{r^2}$, we need to construct a set of more complicated eigenfunctions. The main issue is that if we directly consider the following eigenvalue problem:
\begin{equation*}
(\partial_r^2-\frac{k^2-\frac{1}{4}}{r^2})\psi=\lambda \psi, \quad \psi|_{r=1, R}=0,
\end{equation*}
it can be shown that there is no closed-form solution. In this paper, we instead solve the weighted eigenvalue problem:
\begin{equation*}
(\partial_r^2-\frac{k^2-\frac{1}{4}}{r^2})\psi=\lambda w\psi, \quad \psi|_{r=1, R}=0
\end{equation*}
with $w(r)=\frac{1}{r^2}$. Utilizing change of variables, we obtain a family of closed-form solutions to this boundary value problem, which can be written in the following explicit expressions
\begin{equation*}
\psi_l(r)=(\f{2}{\log R})^{\f12}r^{1/2}\sin (\f{l\pi}{\log R} \log r), \quad \textrm{for} \ l\in \mathbb{N}_+.
\end{equation*}
These functions form an orthonormal basis of the weighted space $L^2_w([1, R])$. Therefore, with the aid of $\{ \psi_l \}_{l=1}^{\infty}$, we can evaluate 
\begin{equation*}
   \begin{split}
      & -\langle \tilde{w}^{l_0}_k, \tilde{\varphi}^{l_0}_k \rangle=
    -\langle w^{-1}\tilde{w}^{l_0}_k, \tilde{\varphi}^{l_0}_k \rangle_w=\sum\limits_{l=1}^{\infty} \f{1}{(\f{l\pi}{\log R} )^2+k^2} \langle w^{-1}\tilde{w}^{l_0}_k, \psi_l \rangle_w \overline{\langle \tilde{w}^{l_0}_k , \psi_l \rangle} 
    =\sum\limits_{l=1}^{\infty}  \frac{|\langle \tilde{w}^{l_0}_k , \psi_l \rangle|^2}{(\f{l\pi}{\log R} )^2+k^2},
    \end{split}
\end{equation*}
where $\tilde{w}^{l_0}_k:=e^{-i kB\f{t}{r^2}}w_k(0)$. Notice that
\begin{equation*}
\begin{split}
    \langle \tilde{w}^{l_0}_k, \psi_l \rangle=&\int_{1}^{R} e^{-ikB\f{t}{r^2}}w_k(0, r)\psi_l(r) dr\\=&\int_{\f{1}{R^2}}^1 e^{-i kBt s} w_k(0, \f{1}{\sqrt{s}}) \psi_l (\f{1}{\sqrt{s}})\f{ds}{2s\sqrt{s}}.
    \end{split}
\end{equation*}
Thus, employing Fourier transform and Plancherel's formula, the desired estimate \eqref{|kB|(partial r phil0+phil0) est} follows.

\begin{comment}
They can be written in the following explicit expressions
 \begin{align*}
        \cup_{l\in\mathbb{Z}}\{\f{2}{\log R}r^{1/2}\sin (\f{l\pi}{\log R} \log r)\} .
     \end{align*}
\end{comment}
 
  It is worth noting that Lemma \ref{nonzero frequency linear part for varphi} also plays an important role in improving the nonlinear transition threshold of Theorem \ref{main-transition threshold}. 

\subsubsection{Space-time estimates for the linearized Navier-Stokes equations}
In Proposition \ref{nonzero nonlinear-original} of Section \ref{3-space time}, we establish the space-time estimates for the linearized 2D Navier-Stokes equation written in the vorticity formulation \eqref{scaling nonlinear}:
\begin{align*}
\partial_t w_k&- \nu(\partial_r^2-\f{k^2-\f14}{r^2})w_k+\f{ikB}{r^2}w_k+\f{1}{r}[ik\sum_{l\in\mathbb{Z}}\partial_r(r^{-\f12}\varphi_l)w_{k-l}-r^{\f12}\partial_r(\sum_{l\in\mathbb{Z}}ilr^{-1}\varphi_lw_{k-l})]=0
\end{align*}
with the stream function $\varphi_k$  satisfying
\begin{align*}
&(\partial_r^2-\f{k^2-\f14}{r^2})\varphi_k=w_k,\quad \varphi_k|_{r=1,R}=0,\quad r\in[1,R],\quad t\geq0.
\end{align*}
To achieve so, we employ the estimates in Section \ref{Section Resolvent estimate} and Section \ref{Enhanced dissipation and invisid damping}.

\noindent (1) The space-time estimate for the zero frequency part (i.e. $k=0$) can directly obtained in view of the equation structure and integration by parts, as detailed in Lemma \ref{zero frequency part}.

\noindent (2) For the nonzero frequency part with $k\in \mathbb{Z}\backslash \{0 \}$, we decompose  the solution $w_k=w_k^l+w_k^{n}$, where $w_k^l$ obeys the homogeneous linear equation
\begin{align*}
&\partial_tw_k^l+L_kw_k^l=0,\quad w_k^l(0)=w_k(0),
\end{align*}
while $w_k^{n}$ is the solution to the inhomogeneous linear equation with zero initial condition
\begin{align*}
&\partial_tw_k^{n}+L_kw_k^{n}+\f{1}{r}[ikf_1-r^{\f12}\partial_{r}(r^{\f12}f_2)]=0,\quad w_k^{n}(0)=0.
\end{align*}
Here we denote the nonlinear forms by $f_1$ and $f_2$, which are expressed as
    \begin{align*}
f_1=\sum_{l\in\mathbb{Z}}\partial_r(r^{-\f12}\varphi_l)w_{k-l},\quad f_2=\sum_{l\in\mathbb{Z}}ilr^{-\f32}\varphi_lw_{k-l}.
\end{align*}
As a result, we also decompose $\varphi_k$ as $\varphi_k=\varphi^{l}_k+\varphi^{n}_k$, where $\varphi_k^l$ and $\varphi_k^n$ fulfill
\begin{equation*}
    (\partial_r^2-\f{k^2-\f14}{r^2})\varphi^{l}_k=w_k^l,\quad \varphi^{l}_k|_{r=1,R}=0,\quad \textrm{with} \ r\in[1,R] \ \textrm{and} \ t\geq0,
\end{equation*}
and 
\begin{equation*}
    (\partial_r^2-\f{k^2-\f14}{r^2})\varphi^{n}_k=w_k^n,\quad \varphi^{n}_k|_{r=1,R}=0,\quad \textrm{with} \ r\in[1,R] \ \textrm{and} \ t\geq0.
\end{equation*}

For our space-time estimate, the first step is to provide a control of $w_k^l(t)$ based on the initial data $w_k^l(0)$. With the aid of the pseudospectrum bounds and enhanced dissipation effect proved in Section \ref{Section Resolvent estimate}, we establish the following inequality in Lemma \ref{nonzero frequency linear part}, which holds for any $k\in \mathbb{Z}$ and $|k|\geq1$:
	\begin{align*}
&\|e^{c'(\nu k^2)^{\f13}|B|^{\f23}R^{-2}t}w_k^l\|_{L^{\infty}L^2}^2+(\nu k^2)^{\f13}|B|^{\f23}R^{-2}\|e^{c'(\nu k^2)^{\f13}|B|^{\f23}R^{-2}t}w_k^l\|_{L^2L^2}^2\\
 &+\nu\|e^{c'(\nu k^2)^{\f13}|B|^{\f23}R^{-2}t}\partial_{r}w_k^l\|_{L^2L^2}^2 +\nu k^2\|e^{c'(\nu k^2)^{\f13}|B|^{\f23}R^{-2}t}\f{w_k^l}{r}\|_{L^2L^2}^2
\lesssim\|w_k(0)\|_{L^2}^2.
	\end{align*}

To proceed, we derive an integrated linear invisid damping estimate, in which $\varphi_k^l$ can be controlled with the initial data of the vorticity
\begin{equation*}
\begin{split}
     &k^2|B|R^{-4}(\|e^{c'(\nu k^2)^{\f13}|B|^{\f23}R^{-2}t}\partial_r{\varphi}^{l}_k\|_{L^2L^2}^2+|k|^2\|e^{c'(\nu k^2)^{\f13}|B|^{\f23}R^{-2}t}\f{{\varphi}^{l}_k}{r}\|_{L^2L^2}^2) \\
     \lesssim& (\log R)^{-2}R^{-4}\|r^2 w_k(0)\|_{L^2}^2+R^6 \|\f{w_k(0)}{r^3}\|_{L^2}^2+(\f{\nu}{|kB|})^{\f23}R^2\|\partial_r w_k(0)\|_{L^2}^2 \\&+\left((\f{\nu}{|kB|})^{\f13} \log R +1 \right)\left( R^{-2}\|rw_k(0)\|_{L^2}^2 + (\f{\nu}{|kB|})^{\f43}k^4 R^2\|\f{w_k(0)}{r}\|_{L^2}^2\right),
\end{split}
\end{equation*}
as detailed in Lemma \ref{nonzero frequency linear part for varphi}. Note that $\nu$ is much smaller than $1$, and $B$ is a fixed constant, so the initial vorticity has an upper bound independent of the viscosity coefficient $\nu$. 

By combining this observation with the estimates for $w_k^l$, we obtain the following inequality
\begin{align*}
&\|e^{c'(\nu k^2)^{\f13}|B|^{\f23}R^{-2}t}w_k^l\|_{L^{\infty}L^2}+|k||B|^{\f12}\|e^{c'(\nu k^2)^{\f13}|B|^{\f23}R^{-2}t}\partial_r{\varphi}^{l}_k\|_{L^2L^2}\\
\lesssim& C(R) (\|\partial_r w_k(0)\|_{L^2}+\|w_k(0)\|_{L^2}).
\end{align*}
It is worth mentioning that the coefficients in front of the terms on both sides of the inequality are independent of viscosity. Therefore, similar estimates can also hold for the inviscid flow (i.e. $\nu=0$).

By applying the resolvent estimates from Proposition \ref{resolvent estimate-main}, we then arrive at the conclusion stated in Proposition \ref{spacetime estimate proposition-nonlinear} for the nonlinear part, which involves the inhomogeneous linear equation with zero initial data. This result allows us to bound $w_k^n$ and $\varphi^{n}_k$ by the nonlinear terms $f_1$ and $f_2$ as follows
\begin{align*}
&\|e^{c'(\nu k^2)^{\f13}|B|^{\f23}R^{-2}t}w_k^n\|_{L^{\infty}L^2}+(\nu k^2)^{\f16}|B|^{\f13}R^{-1}\|e^{c'(\nu k^2)^{\f13}|B|^{\f23}R^{-2}t}w_k^n\|_{L^2L^2}\\
&+\nu^{\f12}\|e^{c'(\nu k^2)^{\f13}|B|^{\f23}R^{-2}t}\partial_rw_k^n\|_{L^2L^2}+(\nu k^2)^{\f12}\|e^{c'(\nu k^2)^{\f13}|B|^{\f23}R^{-2}t}\f{w_k^n}{r}\|_{L^2L^2}\\
&+|B|^{\f12}R^{-2}\left(|k|\|e^{c'(\nu k^2)^{\f13}|B|^{\f23}R^{-2}t}\partial_r\varphi_k^n\|_{L^2L^2}+k^2\|e^{c'(\nu k^2)^{\f13}|B|^{\f23}R^{-2}t}\f{\varphi_k^n}{r}\|_{L^2L^2}\right)\\
\lesssim&\left(|kB|^{-\f12} (\log R)^{\f12} +\nu^{-\f16} |kB|^{-\f13} \right) \cdot\|e^{c'(\nu k^2)^{\f13}|B|^{\f23}R^{-2}t}kf_1\|_{L^2L^2}\\
&+\nu^{-\f12} \|e^{c'(\nu k^2)^{\f13}|B|^{\f23}R^{-2}t}f_2\|_{L^2L^2}.
	\end{align*}

Finally, by combining the aforementioned  three estimates, we derive the crucial proposition of Section \ref{3-space time}, namely Proposition \ref{nonzero nonlinear-original}. This proposition states that the following inequality holds
	\begin{align*}
&\|e^{c'(\nu k^2)^{\f13}|B|^{\f23}R^{-2}t}w_k\|_{L^{\infty}L^2}+(\nu k^2)^{\f16}|B|^{\f13}R^{-1}\|e^{c'(\nu k^2)^{\f13}|B|^{\f23}R^{-2}t}w_k\|_{L^2L^2}\\
&+\nu^{\f12}\|e^{c'(\nu k^2)^{\f13}|B|^{\f23}R^{-2}t}\partial_r w_k\|_{L^2L^2}+(\nu k^2)^{\f12}\|e^{c'(\nu k^2)^{\f13}|B|^{\f23}R^{-2}t}\f{w_k}{r}\|_{L^2L^2}\\
&+|B|^{\f12}R^{-2}\left(|k|\|e^{c'(\nu k^2)^{\f13}|B|^{\f23}R^{-2}t}\partial_r\varphi_k\|_{L^2L^2}+k^2\|e^{c'(\nu k^2)^{\f13}|B|^{\f23}R^{-2}t}\f{\varphi_k}{r}\|_{L^2L^2}\right)\\
\lesssim&\|w_k(0)\|_{L^2}+R^{-2}(\log R)^{-2}\|r^2w_k(0)\|_{L^2}+R^{3}\|\f{w_k(0)}{r^3}\|_{L^2}+(\f{\nu}{|kB|})^{\f13}R\|\partial_rw_k(0)\|_{L^2}\\
&+\left(1+(\f{\nu}{|kB|})^{\f13}\log R\right)^{\f12} \cdot \left(R^{-1}\|rw_k(0)\|_{L^2}+R\|\f{w_k(0)}{r}\|_{L^2}(\f{\nu}{|kB|})^{\f23}k^2\right)\\
&+\left(|kB|^{-\f12} (\log R)^{\f12} +\nu^{-\f16} |kB|^{-\f13} \right) \cdot\|e^{c'(\nu k^2)^{\f13}|B|^{\f23}R^{-2}t}kf_1\|_{L^2L^2}\\
&+\nu^{-\f12} \|e^{c'(\nu k^2)^{\f13}|B|^{\f23}R^{-2}t}f_2\|_{L^2L^2}\\=&:\mathcal{M}_k(0).
	\end{align*}
We refer to it as the space-time estimates for the vorticity $w_k$ and the stream function $\varphi_k$.

\subsubsection{Nonlinear stability} Finally, we close the energy estimate with a bootstrap argument. By utilizing the spacetime estimates established in Proposition \ref{nonzero nonlinear-original} from Section \ref{3-space time}, we derive the result of Proposition \ref{nonzero nonlinear-all cases}:
	\begin{align*}
&\|e^{c'\mu_kt}w_k\|_{L^{\infty}L^2}+\mu_k^{\f12}\|e^{c'\mu_kt}w_k\|_{L^2L^2}+\nu^{\f12}\|e^{c'\mu_kt}w_k'\|_{L^2L^2}+(\nu k^2)^{\f12}\|e^{c'\mu_kt}\f{w_k}{r}\|_{L^2L^2}\\
&+\max\{|B|^{\f12}R^{-2},  (\nu k^2)^{\f12} R^{-2} \}\cdot\left(|k|\|e^{c'\mu_kt}\partial_r\varphi_k\|_{L^2L^2}+k^2\|e^{c'\mu_kt}\f{\varphi_k}{r}\|_{L^2L^2}\right)\\
\lesssim&\mathcal{M}_k(0)+\mu_k^{-\f12}R^{-1}\|e^{c'\mu_kt}kf_1\|_{L^2L^2}+\nu^{-\f12} \|e^{c'\mu_kt}f_2\|_{L^2L^2}.
\end{align*}

Here, we define $\mu_k:=\max\left\{(\nu k^2)^{\f13} |B|^{\f23} R^{-2}, \nu k^2 R^{-2} \right\}$, as it is necessary to use distinct energy estimates for the high-frequency case of $\nu k^2\geq |B|$ and the low-frequency case of $\nu k^2\leq |B|$ separately. In the case of 2D Couette flow \cite{CLWZ-2D-C}, the discussions on different frequencies $k$ correspond to $\nu k^2\geq1$ and $\nu k^2\leq1$. It is important to mention that the frequency classification of $k$ for the 2D TC flow also relies on the rotation parameter $B$.

We then define the energy functional $E_k$ as follows:
\begin{align*}
E_0=&\|w_0\|_{L^{\infty}L^2},\\
E_k=
&\|e^{c'\mu_kt}w_k\|_{L^{\infty}L^2}+\mu_k^{\f12}\|e^{c'\mu_kt}w_k\|_{L^2L^2}+|B|^{\f12}|k|^{\f32}R^{-2}\|e^{c'\mu_kt}\f{\varphi_k}{r^{\f12}}\|_{L^2L^{\infty}} \quad  \textrm{for} \ |k|\ge 1.  
\end{align*}
According to Proposition \ref{nonzero nonlinear-all cases}, it  can be infered that for $k\in \mathbb{Z}\backslash \{0 \}$, it holds
	\begin{align*}
&E_k\lesssim\mathcal{M}_k(0)+\mu_k^{-\f12}R^{-1}\|e^{c'\mu_kt}kf_1\|_{L^2L^2}+\nu^{-\f12} \|e^{c'\mu_kt}f_2\|_{L^2L^2}.
\end{align*}
Utilizing Lemma \ref{f1 f2}, we can control nonlinear forms $f_1, f_2$ by
\begin{equation*}
    \mu_k^{-\f12}R^{-1}\|e^{c'(\nu k^2)^{\f13}|B|^{\f23}R^{-2}t}kf_1\|_{L^2L^2}\lesssim \nu^{-\f12} |B|^{-\f12}R^2  (\f{R}{R-1})^{\f12}(1+\log R)\sum_{l\in\mathbb{Z}}E_lE_{k-l},
\end{equation*}
and
\begin{equation*}
\begin{split}
\nu^{-\f12} \|e^{c'\mu_k t}f_2\|_{L^2L^2}
\leq\nu^{-\f12} |B|^{-\f12}R^2\sum_{l\in\mathbb{Z}\backslash\{0, k\}}E_lE_{k-l}.
\end{split}
\end{equation*}
This allows us to obtain Lemma \ref{E_k k non-zero}, which states that
 \begin{equation*}
       E_k\lesssim
\mathcal{M}_k(0)
+\nu^{-\f12}|B|^{-\f12}R^2 (\f{R}{R-1})^{\f12}(1+\log R)\sum_{l\in\mathbb{Z}}E_lE_{k-l}.
    \end{equation*}
Simultaneously, we can derive analogous  estimates for $E_0$ by introducing $\mathcal{M}_0(0)=\|w_0(0)\|_{L^2}$. Incorporating  all estimates above through a bootstrap argument, we finally arrive at the main theorem of this paper: Theorem \ref{transition threshold}.

\subsubsection{Subcritical transition and transition threshold } 
 
Theorem \ref{main theorem-2} and Theorem \ref{main-transition threshold} follow from Proposition \ref{linear exp decay} and Theorem \ref{transition threshold}, respectively. Here are some further remarks regarding these two main theorems.

On one hand, Theorem \ref{main theorem-2} shows that the linearized Navier-Stokes equations around the 2D Taylor-Couette flow is dynamically stable at any Reynolds number (including the inviscid case), and the vorticity at any time can be controlled by the given initial value in $L^2$ norm. On the other hand, for the full Navier-Stokes equation, we only expect nonlinear asymptotic stability within the range of small perturbation of size $\nu^{\beta(s)}(\beta(s)>0)$ in $H^s$ space. 

This proved in our Theorem \ref{main-transition threshold}, provided the initial perturbation does not exceed $\nu^\f12 |B|^{\f12}R^{-2}$. The exponent $\f12$ of $\nu$ agrees with that of the 2D Couette flow \cite{BWV,CLWZ-2D-C}, and provides an upper bound $\f12$ for the corresponding transition threshold $\beta$. The dependency of the subcritical transition on the rotational speed $B$ is captured by $|B|^{\f12}$. This indicates that the system can achieve global stability under a wider range of initial data via taking account of the rotating effect.

\section{\textbf{Derivation of the perturbative equation}}\label{Derivation of the perturbative equation}

We take the Fourier transform in the $\theta$ direction and denote the Fourier coefficients of $w$ and $\varphi$ by $\hat{w}_k$ and $\hat{\varphi}_k$, respectively. Using this notation, the equation \eqref{pertubation of NS vor polar coordinates} can be rewritten as follows:
\begin{align}\label{pertubation of NS vor-fourier}\partial_t\hat{w}_k&-\nu[(\partial_{r}^2+\f{1}{r}\partial_{r}-\f{k^2}{r^2})]\hat{w}_k+(A+\f{B}{r^2})ik\hat{w}_k\\
\nonumber&+\f{1}{r}[\sum_{l\in\mathbb{Z}}i(k-l)\partial_r\hat{\varphi}_l\hat{w}_{k-l}-\sum_{l\in\mathbb{Z}}il\hat{\varphi}_{l}\partial_{r}\hat{w}_{k-l}]=0
\end{align}
with $\hat{\varphi}_k$ satisfying $(\partial_r^2+\f{1}{r}\partial_r-\f{k^2}{r^2})\hat{\varphi}_k=\hat{w}_k$. 

To eliminate the first derivative $\f{1}{r}\partial_r$, we introduce the weight $r^{\f12}$ and define $w_k$ and $\varphi_k$ as $w_k:=r^{\f12}e^{ikAt}\hat{w}_k$ and $\varphi_k:=r^{\f12} e^{ikAt}\hat{\varphi}_k$. As a consequence, the equation \eqref{pertubation of NS vor-fourier} is transferred into  
\begin{align}
\label{scaling nonlinear}\partial_t w_k&- \nu(\partial_r^2-\f{k^2-\f14}{r^2})w_k+\f{ikB}{r^2}w_k\\
\nonumber&+\f{1}{r}[ik\sum_{l\in\mathbb{Z}}\partial_r(r^{-\f12}\varphi_l)w_{k-l}-r^{\f12}\partial_r(\sum_{l\in\mathbb{Z}}ilr^{-1}\varphi_lw_{k-l})]=0,
\end{align}
where $\varphi_k$  satisfies
\begin{align}
&(\partial_r^2-\f{k^2-\f14}{r^2})\varphi_k=w_k,\quad \varphi_k|_{r=1,R}=0\quad \textrm{with}\ r\in[1,R] \ \textrm{and} \ t\geq0.
\end{align}
Note that the Navier boundary condition \begin{equation*}
    w_k|_{r=1,R}=0 , \  \varphi_k|_{r=1,R}=0 \quad \textrm{for any}\ k\in \mathbb{Z}
 \end{equation*}
 implies that $\partial_t w_k$ and the nonlinear terms in \eqref{scaling nonlinear} vanish on the boundary of interval $[1, R]$. This forces the second derivative of $w_k$ to be zero on the boundary, i.e.
 \begin{equation*}
     \partial_r^2 w_k|_{r=1,R}=0  \quad \textrm{for any}\ k\in \mathbb{Z}.
 \end{equation*}   
 Up to this point, we convert our problem to exploring the dynamics and long-time behaviors of \eqref{scaling nonlinear}. 

\begin{remark}
In below, for two quantities $A$ and $B$, we frequently use $A\lesssim B$ in short to stand for the inequality $A\le C B$ with some universal constant $C>0$,  that is independent of $\nu, k, B, \lambda$ and $R$. Additionally, we also write $A\approx B$ to indicate that both $A\lesssim B$ and $B\lesssim A$ are true.
\end{remark}

\section{\textbf{Resolvent estimates }}\label{Section Resolvent estimate}

 To establish decays of the linearized equation of \eqref{scaling nonlinear} as below
\begin{align*}
\partial_t w_k&- \nu(\partial_r^2-\f{k^2-\f14}{r^2})w_k+\f{ikB}{r^2}w_k=0,
\end{align*}
the key step is to derive the resolvent estimates that will be presented in Proposition \ref{resolvent estimate-main}. More precisely, we want to study the resolvent equation subject to the following (Navier) boundary condition for any $\lambda\in \mathbb{R}$. After  taking the Fourier transform with respect to time $t$, the resolvent equation becomes
\begin{align}\label{vorticity eqn}
&-\nu(\partial_r^2-\f{k^2-\f14}{r^2})w+ikB(\f{1}{r^2}-\lambda)w=F \quad \textrm{with} \ w|_{r=1,R}=0 \ \textrm{and} \ \partial_r^2 w_k|_{r=1,R}=0.
\end{align}

The domain of the operator is defined as
\begin{align*}
D_k=\{w\in H_{loc}^2(\mathbb{R}_{+},dr)\cap L^2(\mathbb{R}_{+},dr):-\nu(\partial_r^2-\f{k^2-\f14}{r^2})w+i\f{kB}{r^2}w\in L^2(\mathbb{R}_{+},dr) \}.
\end{align*}
Note that for any $|k|\geq1$, it holds
\begin{align*}
D_k&=\{w\in L^2(\mathbb{R}_{+},dr):\partial_r^2w,\f{w}{r^2}\in L^2(\mathbb{R}_{+},dr) \}.
\end{align*}
We also introduce 
$$\|f\|_{H^1_r}^2:=\|f'\|_{L^2}^2+\|\f{f}{r}\|_{L^2}^2 \quad \textrm{and}  \quad \|f\|_{H^{-1}_r}:=\sup\limits_{\|g\|_{H^1_r}\le1} |\langle f,g\rangle|.$$
Here $\langle \ , \ \rangle$ represents the canonical inner product in $L^2(\mathbb{R}_{+},dr)$.

In this section, we will present the resolvent  estimates for $w'$, $\f{w}{r}$, $\varphi'$ and $\f{\varphi}{r}$ with respect to $F$ in both $L^2$ norm and $H_r^{-1}$ norm.

\subsection{Coercive estimates of the real part}
We start with the coercive estimates for the real part of $ \langle F, w \rangle$, which will be used to obtain the desired resolvent estimates for $w$.

\begin{lemma}\label{trivial w' lemma}For any $|k|\geq1$ and $w\in D_k$, it holds
  \begin{align}
\label{trivial w'}&\Re\langle F,w\rangle=\nu\|w'\|_{L^2}^2+\nu(k^2-\f14)\|\f{w}{r}\|_{L^2}^2.
\end{align}
\end{lemma}
\begin{proof}Via integration by parts, one can directly check
\begin{align*}
&\Re\langle F,w\rangle=\Re\langle -\nu(\partial_r^2-\f{k^2-\f14}{r^2})w,w\rangle=\nu\|w'\|_{L^2}^2+\nu(k^2-\f14)\|\f{w}{r}\|_{L^2}^2.
\end{align*}
\end{proof}

\subsection{Resolvent estimates}
For equation \eqref{vorticity eqn}, this subsection aims to establish upper bounds for $w$ and $\varphi$ in terms of $F$. We prove the following inequalities.
\begin{proposition}\label{resolvent estimate-main}For any $|k|\geq1$, $\lambda\in\mathbb{R}$ and $w\in D_k$, there exist constants $C,c>0$ independent of $\nu,k,B,\lambda,R$, such that for any $0\leq c'\leq c$, we have
\begin{align*}
  &\nu^{\f23} |kB|^{\f13}\|w'\|_{L^2}+\nu^{\f13} |kB|^{\f23}\|\f{w}{r}\|_{L^2} \leq C\|rF-c'\nu^{\f13} |kB|^{\f23}\f{w}{r}\|_{L^2},\\
&\nu^{\f16}|kB|^{\f56}|k|^{\f12}\left( \|\varphi'\|_{L^2}+|k|\|\f{\varphi}{r}\|_{L^2} \right) 
\leq CR^2\left((\f{\nu}{|kB|})^{\f16} (\log R)^{\f12} +1 \right) \|rF-c'\nu^{\f13} |kB|^{\f23}\f{w}{r}\|_{L^2}.
\end{align*}
Moreover, it also holds
\begin{align*}
\nu\|w\|_{H^1_r}+\nu^{\f23}|kB|^{\f13}\|\f{w}{r}\|_{L^2}\leq& C\|F-c'\nu^{\f13} |kB|^{\f23} R^{-2}w\|_{H^{-1}_r},\\
\nu^{\f{1}{2}}|kB|^{\f12}\|\varphi'\|_{L^2}+\nu^{\f{1}{2}}|k||kB|^{\f12}\|\f{\varphi}{r}\|_{L^2}\leq& C R^2 \|F-c'\nu^{\f13}|kB|^{\f23}R^{-2}w\|_{H^{-1}_r}.
\end{align*}
\end{proposition}
\begin{comment}
\begin{proof}
    This is a direct conclusion of Lemma \ref{resolvent estimate-1}, Lemma \ref{resolvent estimate-2}, Lemma \ref{resolvent estimate-3} and Lemma \ref{resolvent estimate-4}.
\end{proof}
\end{comment}
The proof of this proposition can be separated into four parts. We first derive the resolvent estimates for $w$ from $L^2$ to $L^2$.
\begin{lemma}\label{resolvent estimate-1}For any $|k|\geq1$, $\lambda\in\mathbb{R}$ and $w\in D_k$, there exist constants $C>0$ independent of $\nu,k,B,\lambda,R$, such that the following estimate holds
\begin{align}\label{resolvent estimate L2}
\nu^{\f23} |kB|^{\f13}\|w'\|_{L^2}+\nu^{\f13} |kB|^{\f23}\|\f{w}{r}\|_{L^2}+|kB|\|r(\f{1}{r^2}-\lambda)w\|_{L^2}\leq C \|rF\|_{L^2}.
\end{align}
Moreover, there exist  constants $C,c>0$ independent of $\nu,k,B,\lambda,R$, such that for any $0\leq c'\leq c$, it is also true that
\begin{align*}
  &\nu^{\f23} |kB|^{\f13}\|w'\|_{L^2}+\nu^{\f13} |kB|^{\f23}\|\f{w}{r}\|_{L^2} \leq C\|rF-c'\nu^{\f13} |kB|^{\f23}\f{w}{r}\|_{L^2}.
\end{align*}
\end{lemma}
\begin{proof} We first prove \eqref{resolvent estimate L2}. The second inequality readily follows if one selects $c=\frac{1}{2C}$ since for all $c'\in [0, c]$ it holds
\begin{equation*}
    C \|rF\|_{L^2}\le C\|rF-c'\nu^{\f13} |kB|^{\f23}\f{w}{r}\|_{L^2}+\frac{1}{2}\nu^{\f13} |kB|^{\f23}\|\f{w}{r}\|_{L^2}.
\end{equation*}
 By utilizing Lemma \ref{trivial w' lemma} and Cauchy-Schwarz inequality, we obtain
\begin{align}
\label{basic for w'}\nu\|w'\|_{L^2}^2\leq \|rF\|_{L^2}\|\f{w}{r}\|_{L^2}.
\end{align}

To proceed, we then first prove the statement
\begin{equation}\label{w/r rF ineq}
\nu^{\f13} |kB|^{\f23}\|\f{w}{r}\|_{L^2}\leq C \|rF\|_{L^2}.
\end{equation}
The mathematical discussion can be divided into three cases:
\begin{align*}
\lambda\in(-\infty,\f{1-\delta}{R^2}]\cup[1+\delta,\infty),\quad \lambda\in[\f{1-\delta}{R^2},\f{1}{R^2}]\cup[1,1+\delta],\quad \lambda\in[\f{1}{R^2},1].
\end{align*}
Here $0<\delta\ll1$ is a small constant, which will be determined later.

\begin{enumerate}
\item \textbf{Case of $\lambda\in(-\infty,\f{1-\delta}{R^2}]\cup[1+\delta,\infty)$.} One can check
\begin{align*}
|\Im\langle F,w\rangle|=|kB|\big|\langle(\f{1}{r^2}-\lambda)w,w\rangle\big|=|kB|\big|\langle(1-\lambda r^2)\f{w}{r},\f{w}{r}\rangle\big|\geq |kB|\delta \|\f{w}{r}\|_{L^2}^2.
\end{align*}
Taking $\delta=\f{\nu^{\f13}}{|kB|^{\f13}}$, it yields
\begin{align}
\label{lmabda trival-0}\|rF\|_{L^2}\geq\nu^{\f13}|kB|^{\f23}\|\f{w}{r}\|_{L^2}.
\end{align}

\item \textbf{Case of $\lambda\in[\f{1-\delta}{R^2},\f{1}{R^2}]$ or $\lambda\in[1,1+\delta]$.}

i)$\lambda\in[\f{1-\delta}{R^2},\f{1}{R^2}]$. Noting that $\sqrt{\f{1-\delta}{\lambda}}\in [1 ,R]$, we first obtain
\begin{align*}
|\Im\langle F,w\rangle|=&|kB|\int_1^R(1-\lambda r^2)|\f{w}{r}|^2dr\\
\geq&|kB|\int_1^{\sqrt{\f{1-\delta}{\lambda}}}(1-\lambda r^2)|\f{w}{r}|^2dr\geq|kB|\delta \|\f{w}{r}\|_{L^2(1,\sqrt{\f{1-\delta}{\lambda}})}^2.
\end{align*}
The rest part of $\|\f{w}{r}\|_{L^2}$ can be controlled as below
\begin{align*}
&\|\f{w}{r}\|_{L^2(\sqrt{\f{1-\delta}{\lambda}},R)}^2=\int_{\sqrt{\f{1-\delta}{\lambda}}}^R|\f{w}{r}|^2dr\leq\int_{\sqrt{\f{1-\delta}{\lambda}}}^R\f{1}{r^2}dr\|w\|_{L^{\infty}}^2\\
=&(\sqrt{\f{\lambda}{1-\delta}}-\f{1}{R})\|w\|_{L^{\infty}}^2
\leq(\sqrt{\f{1}{1-\delta}}-1)\f{\|w\|_{L^{\infty}}^2}{R}\lesssim\f{\delta}{R}\|w\|_{L^{\infty}}^2\leq\delta\|\f{w}{r^{\f12}}\|_{L^{\infty}}^2.
\end{align*}
Combining these two bounds, we deduce
  \begin{align*}
\|\f{w}{r}\|_{L^2}^2=&\|\f{w}{r}\|_{L^2(1,\sqrt{\f{1-\delta}{\lambda}})}^2+\|\f{w}{r}\|_{L^2(\sqrt{\f{1-\delta}{\lambda}},R)}^2\\
\lesssim&|kB|^{-1}\delta^{-1}|\Im\langle F,w\rangle|+\delta\|\f{w}{r^{\f12}}\|_{L^{\infty}}^2.
\end{align*}
It is then inferred from Lemma \ref{Appendix A1-1} and Lemma \ref{trivial w' lemma} that
  \begin{align*}
\|\f{w}{r}\|_{L^2}^2\lesssim&|kB|^{-1}\delta^{-1}|\Im\langle F,w\rangle|+\delta\|\f{w}{r^{\f12}}\|_{L^{\infty}}^2\\
\lesssim&|kB|^{-1}\delta^{-1}\|rF\|_{L^2}\|\f{w}{r}\|_{L^2}+\delta\|\f{w}{r}\|_{L^2}\|r^{\f12}\partial_r(\f{w}{r^{\f12}})\|_{L^2}\\
\lesssim&|kB|^{-1}\delta^{-1}\|rF\|_{L^2}\|\f{w}{r}\|_{L^2}+\delta\|\f{w}{r}\|_{L^2}(\|w'\|_{L^2}+\|\f{w}{r}\|_{L^2})\\
\lesssim&|kB|^{-1}\delta^{-1}\|rF\|_{L^2}\|\f{w}{r}\|_{L^2}+\nu^{-\f12}\delta\|rF\|_{L^2}^{\f12}\|\f{w}{r}\|_{L^2}^{\f32}.
\end{align*}
This leads to 
\begin{align*}
\|rF\|_{L^2}\gtrsim\min\{|kB|\delta,\nu\delta^{-2}\}\|\f{w}{r}\|_{L^2}.
\end{align*}
The desired estimate \eqref{w/r rF ineq} follows by picking $\delta=\f{\nu^{\f13}}{|kB|^{\f13}}$.

ii)$\lambda\in[1,1+\delta]$. Observing $\sqrt{\f{1+\delta}{\lambda}}\in [1, R]$, we have
\begin{align*}
|\Im\langle F,w\rangle|=&|kB|\int_1^R(\lambda r^2-1)|\f{w}{r}|^2dr\\
\geq&|kB|\int_{\sqrt{\f{1+\delta}{\lambda}}}^R(\lambda r^2-1)|\f{w}{r}|^2dr\geq|kB|\delta \|\f{w}{r}\|_{L^2(\sqrt{\f{1+\delta}{\lambda}},R)}^2.
\end{align*}
When $r\in (1,\sqrt{\f{1+\delta}{\lambda}})$, we bound $\f{w}{r}$  by the $L^\infty$ norm of $\f{w}{r^{\f12}}$ and it holds
\begin{align*}
&\|\f{w}{r}\|_{L^2(1,\sqrt{\f{1+\delta}{\lambda}})}^2=\int_1^{\sqrt{\f{1+\delta}{\lambda}}}|\f{w}{r}|^2dr
\leq\int_1^{\sqrt{\f{1+\delta}{\lambda}}}\f{1}{r}dr\|\f{w}{r^{\f12}}\|_{L^{\infty}}^2\\
=&(\sqrt{\f{1+\delta}{\lambda}}-1)\|\f{w}{r^{\f12}}\|_{L^{\infty}}^2
\leq(\sqrt{1+\delta}-1)\|\f{w}{r^{\f12}}\|_{L^{\infty}}^2\leq\delta\|\f{w}{r^{\f12}}\|_{L^{\infty}}^2.
\end{align*}
Summing two inequalities above renders
  \begin{align*}
\|\f{w}{r}\|_{L^2}^2=&\|\f{w}{r}\|_{L^2(1,\sqrt{\f{1+\delta}{\lambda}})}^2+\|\f{w}{r}\|_{L^2(\sqrt{\f{1+\delta}{\lambda}},R)}^2\\
\lesssim&\delta\|\f{w}{r^{\f12}}\|_{L^{\infty}}^2+|kB|^{-1}\delta^{-1}|\Im\langle F,w\rangle|.
\end{align*}
In view of Lemma \ref{Appendix A1-1} and Lemma \ref{trivial w' lemma}, we  further deduce
  \begin{align*}
\|\f{w}{r}\|_{L^2}^2\lesssim&|kB|^{-1}\delta^{-1}|\Im\langle F,w\rangle|+\delta\|\f{w}{r^{\f12}}\|_{L^{\infty}}^2\\
\lesssim&|kB|^{-1}\delta^{-1}\|rF\|_{L^2}\|\f{w}{r}\|_{L^2}+\delta\|\f{w}{r}\|_{L^2}\|r^{\f12}\partial_r(\f{w}{r^{\f12}})\|_{L^2}\\
\lesssim&|kB|^{-1}\delta^{-1}\|rF\|_{L^2}\|\f{w}{r}\|_{L^2}+\delta\|\f{w}{r}\|_{L^2}(\|w'\|_{L^2}+\|\f{w}{r}\|_{L^2})\\
\lesssim&|kB|^{-1}\delta^{-1}\|rF\|_{L^2}\|\f{w}{r}\|_{L^2}+\nu^{-\f12}\delta\|rF\|_{L^2}^{\f12}\|\f{w}{r}\|_{L^2}^{\f32}.
\end{align*}
This yields
\begin{align*}
\|rF\|_{L^2}\gtrsim\min\{|kB|\delta,\nu\delta^{-2}\}\|\f{w}{r}\|_{L^2}.
\end{align*}
By setting $\delta=\f{\nu^{\f13}}{|kB|^{\f13}}$, we can optimize the inequality above and thus obtain
\begin{align}
\label{lmabda nontrival-2}\|rF\|_{L^2}\gtrsim \nu^{\f13}|kB|^{\f23}\|\f{w}{r}\|_{L^2}.
\end{align}

\item \textbf{Case of $\lambda\in[\f{1}{R^2},1]$.} Using the fact 
\begin{align*}
|1-\lambda r^2|\leq\delta \Rightarrow 1-\delta \leq \lambda r^2\leq 1+\delta \Rightarrow \sqrt{\f{1-\delta}{\lambda}}\leq r\leq \sqrt{\f{1+\delta}{\lambda}},
\end{align*}
for small $\delta>0$, it can be seen that
\begin{align*}
\sqrt{\f{1+\delta}{\lambda}}-\sqrt{\f{1-\delta}{\lambda}}=\f{2\delta }{\sqrt{\lambda}(\sqrt{1+\delta}+\sqrt{1-\delta})}\lesssim\f{\delta }{\sqrt{\lambda}}.
\end{align*}
Now we choose $r_{-}\in(\sqrt{\f{1-\delta}{\lambda}}-\f{\delta }{\sqrt{\lambda}},\sqrt{\f{1-\delta}{\lambda}})$ and $r_{+}\in(\sqrt{\f{1+\delta}{\lambda}},\sqrt{\f{1+\delta}{\lambda}}+\f{\delta }{\sqrt{\lambda}})$ such that the following inequalities hold
\begin{align}
\label{r-,r+1}|w'(r_{-})|^2\leq\f{\|w'\|_{L^2}^2}{\delta /\sqrt{\lambda}},\quad |w'(r_{+})|^2\leq\f{\|w'\|_{L^2}^2}{\delta /\sqrt{\lambda}}.
\end{align}
In order to control $\|\f{w}{r}\|_{L^2((1, r_{-})\cup(r_{+},R))}$, we examine the inner product of $F$ and $w(\chi_{(1,r_{-})}-\chi_{(r_{+},R)})$. Via integration by parts, we write
\begin{align*}
&\langle F,w(\chi_{(1,r_{-})}-\chi_{(r_{+},R)})\rangle\\
=&-\nu\int_1^{r_{-}}w''\overline{w}dr+\nu\int_{r_{+}}^Rw''\overline{w}dr+\nu (k^2-\f14)[\int_1^{r_{-}}\f{|w|^2}{r^2}dr-\int_{r_{+}}^R\f{|w|^2}{r^2}dr]\\
&+ikB\big(\int_1^{r_{-}}(\f{1}{r^2}-\lambda)|w|^2dr+\int_{r_{+}}^R(\lambda-\f{1}{r^2})|w|^2dr\big)\\
=&\nu\int_1^{r_{-}}|w'|^2dr-\nu\int_{r_{+}}^R|w'|^2dr-\nu[w'(r_{-})\overline{w}(r_{-})+w'(r_{+})\overline{w}(r_{+})]\\
&+\nu (k^2-\f14)[\int_1^{r_{-}}\f{|w|^2}{r^2}dr-\int_{r_{+}}^R\f{|w|^2}{r^2}dr]\\
&+ikB\big(\int_1^{r_{-}}(1-\lambda r^2)\f{|w|^2}{r^2}dr+\int_{r_{+}}^R(\lambda r^2-1)\f{|w|^2}{r^2}dr\big).
\end{align*}
Taking the imaginary part of above equality, it indicates that
\begin{align*}
&|kB|\Big|\int_1^{r_{-}}(1-\lambda r^2)\f{|w|^2}{r^2}dr+\int_{r_{+}}^R(\lambda r^2-1)\f{|w|^2}{r^2}dr\Big|\\
\leq&\|rF\|_{L^2}\|\f{w}{r}\|_{L^2}+\nu\big(|w'(r_{-})\overline{w}(r_{-})|+|w'(r_{+})\overline{w}(r_{+})|\big).
\end{align*}
Thanks to the choice of $r_-$ and $r_+$, we can derive the estimate for $\|\f{w}{r}\|_{L^2\big((1, r{-})\cup(r_{+},R)\big)}$ as below
\begin{align}
\label{w/r 1 to r_+ r_+ to R}
&\|\f{w}{r}\|_{L^2((1, r_{-})\cup(r_{+},R))}^2\\
\nonumber
\leq&|kB|^{-1}\delta^{-1}\Big[\|rF\|_{L^2}\|\f{w}{r}\|_{L^2}+\nu\big(|w'(r_{-})\overline{w}(r_{-})|+|w'(r_{+})\overline{w}(r_{+})|\big)\Big].
\end{align}
According to \eqref{r-,r+1}, the second term on the right-hand side of \eqref{w/r 1 to r_+ r_+ to R} can be bounded by
\begin{align*}
&\nu\big(|w'(r_{-})\overline{w}(r_{-})|+|w'(r_{+})\overline{w}(r_{+})|\big)\leq\f{\nu\lambda^{\f14}}{\delta^{\f12}}\|w'\|_{L^2}\Big(|\overline{w}(r_{-})|+|\overline{w}(r_{+})|\Big)\\
\leq&\f{\nu\lambda^{\f14}}{\delta^{\f12}}(r_{-}^{\f12}+r_{+}^{\f12})\|w'\|_{L^2}\|\f{w}{r^{\f12}}\|_{L^{\infty}}\lesssim\f{\nu\lambda^{\f14}}{\delta^{\f12}}\lambda^{-\f14}\|w'\|_{L^2}\|\f{w}{r^{\f12}}\|_{L^{\infty}}=\f{\nu}{\delta^{\f12}}\|w'\|_{L^2}\|\f{w}{r^{\f12}}\|_{L^{\infty}}.
\end{align*}
Thus, we obtain the following estimate for $\|\frac{w}{r}\|_{L^2}$:
\begin{align*}
&\|\f{w}{r}\|_{L^2}^2=\|\f{w}{r}\|_{L^2((1, r_{-})\cup(r_{+},R))}^2+\|\f{w}{r}\|_{L^2(r_{-},r_{+})}^2\\
\leq&|kB|^{-1}\delta^{-1}\Big(\|rF\|_{L^2}\|\f{w}{r}\|_{L^2}+\f{\nu}{\delta^{\f12}}\|w'\|_{L^2}\|\f{w}{r^{\f12}}\|_{L^{\infty}}\Big)+\int_{r_{-}}^{r_{+}}\f{1}{r}dr\|\f{w}{r^{\f12}}\|_{L^\infty}^2 \\
\lesssim&|kB|^{-1}\delta^{-1}\Big(\|rF\|_{L^2}\|\f{w}{r}\|_{L^2}+\f{\nu}{\delta^{\f12}}\|w'\|_{L^2}\|\f{w}{r^{\f12}}\|_{L^{\infty}}\Big)+\delta\|\f{w}{r^{\f12}}\|_{L^\infty}^2,
\end{align*}
where in the last line we employ Lemma \ref{Appendix A1-2} from the Appendix.

 \noindent Utilizing  Lemma \ref{trivial w' lemma} and Lemma \ref{Appendix A1-1}, it can be further deduced that
 \begin{align*}
\|\f{w}{r}\|_{L^2}^2
\lesssim&|kB|^{-1}\delta^{-1}\Big(\|rF\|_{L^2}\|\f{w}{r}\|_{L^2}+\f{\nu}{\delta^{\f12}}\|w'\|_{L^2}\|\f{w}{r^{\f12}}\|_{L^{\infty}}\Big)+\delta\|\f{w}{r^{\f12}}\|_{L^\infty}^2\\
\lesssim&|kB|^{-1}\delta^{-1}\Big(\|rF\|_{L^2}\|\f{w}{r}\|_{L^2}+\f{\nu^{\f12}}{\delta^{\f12}}\|rF\|_{L^2}^{\f12}\|\f{w}{r}\|_{L^2}^{\f12}\|\f{w}{r^{\f12}}\|_{L^{\infty}}\Big)+\delta\|\f{w}{r^{\f12}}\|_{L^\infty}^2\\
\lesssim&|kB|^{-1}\delta^{-1}\|rF\|_{L^2}\|\f{w}{r}\|_{L^2}+|kB|^{-2}\delta^{-4}\nu\|rF\|_{L^2}\|\f{w}{r}\|_{L^2}+\delta\|\f{w}{r^{\f12}}\|_{L^\infty}^2\\
\lesssim&|kB|^{-1}\delta^{-1}\|rF\|_{L^2}\|\f{w}{r}\|_{L^2}+|kB|^{-2}\delta^{-4}\nu\|rF\|_{L^2}\|\f{w}{r}\|_{L^2}+\delta\|\f{w}{r}\|_{L^2}\|r^{\f12}\partial_r(\f{w}{r^{\f12}})\|_{L^2}\\
\lesssim&|kB|^{-1}\delta^{-1}\|rF\|_{L^2}\|\f{w}{r}\|_{L^2}+|kB|^{-2}\delta^{-4}\nu\|rF\|_{L^2}\|\f{w}{r}\|_{L^2}\\
&+\delta\|\f{w}{r}\|_{L^2}(\|w'\|_{L^2}+\|\f{w}{r}\|_{L^2})\\
\lesssim&|kB|^{-1}\delta^{-1}\|rF\|_{L^2}\|\f{w}{r}\|_{L^2}+|kB|^{-2}\delta^{-4}\nu\|rF\|_{L^2}\|\f{w}{r}\|_{L^2}+\delta\nu^{-\f12}\|rF\|_{L^2}^{\f12}\|\f{w}{r}\|_{L^2}^{\f32}.
\end{align*}
This leads to the conclusion
\begin{align*}
\|rF\|_{L^2}\gtrsim\min\{|kB|\delta,|kB|^{2}\delta^{4}\nu^{-1},\delta^{-2}\nu\}\|\f{w}{r}\|_{L^2}.
\end{align*}
One can obtain an optimized form of the above inequality if  $\delta$ is chosen as $\delta=\frac{\nu^{\frac{1}{3}}}{|kB|^{\frac{1}{3}}}$. It then follows that
\begin{align*}
\|rF\|_{L^2}\gtrsim \nu^{\f13}|kB|^{\f23}\|\f{w}{r}\|_{L^2}.
\end{align*}
\end{enumerate}
Therefore, we arrive at
\begin{align}
\label{k.6-}
\nu^{\f13}|kB|^{\f23}\|\f{w}{r}\|_{L^2}\lesssim\|rF\|_{L^2} \quad \text{for all }  \lambda\in \mathbb{R}.
\end{align}
Notice that Lemma \ref{trivial w' lemma} also implies
  \begin{align*}
&\Re\langle F,w\rangle\geq\nu\|w'\|_{L^2}^2.
\end{align*}
Along with \eqref{k.6-} this yields the desired estimate on $w'$:
\begin{align}
\label{w'-reslovent}\nu^{\f23}|kB|^{\f13}\|w'\|_{L^2}\leq C\|rF\|_{L^2}.
\end{align}
The remaining task is to control $\|r(\f{1}{r^2}-\lambda)w\|_{L^2}$, and we appeal to exploring the imaginary part of $\langle F,r^2(\f{1}{r^2}-\lambda)w\rangle$. Integration by parts gives
\begin{align*}
&\langle F,r^2(\f{1}{r^2}-\lambda)w\rangle=\langle-\nu(\partial_r^2-\f{k^2-\f14}{r^2})w+ikB(\f{1}{r^2}-\lambda)w,r^2(\f{1}{r^2}-\lambda)w\rangle\\
=&\nu\Big(\langle w',[(1-\lambda r^2)w]'\rangle+(k^2-\f14)\langle w,(\f{1}{r^2}-\lambda)w\rangle\Big)+ikB\|r(\f{1}{r^2}-\lambda)w\|_{L^2}^2\\
=&\nu\Big(\langle w',(1-\lambda r^2)w'\rangle-\langle w',2\lambda rw\rangle+(k^2-\f14)\langle w,(\f{1}{r^2}-\lambda)w\rangle\Big)+ikB\|r(\f{1}{r^2}-\lambda)w\|_{L^2}^2.
\end{align*}
Then we take the imaginary part of both sides and rewrite $-\lambda rw = r(\f{1}{r^2}-\lambda)w-\f{w}{r}$ to deduce
\begin{align*}
&|kB|\|r(\f{1}{r^2}-\lambda)w\|_{L^2}^2\leq |\langle F,r^2(\f{1}{r^2}-\lambda)w\rangle|+\nu|\langle w',-2\lambda rw\rangle|\\
\leq&|\langle F,r^2(\f{1}{r^2}-\lambda)w\rangle|+\nu|\langle w',2r(\f{1}{r^2}-\lambda)w\rangle|+2\nu|\langle w',\f{w}{r}\rangle|\\
\leq&\|rF\|_{L^2}\|r(\f{1}{r^2}-\lambda)w\|_{L^2}+2\nu\|w'\|_{L^2}\|r(\f{1}{r^2}-\lambda)w\|_{L^2}+2\nu\|w'\|_{L^2}\|\f{w}{r}\|_{L^2}.
\end{align*}
In view of \eqref{k.6-} and \eqref{w'-reslovent}, the $L^2$ norm of $r(\f{1}{r^2}-\lambda)w$ obeys the following bounds:
\begin{align*}
&|kB|\|r(\f{1}{r^2}-\lambda)w\|_{L^2}^2\\
\lesssim & \f{1}{|kB|}\|rF\|_{L^2}^2+\f{\nu^2}{|kB|}\|w'\|_{L^2}^2+\f{1}{|kB|}\nu^{\f23}|kB|^{\f13}\|w'\|_{L^2}\nu^{\f13}|kB|^{\f23}\|\f{w}{r}\|_{L^2}\\
\lesssim & \f{1}{|kB|}\|rF\|_{L^2}^2+\f{\nu^2}{|kB|}\|w'\|_{L^2}^2+\f{1}{|kB|}\|rF\|_{L^2}^2\lesssim \f{1}{|kB|}\|rF\|_{L^2}^2+\f{\nu^2}{|kB|}\|w'\|_{L^2}^2\\
=&\f{1}{|kB|}\|rF\|_{L^2}^2+\f{\nu}{|kB|}\|rF\|_{L^2}\|\f{w}{r}\|_{L^2}
\lesssim  \f{1}{|kB|}\|rF\|_{L^2}^2+\f{\nu^{\f23}}{|kB|^{\f53}}\|rF\|_{L^2}\nu^{\f13}|kB|^{\f23}\|\f{w}{r}\|_{L^2}\\
\lesssim& (1+(\f{\nu}{|kB|})^{\f23})\f{1}{|kB|}\|rF\|_{L^2}^2\lesssim \f{1}{|kB|}\|rF\|_{L^2}^2,
\end{align*}
provided that $\f{\nu}{|kB|} \ll 1$. This completes the proof.
\end{proof}

Relying on the above $L^2$ estimates for $w$ and $w'$, we can derive the below bounds for $\varphi$ and $\varphi'$.
\begin{lemma}\label{resolvent estimate-2}
For any $|k|\geq1$, $\lambda\in\mathbb{R}$ and $w\in D_k$, there exist constants $C,c>0$ independent of $\nu,k,B,\lambda,R$, such that for any $0\leq c'\leq c$, it holds
\begin{align*}
&\nu^{\f16}|kB|^{\f56}|k|^{\f12}\left( \|\varphi'\|_{L^2}+|k|\|\f{\varphi}{r}\|_{L^2} \right) 
\leq CR^2\left((\f{\nu}{|kB|})^{\f16} (\log R)^{\f12} +1 \right) \|rF-c'\nu^{\f13} |kB|^{\f23}\f{w}{r}\|_{L^2}.
\end{align*}
\end{lemma}
\begin{proof}
Denote $\tilde{F}:=rF-c'\nu^{\f13} |kB|^{\f23}\f{w}{r}$. For $\delta=(\f{\nu}{|kB|})^{\f13}\ll 1$, we set 
$$E=\{r> 0:  |1-\lambda r^2|\le \delta \} \quad \textrm{and} \quad E^c=(0, \infty)\backslash E.$$ 
Employing Lemma \ref{resolvent estimate-1}, Lemma \ref{Appendix A1-1} and Lemma \ref{Appendix A5}, we have
\begin{align*}
    \|\f{w}{r^{\f32}}\|_{L^1([1, R] \cap E)}\le& \|\f{w}{r^{\f12}}\|_{L^\infty} \|\f{1}{r}\|_{L^1(E)} \lesssim
    \delta\|\f{w}{r}\|_{L^2}^{\f12}(\|\f{w}{r}\|_{L^2}+\|w'\|_{L^2})^{\f12} \\
    \lesssim&  \nu^{-\f16}|kB|^{-\f56} \|\tilde{F}\|_{L^2}.
\end{align*}
On the other hand, utilizing  Lemma \ref{resolvent estimate-1} again as well as Lemma \ref{Appendix A6}, it can be inferred that
\begin{align*}
    \|\f{w}{r^{\f32}}\|_{L^1([1, R] \cap E^c)}\le& \|\f{1}{r^{\f12}(1-\lambda r^2)^2}\|_{L^2([1, R]\cap E^c)}\|r(\f{1}{r^2}-\lambda)w\|_{L^2} \\
    \lesssim& \left( (\log R)^{\f12}+\delta^{-\f12} \right) |kB|^{-1} \|\tilde{F}\|_{L^2} \\
    \lesssim&   \nu^{-\f16}|kB|^{-\f56}\left((\f{\nu}{|kB|})^{\f16} (\log R)^{\f12} +1 \right) \|\tilde{F}\|_{L^2}.
\end{align*}
Thus we deduce
\begin{align*}
   \|\f{w}{r^{\f32}}\|_{L^1}=& \|\f{w}{r^{\f32}}\|_{L^1([1, R] \cap E)}+  \|\f{w}{r^{\f32}}\|_{L^1([1, R] \cap E^c)} \\
   \lesssim& \nu^{-\f16}|kB|^{-\f56}\left((\f{\nu}{|kB|})^{\f16} (\log R)^{\f12} +1 \right) \|\tilde{F}\|_{L^2}.
\end{align*}
The controls for $\varphi'$ and $\f{\varphi}{r}$ also follow from Lemma \ref{Appendix A4}. And we have
\begin{align*}
 \|\varphi'\|_{L^2}+|k|\|\f{\varphi}{r}\|_{L^2}\lesssim& |k|^{-\f12}\|r^{\f12}w\|_{L^1}\lesssim |k|^{-\f12}R^2\|\f{w}{r^{\f32}}\|_{L^1} \\
    \lesssim& |k|^{-\f12}R^2\nu^{-\f16}|kB|^{-\f56}\left((\f{\nu}{|kB|})^{\f16} (\log R)^{\f12} +1 \right) \|\tilde{F}\|_{L^2}.
\end{align*}
This completes the proof.
\end{proof}
Now we turn to estimate the $H^{-1}_r$ norm of $w$ through the $H^1_r$ norm of $F$.

\begin{lemma}\label{resolvent estimate-3}For any $|k|\geq1$, $\lambda\in\mathbb{R}$ and $w\in D_k$, there exist a constants $C,c>0$ independent of $\nu,k,B,\lambda$, such that for $0\leq c'\leq c$, the following estimate holds
\begin{align*}
\nu\|w\|_{H^1_r}+\nu^{\f23}|kB|^{\f13}\|\f{w}{r}\|_{L^2}\leq C\|F-c'\nu^{\f13} |kB|^{\f23} R^{-2}w\|_{H^{-1}_r}.
\end{align*}
\end{lemma}
\begin{proof}By Lemma \ref{trivial w' lemma}, one first obtains
\begin{align*}
&\nu\|w'\|_{L^2}^2+\nu(k^2-\f14)\|\f{w}{r}\|_{L^2}^2\\
=&\Re\langle F-c'\nu^{\f13} |kB|^{\f23} R^{-2}w,w\rangle+c'\nu^{\f13} |kB|^{\f23} R^{-2}\|w\|_{L^2}^2,
\end{align*}
which gives
\begin{align}
\label{basic for w'-1-0}\nu\|w\|_{H^1_r}\lesssim \|F-c'\nu^{\f13} |kB|^{\f23}R^{-2}w\|_{H^{-1}_r}+\sqrt{c'}\nu^{\f23} |kB|^{\f13}\|\f{w}{r}\|_{L^2}.
\end{align}
Denote $G:=F-c'\nu^{\f13} |kB|^{\f23}R^{-2}w$. In below, we utilize a similar method as in the proof of Proposition \ref{resolvent estimate-1} to demonstrate
\begin{equation}\label{w over r bound by G}
    \nu^{\f23} |kB|^{\f13}\|\f{w}{r}\|_{L^2}\leq C \|G\|_{H^{-1}_r}.
\end{equation}
To achieve this, we will examine the following three cases:
\begin{align*}
\lambda\in(-\infty,\f{1-\delta}{R^2}]\cup[1+\delta,\infty),\quad \lambda\in[\f{1-\delta}{R^2},\f{1}{R^2}]\cup[1,1+\delta],\quad \lambda\in[\f{1}{R^2},1],
\end{align*}
where $0<\delta\ll1$ is a small constant to be determined later.
\begin{enumerate}
\item \textbf{Case of $\lambda\in(-\infty,\f{1-\delta}{R^2}]\cup[1+\delta,\infty)$.} One can check
\begin{align*}
&|\Im\langle G,w\rangle|=|\Im\langle F,w\rangle|=|kB|\big|\langle(\f{1}{r^2}-\lambda)w,w\rangle\big|=|kB|\big|\langle(1-\lambda r^2)\f{w}{r},\f{w}{r}\rangle\big|\geq |kB|\delta \|\f{w}{r}\|_{L^2}^2.
\end{align*}
Taking $\delta=\f{\nu^{\f13}}{|kB|^{\f13}}$, together with \eqref{basic for w'-1-0}, we obtain
\begin{align*}
C(\|G\|_{H^{-1}_r}+c'\nu^{\f23} |kB|^{\f13}\|\f{w}{r}\|_{L^2})\geq\nu^{\f23} |kB|^{\f13}\|\f{w}{r}\|_{L^2}.
\end{align*}
Choosing $Cc'\leq\f12$, we then prove \eqref{w over r bound by G}.

\item \textbf{Case of $\lambda\in[\f{1-\delta}{R^2},\f{1}{R^2}]$ or $\lambda\in[1,1+\delta]$.}

i) $\lambda\in[\f{1-\delta}{R^2},\f{1}{R^2}]$. Observing that
\begin{align*}
1-\lambda r^2\geq\delta\Longleftrightarrow 1-\delta\geq\lambda r^2\Longleftrightarrow 1\leq r\leq\sqrt{\f{1-\delta}{\lambda}},
\end{align*}
we then deduce
\begin{align*}
|\Im\langle G,w\rangle|=&|\Im\langle F,w\rangle|=|kB|\int_1^R(1-\lambda r^2)|\f{w}{r}|^2dr\\
\geq&|kB|\int_1^{\sqrt{\f{1-\delta}{\lambda}}}(1-\lambda r^2)|\f{w}{r}|^2dr\geq|kB|\delta \|\f{w}{r}\|_{L^2(1,\sqrt{\f{1-\delta}{\lambda}})}^2.
\end{align*}
The remaining part of $\|\f{w}{r}\|_{L^2}$ can be bounded in terms of $\|\f{w}{r^{\f12}}\|_{L^{\infty}}$ as below
\begin{align*}
&\|\f{w}{r}\|_{L^2(\sqrt{\f{1-\delta}{\lambda}},R)}^2=\int_{\sqrt{\f{1-\delta}{\lambda}}}^R|\f{w}{r}|^2dr\leq\int_{\sqrt{\f{1-\delta}{\lambda}}}^R\f{1}{r^2}dr\|w\|_{L^{\infty}}^2\\
=&(\sqrt{\f{\lambda}{1-\delta}}-\f{1}{R})\|w\|_{L^{\infty}}^2
\leq(\sqrt{\f{1}{1-\delta}}-1)\f{\|w\|_{L^{\infty}}^2}{R}\lesssim\f{\delta}{R}\|w\|_{L^{\infty}}^2\leq\delta\|\f{w}{r^{\f12}}\|_{L^{\infty}}^2.
\end{align*}
Combining these two estimates, together with \eqref{basic for w'-1-0},  we get
  \begin{align*}
\|\f{w}{r}\|_{L^2}^2=&\|\f{w}{r}\|_{L^2(1,\sqrt{\f{1-\delta}{\lambda}})}^2+\|\f{w}{r}\|_{L^2(\sqrt{\f{1-\delta}{\lambda}},R)}^2\\
\lesssim&|kB|^{-1}\delta^{-1}(\nu^{-1}\|G\|_{H^{-1}_r}^2+c'\nu^{\f13}|kB|^{\f23}\|\f{w}{r}\|_{L^2}^2)+\delta\|\f{w}{r^{\f12}}\|_{L^{\infty}}^2.
\end{align*}
Applying Lemma \ref{Appendix A1-1} and \eqref{basic for w'-1-0}, we provide the estimate for $\|\f{w}{r}\|_{L^2}$ by
  \begin{align*}
\|\f{w}{r}\|_{L^2}^2\lesssim&|kB|^{-1}\delta^{-1}(\nu^{-1}\|G\|_{H^{-1}_r}^2+c'\nu^{\f13}|kB|^{\f23}\|\f{w}{r}\|_{L^2}^2)+\delta\|\f{w}{r^{\f12}}\|_{L^{\infty}}^2\\
\lesssim&|kB|^{-1}\delta^{-1}(\nu^{-1}\|G\|_{H^{-1}_r}^2+c'\nu^{\f13}|kB|^{\f23}\|\f{w}{r}\|_{L^2}^2)+\delta\|\f{w}{r}\|_{L^2}\|r^{\f12}\partial_r(\f{w}{r^{\f12}})\|_{L^2}\\
\lesssim&|kB|^{-1}\delta^{-1}(\nu^{-1}\|G\|_{H^{-1}_r}^2+c'\nu^{\f13}|kB|^{\f23}\|\f{w}{r}\|_{L^2}^2)\\
&+\delta\|\f{w}{r}\|_{L^2}(\|w'\|_{L^2}+\|\f{w}{r}\|_{L^2})\\
\lesssim&|kB|^{-1}\delta^{-1}(\nu^{-1}\|G\|_{H^{-1}_r}^2+c'\nu^{\f13}|kB|^{\f23}\|\f{w}{r}\|_{L^2}^2)\\
&+\delta\|\f{w}{r}\|_{L^2}(\nu^{-\f12}\|G\|_{H^{-1}_r}+\sqrt{c'}\nu^{\f16}|kB|^{\f13}\|\f{w}{r}\|_{L^2})+\delta\|\f{w}{r}\|_{L^2}^2 \\
\lesssim&|kB|^{-1}\delta^{-1}(\nu^{-1}\|G\|_{H^{-1}_r}^2+c'\nu^{\f13}|kB|^{\f23}\|\f{w}{r}\|_{L^2}^2)+|kB|\delta^3\|\f{w}{r}\|_{L^2}^2+\delta\|\f{w}{r}\|_{L^2}^2.
\end{align*}
\begin{comment}
In other words there exists some constant $C>0$ so that
  \begin{align*}
\|\f{w}{r}\|_{L^2}^2
\leq&C|kB|^{-1}\delta^{-1}\nu^{-1}\|G\|_{H^{-1}_r}^2+Cc'\|\f{w}{r}\|_{L^2}^2+C|kB|\delta^3\|\f{w}{r}\|_{L^2}+C\delta\|\f{w}{r}\|_{L^2}^2.
\end{align*}
\end{comment}
Picking $\delta=\f{\nu^{\f13}}{|kB|^{\f13}}$, we then obtain
  \begin{align*}
&\|\f{w}{r}\|_{L^2}^2
\leq C\nu^{-\f43}|kB|^{-\f23}\|G\|_{H^{-1}_r}^2+Cc'\|\f{w}{r}\|_{L^2}^2+C\nu\|\f{w}{r}\|_{L^2}^2+C\delta\|\f{w}{r}\|_{L^2}^2.
\end{align*}
The additional $\|\frac{w}{r}\|_{L^2}$ terms on the right can be absorbed to the left by choosing a sufficiently small $c'>0$, and noting that $0<\nu, \delta\ll 1$. Hence, we also get the desired estimate \eqref{w over r bound by G}.

ii) $\lambda\in[1,1+\delta]$. Using that 
\begin{align*}
\lambda r^2-1\geq\delta\Longleftrightarrow r^2\geq\f{1+\delta}{\lambda} \Longleftrightarrow \sqrt{\f{1+\delta}{\lambda}}\leq r\leq R,
\end{align*}
we can bound the imaginary part of $\langle G,w\rangle$ as below
\begin{align*}
&|\Im\langle G,w\rangle|=|\Im\langle F,w\rangle|=|kB|\int_1^R(\lambda r^2-1)|\f{w}{r}|^2dr\\
\geq&|kB|\int_{\sqrt{\f{1+\delta}{\lambda}}}^R(\lambda r^2-1)|\f{w}{r}|^2dr\geq|kB|\delta \|\f{w}{r}\|_{L^2(\sqrt{\f{1+\delta}{\lambda}},R)}^2.
\end{align*}
We further control $\|\f{w}{r}\|_{L^2(1,\sqrt{\f{1+\delta}{\lambda}})}$ via
\begin{align*}
\|\f{w}{r}\|_{L^2(1,\sqrt{\f{1+\delta}{\lambda}})}^2=&\int_1^{\sqrt{\f{1+\delta}{\lambda}}}|\f{w}{r}|^2dr
\leq\int_1^{\sqrt{\f{1+\delta}{\lambda}}}\f{1}{r}dr\|\f{w}{r^{\f12}}\|_{L^{\infty}}^2\\
=&(\sqrt{\f{1+\delta}{\lambda}}-1)\|\f{w}{r^{\f12}}\|_{L^{\infty}}^2
\leq(\sqrt{1+\delta}-1)\|\f{w}{r^{\f12}}\|_{L^{\infty}}^2\leq\delta\|\f{w}{r^{\f12}}\|_{L^{\infty}}^2.
\end{align*}
Together with \eqref{basic for w'-1-0}, this implies
  \begin{align*}
\|\f{w}{r}\|_{L^2}^2=&\|\f{w}{r}\|_{L^2(1,\sqrt{\f{1+\delta}{\lambda}})}^2+\|\f{w}{r}\|_{L^2(\sqrt{\f{1+\delta}{\lambda}},R)}^2\\
\lesssim&\delta\|\f{w}{r^{\f12}}\|_{L^{\infty}}^2+|kB|^{-1}\delta^{-1}(\nu^{-1}\|G\|_{H^{-1}_r}^2+c'\nu^{\f13}|kB|^{\f23}\|\f{w}{r}\|_{L^2}^2).
\end{align*}

Thus, we obtain the same bound for $\|\frac{w}{r}\|_{L^2}$ as in the scenario with $\lambda\in[\frac{1-\delta}{R^2},\frac{1}{R^2}]$. Applying Lemma \ref{Appendix A1-1}, inequality \eqref{basic for w'-1-0} and choosing $\delta=\frac{\nu^{\frac13}}{|kB|^{\frac13}}$ and $c'>0$ sufficiently small, we then derive the desired bound \eqref{w over r bound by G}.
\item \textbf{Case of $\lambda\in[\f{1}{R^2},1]$.} We now use
\begin{equation*}
|1-\lambda r^2|\leq\delta \Rightarrow 1-\delta \leq \lambda r^2\leq 1+\delta \Rightarrow \sqrt{\f{1-\delta}{\lambda}}\leq r\leq \sqrt{\f{1+\delta}{\lambda}}.
\end{equation*}
Given $0<\delta\ll1$, it then holds
\begin{align*}
\sqrt{\f{1+\delta}{\lambda}}-\sqrt{\f{1-\delta}{\lambda}}=\f{2\delta }{\sqrt{\lambda}(\sqrt{1+\delta}+\sqrt{1-\delta})}\lesssim\f{\delta }{\sqrt{\lambda}}.
\end{align*}
We further choose $r_{-}\in(\sqrt{\f{1-\delta}{\lambda}}-\f{\delta }{\sqrt{\lambda}},\sqrt{\f{1-\delta}{\lambda}})$ and $r_{+}\in(\sqrt{\f{1+\delta}{\lambda}},\sqrt{\f{1+\delta}{\lambda}}+\f{\delta }{\sqrt{\lambda}})$ satisfying the following inequalities
\begin{align}
\label{r-,r+1-h-decay}|w'(r_{-})|^2\leq\f{\|w'\|_{L^2}^2}{\delta /\sqrt{\lambda}},\quad |w'(r_{+})|^2\leq\f{\|w'\|_{L^2}^2}{\delta /\sqrt{\lambda}}.
\end{align}
For the next step, we will construct an appropriate multiplier. To do so,  we first define a piecewise $C^1$ cutoff function $\rho$ with domain $(1, R)$ as follows
\begin{align*}
\rho(r)=\left\{
\begin{aligned}
&1,\quad r\in(1,r_{-}-\f{\delta }{\sqrt{\lambda}}),\\
&\sin\big(\f{\pi}{2}\f{\sqrt{\lambda}}{\delta}(r_{-}-r)\big),\quad r\in(r_{-}-\f{\delta }{\sqrt{\lambda}},r_{-}),\\
&0,\quad r\in (r_{-},r_{+}),\\
&\sin\big(\f{\pi}{2}\f{\sqrt{\lambda}}{\delta}(r_{+}-r)\big),\quad r\in(r_{+},r_{+}+\f{\delta }{\sqrt{\lambda}}),\\
&-1,\quad r\in(r_{+}+\f{\delta }{\sqrt{\lambda}},R).
\end{aligned}
\right.
\end{align*}
Via integration by parts, we obtain the following expression
\begin{align*}
-\Im\langle G,w\rho\rangle=\Im\langle F,w\rho\rangle=&\Im\langle -\nu\partial_r^2w+ikB(\f{1}{r^2}-\lambda)w,w\rho\rangle\\
=&\Im\langle ikB(\f{1}{r^2}-\lambda)w,w\rho\rangle+\nu\Im\langle w',w\rho'\rangle.
\end{align*}
Taking the imaginary part of above equality leads to
\begin{align*}
&|kB|\Big|\int_1^{r_{-}-\f{\delta }{\sqrt{\lambda}}}(1-\lambda r^2)\f{|w|^2}{r^2}dr+\int_{r_{+}+\f{\delta }{\sqrt{\lambda}}}^R(\lambda r^2-1)\f{|w|^2}{r^2}dr\Big|\\
\leq&\|G\|_{H^{-1}_r}\|w\rho\|_{H^1_r}+\nu\|w'\|_{L^2}\|w\rho'\|_{L^2}.
\end{align*}
Based on our choice of $\rho$ and the definition of the $H_r^1$ norm, it follows that
\begin{align*}
    \|w\rho\|_{L^2}\lesssim \|(w\rho)'\|_{L^2}+\|\f{w}\rho{r}\|_{L^2} \le& \|w'\|_{L^2}+\|w\rho'\|_{L^2}+\|\f{w}{r}\|_{L^2},
\end{align*}
and
\begin{align*}
    \|w\rho'\|_{L^2}\lesssim\f{\sqrt{\lambda}}{\delta}\|w \|_{L^2\big((r_{-}-\f{\delta }{\sqrt{\lambda}},r_{-})\cup(r_{+},r_{+}+\f{\delta }{\sqrt{\lambda}})\big)} 
    \lesssim \delta^{-1} \|\f{w}{r}\|_{L^2\big((r_--\f{\delta}{\sqrt{\lambda}}, r_++\f{\delta}{\sqrt{\lambda}})\big)}.
\end{align*}
Hence, plugging in  \eqref{basic for w'-1-0} and \eqref{r-,r+1-h-decay}, we then derive the estimate 
\begin{align*}
&\|\f{w}{r}\|_{L^2\big((1, r_{-}-\f{\delta }{\sqrt{\lambda}})\cup(r_{+}+\f{\delta }{\sqrt{\lambda}},R)\big)}^2\\
\lesssim&|kB|^{-1}\delta^{-1}\Big[\|G\|_{H^{-1}_r}(\|w'\|_{L^2}+\delta^{-1} \|\f{w}{r}\|_{L^2\big((r_--\f{\delta}{\sqrt{\lambda}}, r_++\f{\delta}{\sqrt{\lambda}})\big)})+\nu\|w'\|_{L^2}\cdot\delta^{-1}\|\f{w}{r}\|_{L^2}\Big]\\
\lesssim&|kB|^{-1}\delta^{-1}\Big[\|G\|_{H^{-1}_r}(\nu^{-1}\|G\|_{H^{-1}_r}+\nu^{-1}\sqrt{c'}\nu^{\f23}|kB|^{\f13}\|\f{w}{r}\|_{L^2}+\delta^{-1} \|\f{w}{r}\|_{L^2\big((r_--\f{\delta}{\sqrt{\lambda}}, r_++\f{\delta}{\sqrt{\lambda}})\big)})\\&\quad\quad\quad\quad\quad+(\|G\|_{H^{-1}_r}+\sqrt{c'}\nu^{\f23}|kB|^{\f13}\|\f{w}{r}\|_{L^2})\cdot \delta^{-1} \|\f{w}{r}\|_{L^2\big((r_--\f{\delta}{\sqrt{\lambda}}, r_++\f{\delta}{\sqrt{\lambda}})\big)}\Big]\\
\lesssim&|kB|^{-1}\delta^{-1}\Big(\nu^{-1}\|G\|_{H^{-1}_r}^2+c'\nu^{\f13}|kB|^{\f23}\|\f{w}{r}\|_{L^2}^2+\nu\delta^{-2} \|\f{w}{r}\|_{L^2\big((r_--\f{\delta}{\sqrt{\lambda}}, r_++\f{\delta}{\sqrt{\lambda}})\big)}\Big).
\end{align*}
Meanwhile, it can be deduced from Lemma \ref{Appendix A1-2}, Lemma \ref{Appendix A1-1} and \eqref{basic for w'-1-0} that
\begin{align*}
    \|\f{w}{r}\|_{L^2\big((r_--\f{\delta}{\sqrt{\lambda}}, r_++\f{\delta}{\sqrt{\lambda}})\big)}^2\le& \int_{r_{-}-\f{\delta }{\sqrt{\lambda}}}^{r_{+}+\f{\delta }{\sqrt{\lambda}}}\f{1}{r}dr\|\f{w}{r^{\f12}}\|_{L^\infty}^2
    \lesssim \delta\|\f{w}{r}\|_{L^2}\|r^{\f12}\partial_r(\f{w}{r^{\f12}})\|_{L^2} \\
     \lesssim& \delta\|\f{w}{r}\|_{L^2}(\|w'\|_{L^2}+\|\f{w}{r}\|_{L^2})\\
     \lesssim& \delta\|\f{w}{r}\|_{L^2}(\nu^{-1}\|G\|_{H^{-1}_r}+\nu^{-1}\sqrt{c'}\nu^{\f23}|kB|^{\f13}\|\f{w}{r}\|_{L^2}+\|\f{w}{r}\|_{L^2})\\
     \lesssim& \f{1}{\sqrt{c'}}\nu^{-2}\delta^{2}\|G\|_{H^{-1}_r}^2+\sqrt{c'}(\nu^{-\f13}|kB|^{\f13}\delta+1)\|\f{w}{r}\|_{L^2}^2+\delta\|\f{w}{r}\|_{L^2}^2.
\end{align*}
Taking $\delta=\f{\nu^{\f13}}{|kB|^{\f13}}$, we can then bound the $L^2$ norm of $\f{w}{r}$ as below
\begin{align*}
\|\f{w}{r}\|_{L^2}^2=&\|\f{w}{r}\|_{L^2\big((1, r_{-}-\f{\delta }{\sqrt{\lambda}})\cup(r_{+}+\f{\delta }{\sqrt{\lambda}},R)\big)}^2+\|\f{w}{r}\|_{L^2\big(( r_{-}-\f{\delta }{\sqrt{\lambda}},r_{+}+\f{\delta }{\sqrt{\lambda}})\big)}^2\\
\lesssim&|kB|^{-1}\delta^{-1}\Big(\nu^{-1}\|G\|_{H^{-1}_r}^2+c'\nu^{\f13}|kB|^{\f23}\|\f{w}{r}\|_{L^2}^2\Big)\\
&+ \|\f{w}{r}\|_{L^2\big((r_--\f{\delta}{\sqrt{\lambda}}, r_++\f{\delta}{\sqrt{\lambda}})\big)}^2\\
\lesssim& (1+\f{1}{\sqrt{c'}})\nu^{-\f43}|kB|^{-\f23}\|G\|_{H^{-1}_r}^2+(c'+\sqrt{c'}+\delta)\|\f{w}{r}\|_{L^2}^2.
\end{align*}
\begin{comment}
take $\delta=\f{\nu^{\f13}}{|kB|^{\f13}}$, combining with Lemma \ref{Appendix A1-2}, we have
\begin{align*}
\|\f{w}{r}\|_{L^2}^2
\lesssim&|kB|^{-1}\delta^{-1}\Big(\nu^{-1}\|G\|_{H^{-1}_r}^2+c'\nu^{\f13}|kB|^{\f23}\|\f{w}{r}\|_{L^2}^2\Big)+\delta\|\f{w}{r^{\f12}}\|_{L^\infty}^2\\
\lesssim&\nu^{-\f43}|kB|^{-\f23}\|G\|_{H^{-1}_r}^2+c'\|\f{w}{r}\|_{L^2}^2+\delta\|\f{w}{r^{\f12}}\|_{L^\infty}^2.
\end{align*}
 Together with Lemma \ref{Appendix A1-1} and \eqref{basic for w'-1-0}, this gives
 \begin{align*}
\|\f{w}{r}\|_{L^2}^2
\lesssim&\nu^{-\f43}|kB|^{-\f23}\|G\|_{H^{-1}_r}^2+c'\|\f{w}{r}\|_{L^2}^2+\delta\|\f{w}{r}\|_{L^2}\|r^{\f12}\partial_r(\f{w}{r^{\f12}})\|_{L^2}\\
\lesssim&\nu^{-\f43}|kB|^{-\f23}\|G\|_{H^{-1}_r}^2+c'\|\f{w}{r}\|_{L^2}^2
+\delta\|\f{w}{r}\|_{L^2}(\|w'\|_{L^2}+\|\f{w}{r}\|_{L^2})\\
\lesssim&\nu^{-\f43}|kB|^{-\f23}\|G\|_{H^{-1}_r}^2+c'\|\f{w}{r}\|_{L^2}^2+\delta\|\f{w}{r}\|_{L^2}\nu^{-1}(\|G\|_{H^{-1}_r}+c'\nu^{\f23}|kB|^{\f13}\|\f{w}{r}\|_{L^2})\\
\lesssim&\nu^{-\f43}|kB|^{-\f23}\|G\|_{H^{-1}_r}^2+c'\|\f{w}{r}\|_{L^2}^2\leq C\nu^{-\f43}|kB|^{-\f23}\|G\|_{H^{-1}_r}^2+Cc'\|\f{w}{r}\|_{L^2}^2.
\end{align*}
\end{comment}
Picking $c'>0$ sufficiently small and noting that $0<\delta\ll 1$, we can proceed to the conclusion
\begin{align}
\label{lmabda nontrival-3-h-decay}\|G\|_{H^{-1}_r}\gtrsim \nu^{\f23}|kB|^{\f13}\|\f{w}{r}\|_{L^2}.
\end{align}
\end{enumerate}
We have therefore proved \eqref{w over r bound by G} in all three scenarios. This together with \eqref{basic for w'-1-0} yields the desired bounds for $\|w\|_{H^{1}_r}$:
\begin{align*}
\nu\|w\|_{H^{1}_r}\lesssim \|G\|_{H^{-1}_r}+\nu^{\f23} |kB|^{\f13}\|\f{w}{r}\|_{L^2}\lesssim \|G\|_{H^{-1}_r} .
\end{align*}
\end{proof}
For future use, we also establish the resolvent estimate for $\varphi'$.
\begin{lemma}\label{resolvent estimate-4}For any $|k|\geq 1$, $\lambda\in\mathbb{R}$ and $w\in D_k$, there exists a constant $C>0$ independent of $\nu,k,B,\lambda, R$,  such that there holds
\begin{align*}
\nu^{\f{1}{2}}|kB|^{\f12}\|\varphi'\|_{L^2}+\nu^{\f{1}{2}}|k||kB|^{\f12}\|\f{\varphi}{r}\|_{L^2}\leq C R^2 \|F-c'\nu^{\f13}|kB|^{\f23}R^{-2}w\|_{H^{-1}_r}.
\end{align*}
\end{lemma}
\begin{proof}
If $\nu|k|^3> |B|$, by Lemma \ref{resolvent estimate-3} and Lemma \ref{Appendix A4} we obtain
\begin{equation*}
\begin{split}
     \|\varphi'\|_{L^2}+|k|\|\f{\varphi}{r}\|_{L^2}\lesssim& |k|^{-1} \|rw\|_{L^2}\lesssim  |k|^{-1} R^2\|\f{w}{r}\|_{L^2} \\\lesssim&  |k|^{-1} \nu^{-\f23} |kB|^{-\f13} R^2\|G\|_{H^{-1}_r}\\\le& \nu^{-\f12}|kB|^{-\f12} R^2\|G\|_{H^{-1}_r}.
\end{split}
\end{equation*}
Now we work under the case $\nu|k|^3\leq |B|$. Let $$E=\{r> 0: \  |1-\lambda r^2|\le \delta \} \ \textrm{and} \ E^c=(0, \infty)\backslash E$$ with $\delta=(\f{\nu}{|kB|})^{\f13}\ll 1$. Denote $a= \|\varphi'\|_{L^2}+|k|\|\f{\varphi}{r}\|_{L^2}$. Utilizing Lemma \ref{Appendix A4} from the Appendix, we immediately have
\begin{equation}\label{estimate for a}
    a^2 \lesssim |\langle w, \varphi\rangle |\le |\int_{E\cap [1, R]}w\overline{\varphi}  dr|+|\int_{E^c\cap [1, R]}w\overline{\varphi}  dr|.
\end{equation}
Employing Lemma \ref{resolvent estimate-3}, Lemma \ref{Appendix A1-1} and Lemma \ref{Appendix A5}, the first term can be controlled through
\begin{align*}
    |\int_{E\cap [1, R]}w\overline{\varphi}  dr| 
    \le& \|r^2 \f{w}{r^{\f12}}\|_{L^\infty}\|\f{\varphi}{r^{\f{5}{2}}}\|_{L^\infty} \|\f{1}{r}\|_{L^1(E)} \\
    \le& R^2 \|\f{w}{r^{\f12}}\|_{L^\infty}\|\f{\varphi}{r^{\f{1}{2}}}\|_{L^\infty} \|\f{1}{r}\|_{L^1(E)} \\
    \lesssim& \delta R^2(\|\f{w}{r}\|_{L^2} \|r^{\f12}(\f{w}{r^{\f12}})'\|_{L^2})^{\f12}  (\|\f{\varphi}{r}\|_{L^2} \|r^{\f12}(\f{\varphi}{r^{\f{1}{2}}})'\|_{L^2})^{\f12} \\
    \lesssim& \delta \nu^{-\f{5}{6}}|kB|^{-\f{1}{6}}R^2\|G\|_{H^{-1}_r} \times a.
\end{align*}
To bound the second term on the right of \eqref{estimate for a}, we use a piecewise $C^1$  cut-off function $\chi: (0, \infty) \to \mathbb{R}$ from Lemma \ref{Appendix A5}, which is defined as
\begin{align*}
\chi(r)=\left\{
\begin{aligned}
&\f{1}{\f{1}{r^2}-\lambda},\quad \text{if} \ r\in E^c,\\
&-\f{r_+^2+r_-^2}{\delta (r_+-r_-)}(r-r_-)+\f{r_-^2}{\delta},\quad  \text{if} \ r\in E
\end{aligned}
\right.
\end{align*}
with $r_{\pm}=\sqrt{\f{1\pm\delta}{\lambda}}$. Conducting integration by parts allows us to write
\begin{equation*}
\begin{split}
   \langle F, \chi \varphi\rangle 
   =&\nu \langle w', (\chi \varphi)'\rangle+\nu (k^2-\f{1}{4})\langle  \f{w}{r}, \f{\chi \varphi}{r} \rangle\\&+ikB(\int_{E\cap [1, R]}(\f{1}{r^2}-\lambda)\chi w\overline{\varphi} dr+\int_{E^c \cap [1, R]}(\f{1}{r^2}-\lambda)\chi w\overline{\varphi} dr ).
   \end{split}
\end{equation*}
Due to the fact that $\chi (\f{1}{r^2}-\lambda)=1$ on $E^c$, it can be inferred that
\begin{equation}\label{int w varphi r-4 over Ec}
\begin{split}
        |\int_{E^c\cap [1, R]}w\overline{\varphi} dr|
        \le& |kB|^{-1} \large(\|G\|_{H^{-1}_r}\|\chi \varphi\|_{H^1_r}+\nu \|w'\|_{L^2}\|(\chi \varphi)'\|_{L^2}\\&+\nu k^2 \|  \f{w}{r}\|_{L^2} \|\f{\chi \varphi}{r}\|_{L^2}\large)+\| \f{w}{r}\|_{L^2}\|(1-\lambda r^2)\f{\chi\varphi}{r}\|_{L^2(E\cap [1, R])}.
        \end{split}
\end{equation}
Applying Lemma \ref{Appendix A1-1} and Lemma \ref{Appendix A5}, we further have
\begin{align*}
    \|(\chi \varphi)'\|_{L^2}=&\|(r^2\f{\chi \varphi}{r^2})'\|_{L^2}\le \|r^{\f12}(\f{\chi}{r^2})'\|_{L^2}\|r^2\f{\varphi}{r^{\f{1}{2}}}\|_{L^\infty}+\|\f{\chi}{r^2}\|_{L^\infty}\|(r^2\varphi)'\|_{L^2} \\
    \lesssim& \delta^{-\f{3}{2}} R^2(\|\f{\varphi}{r}\|_{L^2} \|r^{\f12}(\f{\varphi}{r^{\f{1}{2}}})'\|_{L^2})^{\f12}+\delta^{-1} R^2(\|\varphi'\|_{L^2}+\|\f{\varphi}{r}\|_{L^2}) \lesssim \delta^{-\f32} R^2 a, \\
    \|\f{\chi\varphi}{r}\|_{L^2} \le& \|r^2 \f{\chi}{r^2}\|_{L^\infty}\|\f{\varphi}{r}\|_{L^2}\lesssim \delta^{-1}R^2a  
\end{align*}
and
\begin{align*}
    \|(1-\lambda r^2)\f{\chi\varphi}{r}\|_{L^2(E\cap [1, R])}\le& \|(1-\lambda r^2)\f{\chi}{r^{\f{5}{2}}}\|_{L^2(E)} \|r^2 \f{\varphi}{r^{\f{1}{2}}}\|_{L^\infty} \\
    \lesssim& \delta^{\f12}R^2 (\|\f{\varphi}{r}\|_{L^2} \|r^{\f12}(\f{\varphi}{r^{\f{1}{2}}})'\|_{L^2})^{\f12} \lesssim \delta^{\f12}R^2 a.
\end{align*}
Plugging the above inequalities into \eqref{int w varphi r-4 over Ec}, together with Lemma \ref{resolvent estimate-3}, we then deduce
\begin{align*}
     |\int_{E^c\cap [1, R]}w\overline{\varphi}  dr|
        \le& \Large(\delta^{-\f32}|kB|^{-1}+\delta^{-1} k^2 \nu^{\f13} |kB|^{-\f{4}{3}} +\delta^{\f12} \nu^{-\f23}|kB|^{-\f13}\Large)R^2 \|G\|_{H^{-1}_r}\times a.
\end{align*}
Therefore, the term $a= \|\varphi'\|_{L^2}+|k|\|\f{\varphi}{r}\|_{L^2}$ now obeys 
\begin{align*}
    a^2\lesssim \Large(\delta \nu^{-\f{5}{6}}|kB|^{-\f{1}{6}} +\delta^{-\f32}|kB|^{-1}+\delta^{-1} k^2 \nu^{\f13} |kB|^{-\f{4}{3}} +\delta^{\f12} \nu^{-\f23}|kB|^{-\f13} \Large)R^2 \|G\|_{H^{-1}_r} \times a.
\end{align*}
Substituting $\delta=(\f{\nu}{|kB|})^{\f13}$ and noticing $\nu|k|^3\leq |B|$, we thus arrive at
\begin{equation*}
    a\lesssim \nu^{-\f12}|kB|^{-\f12}\left(1+(\f{\nu |k|^3} {|B|})^{\f12} \right) R^2\|G\|_{H^{-1}_r}\lesssim 
    \nu^{-\f12}|kB|^{-\f12} R^2\|G\|_{H^{-1}_r}.
\end{equation*}
This completes the proof of this lemma.
\end{proof}

\section{\textbf{Enhanced dissipation and invisid damping}}\label{Enhanced dissipation and invisid damping}

\subsection{Pseudospectral bound} 
\begin{comment}
Recall in \cite{Pa} an operator $L$ in a Hilbert space $H$ is called accretive if
\begin{align*}
\Re\langle Lf, f\rangle\geq0,\quad \textrm{for any}\  f\in D(L).
\end{align*}
The operator $L$ is called m-accretive if in addition all $\Re\lambda<0$ belong to the resolvent set of $L$ (see \cite{Kato} for more details).  The \textbf{pseudospectral bound} of $L$ is defined to be
\begin{align}\label{pseudospectral bound}
\Psi(L)=\inf\{\|(L-i\lambda)f\|:f\in D(L),\lambda\in\mathbb{R},\|f\|=1\}.
\end{align}
\end{comment}
The direct implication of the resolvent estimates is to provide  controls for the semigroup by using pseudospectral bounds.

As in \cite{Pa}, we call an operator $L$ in a Hilbert space $H$ is accretive if
\begin{align*}
\Re\langle Lf, f\rangle\geq0,\quad \textrm{for any}\ f\in D(L).
\end{align*}
The operator $L$ is said to be m-accretive if, in addition, all $\Re\lambda<0$ belong to the resolvent set of $L$ (see \cite{Kato} for more details). The \textbf{pseudospectral bound} of $L$ is defined as
\begin{align}\label{pseudospectral bound}
\Psi(L)=\inf\{\|(L-i\lambda)f\|:f\in D(L),\lambda\in\mathbb{R},\|f\|=1\}.
\end{align}
Consider the operator $L_k$ 
\begin{align*}
L_kw:=&-\nu(\partial_r^2-\f{k^2-\f14}{r^2})w+i\f{kB}{r^2}w
\end{align*}
in the domain
\begin{align*}
D_k=\{w\in H_{loc}^2(\mathbb{R}_{+},dr)\cap L^2(\mathbb{R}_{+},dr):-\nu(\partial_r^2-\f{k^2-\f14}{r^2})w+i\f{kB}{r^2}w\in L^2(\mathbb{R}_{+},dr) \}.
\end{align*}
Note that $-\partial_r^2$ is an operator with the compact resolvent.  Since $L_k$ is a relatively compact perturbation of  $-\nu\partial_r^2$ in the domain $D_k$, it is hence clear that the operators $L_k$ also has the compact resolvent, this indicates that the operator $L_k$ has only point spectrum.

In Lemma \ref{trivial w' lemma} we have obtained
\begin{align*}
\Re\langle L_kw,w\rangle_{L^2}=\nu\|w'\|_{L^2}^2+\nu(k^2-\f14)\|\f{w}{r}\|_{L^2}^2\geq0.
\end{align*}
The above inequality indicates $L_k$ being accretive, and furthermore m-accretive. Recall that in Proposition \ref{resolvent estimate-main} we establish the following resolvent estimates
\begin{align*}
&\|r(L_k-i\lambda)w\|_{L^2}\geq C(\nu k^2)^{\f13}|B|^{\f23}\|\f{w}{r}\|_{L^2}.
\end{align*}
In view of the fact $r\in [1, R]$, this provides a lower bound for the pseudospectrum of $L_k$:
\begin{align*}
&\|(L_k-i\lambda)w\|_{L^2}\geq C(\nu k^2)^{\f13}|B|^{\f23}R^{-2}\|w\|_{L^2}.
\end{align*} 

We summarize it into
\begin{lemma}\label{peudospectral bound}Let $\Psi$ be defined as in (\ref{pseudospectral bound}). There exists some $C>0$ independent of $\nu,k, B, R$  such that
\begin{align*}
\Psi(L_k(L^2\rightarrow L^2))\geq C(\nu k^2)^{\f13}|B|^{\f23}R^{-2}.
\end{align*}
\end{lemma}

\subsection{Semigroup bound and enhanced dissipation} 
It is convenient to use the below space-time norm
\begin{align*}
	\|g\|_{L^pL^2}:=\Big\|\|g\|_{L^2_r([1, R])}\Big\|_{L^p_t([0, \infty)}.
	\end{align*}

To obtain decaying semigroup bounds from pseudospectral bounds, we appeal to the following Gearhart-Pr$\ddot{u}$ss type lemma established by Wei in \cite{Wei}. (See also \cite{HS} by Helffer \and 
Sj\"{o}strand for relevant reference.)

\begin{lemma}\cite{Wei} \label{GP lemma}Let $L$ be a m-accretive operator in a Hilbert space $X$. Then it holds
	\begin{align*}
	\|e^{-tL}\|_X\leq e^{-t\Psi(L)+\f{\pi}{2}} \quad \textrm{for any} \ t\geq0.
	\end{align*}
\end{lemma}

We proceed to studying the homogeneous linear equation
\begin{align}
\label{homogeneous linear equation}
&\partial_tw_k^l+L_kw_k^l=0 \quad \textrm{with} \ w_k^l(0)=w_k(0).
\end{align}
Utilizing semigroup theory, we can express $w_k^l$ as
\begin{align*}
w_k^l(t)=e^{-t L_k}w_k(0).
\end{align*}

We start to derive the space-time estimate for $w_k^l$.
\begin{proposition}\label{linear exp decay}Let $w_k^l$ be the solution to \eqref{homogeneous linear equation} with $w_k(0)\in L^2$. Then for any $k\in \mathbb{Z}$ and $|k|\geq1$, there exist constants $C,c>0$ being independent of $\nu,k,B,R$, such that the following inequality holds
	\begin{align}
	\label{linear exp decay-L2}&\|w_k^l(t)\|_{L^2}\leq Ce^{-c(\nu k^2)^{\f13}|B|^{\f23}R^{-2}t}\|w_k(0)\|_{L^2}, \quad \textrm{for any} \ t\geq0.
	\end{align}
	Moreover, for any $c'\in(0,c)$, we have
	\begin{align}
	\label{linear exp decay-L2L2}&(\nu k^2)^{\f13}|B|^{\f23}R^{-2}\|e^{c'(\nu k^2)^{\f13}|B|^{\f23}R^{-2}t}w_k^l(t)\|_{L^2L^2}^2\leq C\|w_k(0)\|_{L^2}^2.
	\end{align}
\end{proposition}
\begin{proof}
The semigroup bounds \eqref{linear exp decay-L2} readily follows from  Lemma \ref{peudospectral bound} and Lemma \ref{GP lemma}. Hence for any $c'\in(0,c)$, multiplying $e^{c'(\nu k^2)^{\f13}|B|^{\f23}R^{-2}t}$ on both sides of \eqref{linear exp decay-L2} yields
	\begin{align*} &2c'(\nu k^2)^{\f13}|B|^{\f23}R^{-2}\|e^{c'(\nu k^2)^{\f13}|B|^{\f23}R^{-2}t}w_k^l(t)\|_{L^2}^2\\
 \leq&2Cc'(\nu k^2)^{\f13}|B|^{\f23}R^{-2}e^{-2(c-c')(\nu k^2)^{\f13}|B|^{\f23}R^{-2}t}\|w_k(0)\|_{L^2}^2.
	\end{align*}
	Then we integrate above inequality with respect to $t$ to deduce
	\begin{align*} &2c'(\nu k^2)^{\f13}|B|^{\f23}R^{-2}\|e^{c'(\nu k^2)^{\f13}|B|^{\f23}R^{-2}t}w_k^l\|_{L^2L^2}^2\\
 \leq&\int_0^{\infty}2Cc'(\nu k^2)^{\f13}|B|^{\f23}R^{-2}e^{-2(c-c')(\nu k^2)^{\f13}|B|^{\f23}R^{-2}t}dt\|w_k(0)\|_{L^2}^2\lesssim\|w_k(0)\|_{L^2}^2.
	\end{align*}
This implies the estimate in \eqref{linear exp decay-L2L2}.
\end{proof}

\subsection{Homogeneous linear equation with nonzero initial data}

Based on the results of Proposition \ref{linear exp decay} and the structure of the linear equation itself, we can obtain the following $L^{\infty}L^2$ energy estimate for the homogeneous linear equation with respect to $w_k^l$.

\begin{lemma}\label{nonzero frequency linear part}Let $w_k^l$ be the solution to \eqref{homogeneous linear equation} with initial data $w_k(0)\in L^2$. Considering $c$ to be the same as in Proposition \ref{linear exp decay} and $c'\in(0,c)$ . For any $k\in \mathbb{Z}$ and $|k|\geq1$, it holds
	\begin{align}
	\label{nonzero linear part}&\|e^{c'(\nu k^2)^{\f13}|B|^{\f23}R^{-2}t}w_k^l\|_{L^{\infty}L^2}^2+(\nu k^2)^{\f13}|B|^{\f23}R^{-2}\|e^{c'(\nu k^2)^{\f13}|B|^{\f23}R^{-2}t}w_k^l\|_{L^2L^2}^2\\
\nonumber &+\nu\|e^{c'(\nu k^2)^{\f13}|B|^{\f23}R^{-2}t}\partial_{r}w_k^l\|_{L^2L^2}^2 +\nu k^2\|e^{c'(\nu k^2)^{\f13}|B|^{\f23}R^{-2}t}\f{w_k^l}{r}\|_{L^2L^2}^2\\
\nonumber \lesssim&\|w_k(0)\|_{L^2}^2.
	\end{align}
\end{lemma}
\begin{proof}
	We first conduct the integration by parts and get 
	\begin{align*}
	&\Re\langle\partial_tw_k^l-\nu (\partial_r^2-\f{k^2-\f14}{r^2})w_k^l+\f{ikB}{r^2}w_k^l,w_k^l\rangle\\ =&\f12\partial_t\|w_k^l\|_{L^2}^2+\nu\|\partial_rw_k^l\|_{L^2}^2+\nu(k^2-\f14)\|\f{w_k^l}{r}\|_{L^2}^2=0.
	\end{align*}
By multiplying $e^{2c'(\nu k^2)^{\f13}|B|^{\f23}R^{-2}t}$ on both sides, we deduce
	\begin{align*} &\partial_t\|e^{c'(\nu k^2)^{\f13}|B|^{\f23}R^{-2}t}w_k^l\|_{L^2}^2\\&+2\nu\Big(\|e^{c'(\nu k^2)^{\f13}|B|^{\f23}R^{-2}t}\partial_{r}w_k^l\|_{L^2}^2
 +(k^2-\f14)\|e^{c'(\nu k^2)^{\f13}|B|^{\f23}R^{-2}t}\f{w_k^l}{r}\|_{L^2}^2\Big)\\
	\nonumber\leq &2c'(\nu k^2)^{\f13}|B|^{\f23}R^{-2}\|e^{c'(\nu k^2)^{\f13}|B|^{\f23}R^{-2}t}w_k^l\|_{L^2}^2.
	\end{align*}
	This implies a space-time estimate for $w_k^l$:
	\begin{align*} &\|e^{c'(\nu k^2)^{\f13}|B|^{\f23}R^{-2}t}w_k^l\|_{L^{\infty}L^2}^2\\&+2\nu\Big(\|e^{c'(\nu k^2)^{\f13}|B|^{\f23}R^{-2}t}\partial_{r}w_k^l\|_{L^2L^2}^2
 +(k^2-\f14)\|e^{c'(\nu k^2)^{\f13}|B|^{\f23}R^{-2}t}\f{w_k^l}{r}\|_{L^2L^2}^2\Big)\\
	\leq & 2c'(\nu k^2)^{\f13}|B|^{\f23}R^{-2}\|e^{c'(\nu k^2)^{\f13}|B|^{\f23}R^{-2}t}w_k^l\|_{L^2L^2}^2+\|w_k(0)\|_{L^2}^2.
	\end{align*}
	Combining with Proposition \ref{linear exp decay}, we arrive at
	\begin{align*}	&\|e^{c'(\nu k^2)^{\f13}|B|^{\f23}R^{-2}t}w_k^l\|_{L^{\infty}L^2}^2+(\nu k^2)^{\f13}|B|^{\f23}R^{-2}\|e^{c'(\nu k^2)^{\f13}|B|^{\f23}R^{-2}t}w_k^l\|_{L^2L^2}^2\\
&+\nu\|e^{c'(\nu k^2)^{\f13}|B|^{\f23}R^{-2}t}\partial_{r}w_k^l\|_{L^2L^2}^2 +\nu(k^2-\f14)\|e^{c'(\nu k^2)^{\f13}|B|^{\f23}R^{-2}t}\f{w_k^l}{r}\|_{L^2L^2}^2\lesssim\|w_k(0)\|_{L^2}^2.
	\end{align*}
	This completes the proof of Lemma \ref{nonzero frequency linear part}.
\end{proof}

\subsection{Inhomogeneous linear equation with zero initial data} We then derive the space-time estimates for inhomogeneous equations with zero initial conditions, which will be frequently used in later sections. The derivation of this lemma  heavily relies on established resolvent estimates in Proposition \ref{resolvent estimate-main} from Section \ref{Section Resolvent estimate}.

\begin{proposition}\label{spacetime estimate proposition-nonlinear}
Given functions $h_1(t, r), h_2(t, r)$ and $ g(r)$, assume $w$ solves
\begin{align}\label{spacetime estimate lemma eqn}
\left\{
\begin{aligned}
&\partial_t w+L_k w=h_1-g \partial_{r}h_2,\\
&w|_{t=0}=0,\quad w|_{r=1,R}=0,
\end{aligned}
\right.
\end{align}
and let the stream function $\varphi$  satisfy
\begin{align*}
&(\partial_r^2-\f{k^2-\f14}{r^2})\varphi=w,\quad \varphi|_{r=1,R}=0 
\quad \textrm{with}  \ r\in[1,R] \ \textrm{and}  \ t\geq0.
\end{align*}
With $c$ being the same constant as in Proposition \ref{linear exp decay}, then for any $c'\in [0, c]$, the following inequalities hold for $w$ and $\varphi$:
\begin{equation*}
    \begin{split}
&\|e^{c'(\nu k^2)^{\f13}|B|^{\f23}R^{-2}t}w\|_{L^{\infty}L^2}^2+\nu^{\f13} |kB|^{\f23}\|e^{c'(\nu k^2)^{\f13}|B|^{\f23}R^{-2}t}\f{w}{r}\|_{L^2L^2}^2\\
&+\nu\|e^{c'(\nu k^2)^{\f13}|B|^{\f23}R^{-2}t}\partial_{r}w\|_{L^2L^2}^2+\nu k^2\|e^{c'(\nu k^2)^{\f13}|B|^{\f23}R^{-2}t}\f{w}{r}\|_{L^2L^2}^2\\
\lesssim&\nu^{-\f13} |kB|^{-\f23}\|e^{c'(\nu k^2)^{\f13}|B|^{\f23}R^{-2}t}rh_1\|_{L^2L^2}^2+\nu^{-1} \|e^{c'(\nu k^2)^{\f13}|B|^{\f23}R^{-2}t}(|g|+r|\partial_r g|)h_2\|_{L^2L^2}^2,
    \end{split}
\end{equation*}
and
\begin{align*}
&\nu^{\f16}|kB|^{\f56}|k|^{\f12}R^{-2}(\|e^{c'(\nu k^2)^{\f13}|B|^{\f23}R^{-2}t}\partial_r{\varphi}\|_{L^2L^2}+|k|\|e^{c'(\nu k^2)^{\f13}|B|^{\f23}R^{-2}t}\f{{\varphi}}{r}\|_{L^2L^2}) \\
    \lesssim& \left((\f{\nu}{|kB|})^{\f16} (\log R)^{\f12} +1 \right) \|e^{c'(\nu k^2)^{\f13}|B|^{\f23}R^{-2}t}rh_1\|_{L^2L^2}\\&+\nu^{-\f13} |kB|^{\f13}\|e^{c'(\nu k^2)^{\f13}|B|^{\f23}R^{-2}t}(|g|+r|\partial_r g|)h_2\|_{L^2L^2}.
\end{align*}
\end{proposition}
\begin{proof}
We introduce the weighted quantities
$$\tilde{w}:=e^{c'(\nu k^2)^{\f13}|B|^{\f23}R^{-2}t}w_k, \quad \tilde{\varphi}:=e^{c'(\nu k^2)^{\f13}|B|^{\f23}R^{-2}t}\varphi$$ and  
$$\tilde{h}_j=e^{c'(\nu k^2)^{\f13}|B|^{\f23}R^{-2}t}h_j \quad \textrm{for} \ j=1,2.$$ Via a direct check, we can see that
\begin{equation}\label{spacetime estimate lemma eqn tilde w}
    \partial_t \tilde{w}+(L_k-c'(\nu k^2)^{\f13}|B|^{\f23}R^{-2}) \tilde{w}=\f{1}{r}[ik\tilde{h}_1-r^{\f12}\partial_{r}(r^{\f12}\tilde{h}_2)].
\end{equation}
	We then take the Fourier transform in $t$ and define\footnote{Note that $\varphi|_{r=1, R}=0$ in \eqref{spacetime estimate lemma eqn} guarantees the Fourier transforms here being well-defined for $\lambda\in (0,\infty)$.}
	\begin{align*} &\hat{w}(\lambda,r)=\int_0^{\infty}\tilde{w}(t,r)e^{-it\lambda}dt,\quad \hat{\varphi}(\lambda,r)=\int_0^{\infty}\tilde{\varphi}(t,r)e^{-it\lambda}dt,  \\&H_j(\lambda,r)=\int_0^{\infty}\tilde{h}_j(t,r)e^{-it\lambda}dt \quad \textrm{for} \ j=1,2.
	\end{align*}
The inhomogeneous equation \eqref{spacetime estimate lemma eqn tilde w} can thus be transferred into the form
	\begin{align*} (i\lambda+L_k-c'(\nu k^2)^{\f13}|B|^{\f23}R^{-2})\hat{w}(\lambda,r)=H_1-g\partial_{r}H_2.
	\end{align*}
	We further decompose $\hat{w}$ as $\hat{w}=\hat{w}^{1}+\hat{w}^{2}$, where $\hat{w}^{1}$ and $\hat{w}^{2}$ solve
	\begin{align*}
	\Big(i\lambda+L_k-c'(\nu k^2)^{\f13}|B|^{\f23}R^{-2}\Big)\hat{w}^{1}=H_1,
 	\end{align*}
and
\begin{align*}
	\Big(i\lambda+L_k-c'(\nu k^2)^{\f13}|B|^{\f23}R^{-2}\Big)\hat{w}^{2}=-g\partial_{r}H_2.
	\end{align*}
	We also define the corresponding stream functions $\hat{\varphi}^1$ and  $\hat{\varphi}^2$ via the below linear elliptic equations
	\begin{equation*}
	    (\partial_r^2-\f{k^2-\f14}{r^2})\hat{\varphi}^j=\hat{w}^j,\quad \hat{\varphi}^j|_{r=1,R}=0,\quad r\in[1,R], \quad \textrm{for} \ j=1,2.
	\end{equation*}
 This enables us to write $\hat{\varphi}=\hat{\varphi}^1+\hat{\varphi}^2$. In view of Proposition \ref{resolvent estimate-main}, we then obtain the estimate for $w^1$ and $\varphi^1$:
	\begin{align*}
&\nu^{\f23} |kB|^{\f13}\|\partial_{r}\hat{w}^{1}\|_{L^2}+\nu^{\f13} |kB|^{\f23}\|\f{\hat{w}^{1}}{r}\|_{L^2}\leq C \|rH_1\|_{L^2},
	\end{align*}
	\begin{align*}
	   \nu^{\f16}|kB|^{\f56}|k|^{\f12}\left( \|\partial_r \hat{\varphi}^1\|_{L^2}+|k|\|\f{\hat{\varphi}^1}{r}\|_{L^2} \right) \leq CR^2\left((\f{\nu}{|kB|})^{\f16} (\log R)^{\f12} +1 \right) \|rH_1\|_{L^2},
	\end{align*}
	as well as the control for $w^2$ and $\varphi^2$:
	\begin{align*} &\nu\|\partial_{r}\hat{w}^{2}\|_{L^2}+\nu^{\f23} |kB|^{\f13}\|\f{\hat{w}^{2}}{r}\|_{L^2}\leq  C\|g\partial_r H_2\|_{H^{-1}_r}\lesssim\|(|g|+r|g'|)H_2\|_{L^2},
	\end{align*}
	\begin{align*}
	    \nu^{\f12}|kB|^{\f12}R^{-2}(\|\partial_r\hat{\varphi}^{ 2}\|_{L^2}+|k|\|\f{\hat{\varphi}^{ 2}}{r}\|_{L^2})\leq  C\|g\partial_r H_2\|_{H^{-1}_r}\lesssim  \|(|g|+r|g'|)H_2\|_{L^2}.
	\end{align*}
 Here  we utilize the definition of $H^{-1}_r$ norm to bound 
 \begin{equation*}
     \begin{split}
         \|g\partial_r H_2\|_{H^{-1}}:=&\sup\limits_{\|f\|_{H^1_r}\le 1} |\langle g\partial_r H_2, f\rangle| =\sup\limits_{\|f\|_{H^1_r}\le 1} |\langle  H_2, (gf)' \rangle| \\
         \le&  \sup\limits_{\|f\|_{H^1_r}\le 1} |\langle  gH_2, f' \rangle|+\sup\limits_{\|f\|_{H^1_r}\le 1} |\langle  rg' H_2, \f{f}{r} \rangle| \\
         \lesssim& \|(|g|+r|g'|)H_2\|_{L^2}.
     \end{split}
 \end{equation*}
	Combining the above inequalities, we arrive at
	\begin{align*}
	\nu^{\f23} |kB|^{\f13}\|\partial_{r}\hat{w}\|_{L^2}+\nu^{\f13} |kB|^{\f23}\|\f{\hat{w}}{r}\|_{L^2}
	\lesssim \|rH_1\|_{L^2}+\nu^{-\f13} |kB|^{\f13}\|(|g|+r|g'|)H_2\|_{L^2},
	\end{align*}
	and
	\begin{align*}
	    &\nu^{\f16}|kB|^{\f56}R^{-2}(\|\partial_r\hat{\varphi}\|_{L^2}+|k|\|\f{\hat{\varphi}}{r}\|_{L^2})\\\lesssim&
	    \left((\f{\nu}{|kB|})^{\f16} (\log R)^{\f12} +1 \right) \|rH_1\|_{L^2}+\nu^{-\f13} |kB|^{\f13}\|(|g|+r|g'|)H_2\|_{L^2}.
	\end{align*}
	According to the Plancherel's theorem, we have the following equivalence relations based on $L^2$ norms
	\begin{align*}
&\|\partial_{r}\tilde{w}\|_{L^2L^2}\approx\big\|\|\partial_{r}\hat{w}\|_{L^2}\big\|_{L^2(\mathbb{R})}, \quad \|\f{\tilde{w}}{r}\|_{L^2L^2}\approx\big\|\|\f{\hat{w}}{r}\|_{L^2}\big\|_{L^2(\mathbb{R})}, \\
&\|\partial_r\tilde{\varphi}\|_{L^2L^2}+|k|\|\f{\tilde{\varphi}}{r}\|_{L^2L^2}\approx \big\|\|\partial_r\hat{\varphi}\|_{L^2}\big\|_{L^2(\mathbb{R})}+|k|\big\|\|\f{\hat{\varphi}}{r}\|_{L^2}\big\|_{L^2(\mathbb{R})},\\
&\|r\tilde{h}_1\|_{L^2L^2}\approx\big\|\|rH_1\|_{L^2}\big\|_{L^2(\mathbb{R})}, \quad \|(|g|+r|g'|)\tilde{h}_2\|_{L^2L^2}\approx\big\|\|(|g|+r|g'|)H_2\|_{L^2}\big\|_{L^2(\mathbb{R})}.
	\end{align*}
Plugging in all estimates above, we thus deduce
	\begin{align}
\label{spacetime estimate lemma w L2L2}&\nu^{\f23} |kB|^{\f13}\|e^{c'(\nu k^2)^{\f13}|B|^{\f23}R^{-2}t}\partial_{r}w_k\|_{L^2L^2}+\nu^{\f13} |kB|^{\f23}\|e^{c'(\nu k^2)^{\f13}|B|^{\f23}R^{-2}t}\f{w_k}{r}\|_{L^2L^2}	\\
\nonumber\lesssim&\|e^{c'(\nu k^2)^{\f13}|B|^{\f23}R^{-2}t}rh_1\|_{L^2L^2}+\nu^{-\f13} |kB|^{\f13}\|e^{c'(\nu k^2)^{\f13}|B|^{\f23}R^{-2}t}(|g|+r|g'|)h_2\|_{L^2L^2},
	\end{align}
and
\begin{align*}
    &\nu^{\f16}|kB|^{\f56}R^{-2}(\|e^{c'(\nu k^2)^{\f13}|B|^{\f23}R^{-2}t}\partial_r{\varphi}_k\|_{L^2L^2}+|k|\|e^{c'(\nu k^2)^{\f13}|B|^{\f23}R^{-2}t}\f{{\varphi}_k}{r}\|_{L^2L^2}) \\
    \lesssim& \left((\f{\nu}{|kB|})^{\f16} (\log R)^{\f12} +1 \right) \|e^{c'(\nu k^2)^{\f13}|B|^{\f23}R^{-2}t}kh_1\|_{L^2L^2}+\nu^{-\f13} |kB|^{\f13}\|e^{c'(\nu k^2)^{\f13}|B|^{\f23}R^{-2}t}h_2\|_{L^2L^2}.
\end{align*}
	Applying the integration by parts, we also get
	\begin{align*}
	0=&\Re\langle\partial_tw- \nu(\partial_{r}^2-\f{k^2-\f14}{r^2})w+\f{ikB}{r^2}w-h_1+g\partial_{r}h_2,e^{2c'(\nu k^2)^{\f13}|B|^{\f23}R^{-2}t}w\rangle\\
=&\f12\partial_t\|e^{c'(\nu k^2)^{\f13}|B|^{\f23}R^{-2}t}w\|_{L^2}^2+\nu\|e^{c'(\nu k^2)^{\f13}|B|^{\f23}R^{-2}t}\partial_{r}w\|_{L^2}^2+\nu(k^2-\f14)\|e^{c'(\nu k^2)^{\f13}|B|^{\f23}R^{-2}t}\f{w}{r}\|_{L^2}^2\\
&-c'(\nu k^2)^{\f13}|B|^{\f23}R^{-2}\|e^{c'(\nu k^2)^{\f13}|B|^{\f23}R^{-2}t}w\|_{L^2}^2-\Re\langle h_1-g\partial_{r}h_2,e^{2c'(\nu k^2)^{\f13}|B|^{\f23}R^{-2}t}w\rangle,
	\end{align*}
\begin{comment}
which implies
	\begin{align*}
&c'(\nu k^2)^{\f13}|B|^{\f23}R^{-2}\|e^{c'(\nu k^2)^{\f13}|B|^{\f23}R^{-2}t}w\|_{L^2}^2\\
=&\f12\partial_t\|e^{c'(\nu k^2)^{\f13}|B|^{\f23}R^{-2}t}w\|_{L^2}^2+\nu\|e^{c'(\nu k^2)^{\f13}|B|^{\f23}R^{-2}t}\partial_{r}w\|_{L^2}^2+\nu(k^2-\f14)\|e^{c'(\nu k^2)^{\f13}|B|^{\f23}R^{-2}t}\f{w}{r}\|_{L^2}^2\\ &-\Re\langle h_1-g\partial_{r}h_2,e^{2c'(\nu k^2)^{\f13}|B|^{\f23}R^{-2}t}w\rangle.
	\end{align*}
\end{comment}
	By employing the Cauchy-Schwarz inequality, we further obtain
	\begin{align*}
&\partial_t\|e^{c'(\nu k^2)^{\f13}|B|^{\f23}R^{-2}t}w\|_{L^2}^2+2\nu\|e^{c'(\nu k^2)^{\f13}|B|^{\f23}R^{-2}t}\partial_{r}w\|_{L^2}^2+2\nu(k^2-\f14)\|e^{c'(\nu k^2)^{\f13}|B|^{\f23}R^{-2}t}\f{w}{r}\|_{L^2}^2\\
=&2c'(\nu k^2)^{\f13}|B|^{\f23}R^{-2}\|e^{c'(\nu k^2)^{\f13}|B|^{\f23}R^{-2}t}w\|_{L^2}^2+2\Re \langle h_1-g\partial_{r}h_2,e^{2c'(\nu k^2)^{\f13}|B|^{\f23}R^{-2}t}w\rangle\\
\lesssim&c'(\nu k^2)^{\f13}|B|^{\f23}R^{-2}\|e^{c'(\nu k^2)^{\f13}|B|^{\f23}R^{-2}t}w\|_{L^2}^2+\|e^{c'(\nu k^2)^{\f13}|B|^{\f23}R^{-2}t}rh_1\|_{L^2}\|e^{c'(\nu k^2)^{\f13}|B|^{\f23}R^{-2}t}\f{w}{r}\|_{L^2}\\
&+\|e^{c'(\nu k^2)^{\f13}|B|^{\f23}R^{-2}t}(|g|+r|g'|)h_2\|_{L^2}\|e^{c'(\nu k^2)^{\f13}|B|^{\f23}R^{-2}t}w\|_{H^1_r}.
	\end{align*}
	Together with \eqref{spacetime estimate lemma w L2L2}, the above inequality yields
	\begin{align*}
&\|e^{c'(\nu k^2)^{\f13}|B|^{\f23}R^{-2}t}w\|_{L^{\infty}L^2}^2+\nu\|e^{c'(\nu k^2)^{\f13}|B|^{\f23}R^{-2}t}\partial_{r}w\|_{L^2L^2}^2+\nu k^2\|e^{c'(\nu k^2)^{\f13}|B|^{\f23}R^{-2}t}\f{w}{r}\|_{L^2L^2}^2\\
&+\nu^{\f13} |kB|^{\f23}\|e^{c'(\nu k^2)^{\f13}|B|^{\f23}R^{-2}t}\f{w}{r}\|_{L^2L^2}^2\\
\lesssim&\nu^{-\f13} |kB|^{-\f23}\Big(\|e^{c'(\nu k^2)^{\f13}|B|^{\f23}R^{-2}t}rh_1\|_{L^2L^2}+\nu^{-\f13} |kB|^{\f13}\|e^{c'(\nu k^2)^{\f13}|B|^{\f23}R^{-2}t}(|g|+r|g'|)h_2\|_{L^2L^2}\Big)^2\\
\lesssim&\nu^{-\f13} |kB|^{-\f23}\|e^{c'(\nu k^2)^{\f13}|B|^{\f23}R^{-2}t}rh_1\|_{L^2L^2}^2+\nu^{-1} \|e^{c'(\nu k^2)^{\f13}|B|^{\f23}R^{-2}t}(|g|+r|g'|)h_2\|_{L^2L^2}^2.
	\end{align*}
This finishes the proof of Proposition \ref{spacetime estimate proposition-nonlinear}.
\end{proof}

\subsection{Integrated invisid damping}  This subsection is devoted to providing estimates for the stream function $\varphi_k^l$ in terms of initial vorticity $w_k(0)$, with $\varphi_k^l$ admitting the linear elliptic equation
\begin{equation}\label{def linear phi k l} 
	    (\partial_r^2-\f{k^2-\f14}{r^2})\varphi_k^l=w_k^l,\quad \varphi_k^l|_{r=1,R}=0,\quad r\in[1,R].
	\end{equation}
 For inviscid fluids governed by the Euler equation, there is the concept of so-called inviscid damping. In fact, the effect of inviscid damping extends to viscous fluids as well. Lemma \ref{nonzero frequency linear part for varphi} presents a space-time version of linear inviscid damping around 2D TC flow.  

\begin{lemma}\label{nonzero frequency linear part for varphi}
Let $w_k^l$ be the solution to \eqref{homogeneous linear equation} with initial data $w_k(0)\in L^2$ and $\varphi_k^l$ solves \eqref{def linear phi k l}. Considering $c$ to be the same as in Proposition \ref{linear exp decay} and $c'\in(0,c)$. For any $k\in \mathbb{Z}$ and $|k|\geq1$, the following inequality holds 
\begin{equation}\label{nonzero frequency linear part for varphi sum}
\begin{split}
     &k^2|B|R^{-4}(\|e^{c'(\nu k^2)^{\f13}|B|^{\f23}R^{-2}t}\partial_r{\varphi}^{l}_k\|_{L^2L^2}^2+|k|^2\|e^{c'(\nu k^2)^{\f13}|B|^{\f23}R^{-2}t}\f{{\varphi}^{l}_k}{r}\|_{L^2L^2}^2) \\
     \lesssim& (\log R)^{-2}R^{-4}\|r^2 w_k(0)\|_{L^2}^2+R^6 \|\f{w_k(0)}{r^3}\|_{L^2}^2+(\f{\nu}{|kB|})^{\f23}R^2\|\partial_r w_k(0)\|_{L^2}^2 \\&+\left((\f{\nu}{|kB|})^{\f13} \log R +1 \right)\left( R^{-2}\|rw_k(0)\|_{L^2}^2 + (\f{\nu}{|kB|})^{\f43}k^4 R^2\|\f{w_k(0)}{r}\|_{L^2}^2\right).
\end{split}
\end{equation}
\end{lemma}

\begin{remark}
Note that $\nu$ is much smaller than $1$, and $B$ is a fixed constant, so \eqref{nonzero frequency linear part for varphi sum} can also be written in the following expression
\begin{align*}
&|k||B|^{\f12}(\|e^{c'(\nu k^2)^{\f13}|B|^{\f23}R^{-2}t}\partial_r{\varphi}^{l}_k\|_{L^2L^2}+|k|\|e^{c'(\nu k^2)^{\f13}|B|^{\f23}R^{-2}t}\f{{\varphi}^{l}_k}{r}\|_{L^2L^2})\\
\lesssim& C(R) (\|\partial_r w_k(0)\|_{L^2}+\|w_k(0)\|_{L^2}).
\end{align*}
We can see that both sides of the above inequality are independent of viscosity coefficient $\nu$, which indicates that the result is also valid for the Euler equation. This type of estimate is called inviscid damping, and here we establish its integrated form. Thus we refer to it as the integrated invisid damping estimates.
\end{remark}

In order to prove \eqref{nonzero frequency linear part for varphi sum}, we perform a suitable decomposition of $\varphi^l_k$. Define
\begin{equation}
\begin{split}
    &\tilde{w}^{l_0}_k(t,r):=e^{-ikB\f{t}{r^2}} w_k(0), \quad \tilde{w}^{l_1}_k(t, r):=e^{-c(\nu k^2)^{\f13}|B|^{\f23}R^{-2} t} \tilde{w}^{l_0}_k
\end{split}
\end{equation}
and denote $\tilde{w}^{l_2}_k$ to be the solution to the inhomogeneous linear equation equation with zero initial conditions as below: 
\begin{equation}
    \partial_t \tilde{w}^{l_2}_k+L_k \tilde{w}^{l_2}_k=\nu\partial_r^2 \tilde{w}^{l_1}_k-\nu(k^2-\f{1}{4})\f{\tilde{w}^{l_1}_k}{r^2}+c(\nu k^2)^{\f13}|B|^{\f23}R^{-2} \tilde{w}^{l_1}_k, \quad \tilde{w}^{l_2}_k(0)=0.
\end{equation}
As a result, the corresponding stream functions can be defined by
\begin{equation}\label{tilde varphi_k^lj}
     (\partial_r^2-\f{k^2-\f14}{r^2})\tilde{\varphi}^{l_j}_k=w_k^{l_j},\quad \tilde{\varphi}^{l_j}_k|_{r=1,R}=0,\quad r\in[1,R],\quad t\geq0, \quad j=1,2.
\end{equation}
Then it is easily seen that $w^l_k=\tilde{w}^{l_1}_k+\tilde{w}^{l_2}_k$ and $\varphi^l_k=\tilde{\varphi}^{l_1}_k+\tilde{\varphi}^{l_2}_k$. 

We start with the estimates for $\tilde{\varphi}^{l_1}_k$.
\begin{lemma}\label{nonzero frequency linear part for varphi l1 part}
Let $w_k^l$ be the solution to \eqref{homogeneous linear equation} with initial data $w_k(0)\in L^2$. Considering $c$ to be the same as in Proposition \ref{linear exp decay}, $c'\in(0,c)$ and $\tilde{\varphi}^{l_1}_k$ defined in \eqref{tilde varphi_k^lj}. For any $k\in \mathbb{Z}$ and $|k|\geq1$, the following inequality holds 
\begin{equation}
     |kB||k|(\log R)^2 \left( \|e^{c'(\nu k^2)^{\f13}|B|^{\f23}R^{-2} t}\partial_r \tilde{\varphi}^{l_1}_k\|_{L^2L^2}^2+k^2\|e^{c'(\nu k^2)^{\f13}|B|^{\f23}R^{-2} t}\f{\tilde{\varphi}^{l_1}_k}{r}\|^2_{L^2L^2}\right)\lesssim \|r^{2}w_k(0)\|_{L^2}^2
\end{equation}
\end{lemma}
\begin{proof}
Noting that $\tilde{\varphi}^{l_1}_k=e^{-c(\nu k^2)^{\f13}|B|^{\f23}R^{-2} t}\tilde{\varphi}^{l_0}_k$ and $c'\in(0,c)$, to obtain the above conclusion it suffices to prove 
\begin{equation}
     |kB||k|(\log R)^2\left( \|\partial_r \tilde{\varphi}^{l_0}_k\|_{L^2L^2}^2+k^2\|\f{\tilde{\varphi}^{l_0}_k}{r}\|^2_{L^2L^2}\right)\lesssim \|r^{2}w_k(0)\|_{L^2}^2.
\end{equation}
Here we develop some new ideas. With $w=1/r^2$, we introduce the weighted Hilbert space $L^2_w([1, R])$ with inner product 
\begin{equation*}
    \langle f, g\rangle_w =\int_1^R f\overline{g} w dr.
\end{equation*}
We find a set of explicit orthonormal basis for $L^2_w([1, R])$ as below:
\begin{equation*}
    \psi_l(r)=(\f{2}{\log R})^{\f12}r^{1/2}\sin (\f{l\pi}{\log R} \log r), \quad \textrm{for all} \ l\in \mathbb{N}^+.
\end{equation*}
Here $\psi_l$ satisfies the following property
\begin{equation}
    \left(\partial_r^2-\f{k^2-\f{1}{4}}{r^2}\right)\psi_l(r)=-\left((\f{l\pi}{\log R})^2+k^2\right) w(r) \psi_l(r), \quad \psi_l|_{r=1, R}=0.
\end{equation}
As far as we know, our explicit constructions above are new. For the sake of convenience, we denote $\alpha:=(\f{2}{\log R})^{\f12}$, $\beta:=\f{\pi}{\log R}$ and $\lambda_{k,l}:=(\beta l)^2+k^2$. Since
\begin{equation*}
   \begin{split}
        \langle \tilde{\varphi}^{l_0}_k, \psi_l \rangle_w=&  \langle \tilde{\varphi}^{l_0}_k, w\psi_l \rangle = -\lambda_{k ,l}^{-1}\langle \tilde{\varphi}^{l_0}_k,  \left(\partial_r^2-\f{k^2-\f{1}{4}}{r^2}\right)\psi_l \rangle 
        =-\lambda_{k ,l}^{-1}\langle \tilde{w}^{l_0}_k,  \psi_l \rangle,
   \end{split}
\end{equation*}
we can write
\begin{equation}\label{inner product of w l0 and phi l0}
   \begin{split}
      & -\langle \tilde{w}^{l_0}_k, \tilde{\varphi}^{l_0}_k \rangle=
    -\langle w^{-1}\tilde{w}^{l_0}_k, \tilde{\varphi}^{l_0}_k \rangle_w\\=&\sum\limits_{l=1}^{\infty} \lambda_{k, l}^{-1} \langle w^{-1}\tilde{w}^{l_0}_k, \psi_l \rangle_w \overline{\langle \tilde{w}^{l_0}_k , \psi_l \rangle} 
    =\sum\limits_{l=1}^{\infty} \lambda_{k, l}^{-1} |\langle \tilde{w}^{l_0}_k , \psi_l \rangle|^2.
    \end{split}
\end{equation}
In view of the definition $\tilde{w}^{l_0}_k=e^{-i kB\f{t}{r^2}}w_k(0)$, we have the expression
\begin{equation*}
\begin{split}
    \langle \tilde{w}^{l_0}_k, \psi_l \rangle=&\int_{1}^{R} e^{-ikB\f{t}{r^2}}w_k(0, r)\psi_l(r) dr\\=&\int_{\f{1}{R^2}}^1 e^{-i kBt s} w_k(0, \f{1}{\sqrt{s}}) \psi_l (\f{1}{\sqrt{s}})\f{ds}{2s\sqrt{s}} 
    :=\mathcal{F}[W_{k, l}](kBt),
    \end{split}
\end{equation*}
where $\mathcal{F}$ represents the canonical Fourier transform over $\mathbb{R}$, and $W^1_{k, l}(r)=\f{1}{2s\sqrt{s}}w_k(0, \f{1}{\sqrt{s}}) \psi_l (\f{1}{\sqrt{s}})$.\footnote{According to the boundary condition of $w^l_k$, the Fourier transform of $W^1_{k, l}(r)$ is well-defined.}

\noindent Note that the Plancherel’s formula implies
\begin{equation*}
    \begin{split}
    \int_{\mathbb{R}} |\langle \tilde{w}^{l_0}_k , \psi_l \rangle|^2dt=&\f{2\pi}{|kB|} \int_{\f{1}{R^2}}^1 |W_{k, l}(r)|^2 dr \\
    =&\f{2\pi}{|kB|}\int_{1}^R 2r^3 |w_k(0 ,r)|^2 \psi_l^2(r) dr\le \f{4\pi \alpha^2}{|kB|} \|r^{2}w_k(0)\|_{L^2}^2.
        \end{split}
\end{equation*}
This together with \eqref{inner product of w l0 and phi l0} yields
\begin{align*}
    \int_{\mathbb{R}} -\langle \tilde{w}^{l_0}_k, \tilde{\varphi}^{l_0}_k  \rangle dt\leq&\f{4\alpha^2\pi}{|kB|} \|r^{2}w_k(0)\|_{L^2}^2\sum\limits_{l=1}^{\infty} \f{1}{(\beta l)^2+k^2} \\
    \le& \f{4\alpha^2 \pi}{|kB|} \|r^{2}w_k(0)\|_{L^2}^2\int_{0}^{+\infty} \f{dy}{(\beta y)^2+k^2} = \f{2\alpha^2 \beta \pi^2}{|kB|\cdot |k|} \|r^{2}w_k(0)\|_{L^2}^2.
\end{align*}
On the other hand, observing
\begin{equation*}
    \begin{split}
        -\langle \tilde{w}^{l_0}_k, \tilde{\varphi}^{l_0}_k  \rangle
        =\|\partial_r \tilde{\varphi}^{l_0}_k\|_{L^2}^2+(k^2-\f{1}{4})\|\f{\tilde{\varphi}^{l_0}_k}{r}\|_{L^2}^2\ge 0 \quad \textrm{for any} \ t\in \mathbb{R},
    \end{split}
\end{equation*} 
we infer that
\begin{equation*}
    \begin{split}
   &\|\partial_r \tilde{\varphi}^{l_0}_k\|_{L^2L^2}^2+(k^2-\f{1}{4})\|\f{\tilde{\varphi}^{l_0}_k}{r}\|_{L^2L^2}^2 \le\f{2\alpha^2 \beta \pi^2}{|kB|\cdot |k|} \|r^{2}w_k(0)\|_{L^2}^2\\
   \lesssim& |kB|^{-1}|k|^{-1}(\log R)^{-2} \|r^{2}w_k(0)\|_{L^2}^2.
    \end{split}
\end{equation*}
This completes the proof of this lemma.
\end{proof}

To control the $L^2L^2$ norm of $\tilde{\varphi}^{l_2}_k$, we utilize resolvent estimates in Section \ref{Section Resolvent estimate} together with Lemma \ref{nonzero frequency linear part}. 

\begin{lemma}\label{nonzero frequency linear part for varphi l2 part}
Let $w_k^l$ be the solution to \eqref{homogeneous linear equation} with initial data $w_k(0)\in L^2$. Considering $c$ to be the same as in Proposition \ref{linear exp decay} and $c'\in(0,c)$. For any $k\in \mathbb{Z}$ and $|k|\geq1$, there holds 
\begin{align*}
     &k^2|B|R^{-4}(\|e^{c'(\nu k^2)^{\f13}|B|^{\f23}R^{-2}t}\partial_r{\varphi}^{l_2}_k\|_{L^2L^2}^2+|k|^2\|e^{c'(\nu k^2)^{\f13}|B|^{\f23}R^{-2}t}\f{{\varphi}^{l_2}_k}{r}\|_{L^2L^2}^2) \\
     \lesssim& \left((\f{\nu}{|kB|})^{\f13} \log R +1 \right)\left( R^{-2}\|rw_k(0)\|_{L^2}^2 + (\f{\nu}{|kB|})^{\f43}k^4 R^2\|\f{w_k(0)}{r}\|_{L^2}^2\right)\\
     &+R^6 \|\f{w_k(0)}{r^3}\|_{L^2}^2+(\f{\nu}{|kB|})^{\f23}R^2\|\partial_r w_k(0)\|_{L^2}^2.
\end{align*}
\end{lemma}
\begin{proof}
First recall that $\tilde{w}^{l_2}_k$ satisfies the inhomogeneous linear equation
\begin{equation}
    \partial_t \tilde{w}^{l_2}_k+L_k \tilde{w}^{l_2}_k=\nu\partial_r^2 \tilde{w}^{l_1}_k-\nu(k^2-\f{1}{4})\f{\tilde{w}^{l_1}_k}{r^2}+c'(\nu k^2)^{\f13}|B|^{\f23}R^{-2} \tilde{w}^{l_1}_k, \quad \tilde{w}^{l_2}_k(0)=0,
\end{equation}
with the stream function $\tilde{\varphi}^{l_2}_k$ solving the boundary value problem
\begin{equation*}
     (\partial_r^2-\f{k^2-\f14}{r^2})\tilde{\varphi}^{l_2}_k=w_k,\quad \tilde{\varphi}^{l_2}_k|_{r=1,R}=0,\quad r\in[1,R],\quad t\geq0.
\end{equation*}
Employing Proposition \ref{spacetime estimate proposition-nonlinear} with $h_1=-\nu(k^2-\f{1}{4})\f{\tilde{w}^{l_1}_k}{r^2}+c(\nu k^2)^{\f13}|B|^{\f23}R^{-2} \tilde{w}^{l_1}_k$, $h_2=v\partial_r \tilde{w}^{l_1}_k$ and $g(r)=1$ yields
\begin{align*}
     &\nu^{\f16}|kB|^{\f56}|k|^{\f12}R^{-2}(\|e^{c'(\nu k^2)^{\f13}|B|^{\f23}R^{-2}t}\partial_r{\varphi}^{l_2}_k\|_{L^2L^2}+|k|\|e^{c'(\nu k^2)^{\f13}|B|^{\f23}R^{-2}t}\f{{\varphi}^{l_2}_k}{r}\|_{L^2L^2}) \\
    \lesssim& \left((\f{\nu}{|kB|})^{\f16} (\log R)^{\f12} +1 \right) \|e^{c'(\nu k^2)^{\f13}|B|^{\f23}R^{-2}t}rh_1\|_{L^2L^2}+\nu^{-\f13} |kB|^{\f13}\|e^{c'(\nu k^2)^{\f13}|B|^{\f23}R^{-2}t}h_2\|_{L^2L^2}.
\end{align*}
Notice that
\begin{align*}
    &\|e^{c'(\nu k^2)^{\f13}|B|^{\f23}R^{-2}t}rh_1\|_{L^2L^2}\\
    =&\|e^{-(c-c')(\nu k^2)^{\f13}|B|^{\f23}R^{-2} t} \left(\nu(k^2-\f{1}{4})\f{w_k(0)}{r}-c(\nu k^2)^{\f13}|B|^{\f23}R^{-2} rw_k(0) \right)\|_{L^2L^2} \\
    \lesssim& \nu k^2 (\nu k^2)^{\f16}|B|^{\f13}R^{-1}\|\f{w_k(0)}{r}\|_{L^2}+(\nu k^2)^{\f16}|B|^{\f13}R^{-1}\| rw_k(0) \|_{L^2},
\end{align*}
and
\begin{align*}
    &\|e^{c'(\nu k^2)^{\f13}|B|^{\f23}R^{-2}t}h_2\|_{L^2L^2}\\
    =& \|e^{-(c-c')(\nu k^2)^{\f13}|B|^{\f23}R^{-2} t} \nu \partial_r\left(e^{-ikB\f{t}{r^2}}w_k(0)\right)\|_{L^2L^2} \\
    \lesssim&  \|e^{-(c-c')(\nu k^2)^{\f13}|B|^{\f23}R^{-2} t} \nu kBt \f{w_k(0)}{r^3}\|_{L^2L^2}+\|e^{-(c-c')(\nu k^2)^{\f13}|B|^{\f23}R^{-2} t} \nu \partial_r w_k(0)\|_{L^2L^2} \\
    \lesssim&  \nu |kB| [(\nu k^2)^{\f13}|B|^{\f23}R^{-2}]^{-\f32}\|\f{w_k(0)}{r^3}\|_{L^2}+\nu (\nu k^2)^{\f16}|B|^{\f13}R^{-1}\|\partial_r w_k(0)\|_{L^2}.
\end{align*}
Thus we arrive at 
\begin{align*}
     &k^2|B|R^{-4}(\|e^{c'(\nu k^2)^{\f13}|B|^{\f23}R^{-2}t}\partial_r{\varphi}^{l_2}_k\|_{L^2L^2}^2+|k|^2\|e^{c'(\nu k^2)^{\f13}|B|^{\f23}R^{-2}t}\f{{\varphi}^{l_2}_k}{r}\|_{L^2L^2}^2) \\
     \lesssim& \left((\f{\nu}{|kB|})^{\f13} \log R +1 \right)\left( R^{-2}\|rw_k(0)\|_{L^2}^2 + (\f{\nu}{|kB|})^{\f43}k^4 R^2\|\f{w_k(0)}{r}\|_{L^2}^2\right)\\
     &+R^6 \|\f{w_k(0)}{r^3}\|_{L^2}^2+(\f{\nu}{|kB|})^{\f23}R^2\|\partial_r w_k(0)\|_{L^2}^2.
\end{align*}
\end{proof}
Therefore, the desired bound \eqref{nonzero frequency linear part for varphi sum} can be readily deduced from  Lemma \ref{nonzero frequency linear part for varphi l1 part} and Lemma \ref{nonzero frequency linear part for varphi l2 part}. This completes the proof of Lemma \ref{nonzero frequency linear part for varphi}.

\section{\textbf{Space-time estimates for the linearized Navier-Stokes equations}}\label{3-space time}
In this section, we establish the space-time estimates for the linearized 2D Navier-Stokes equation in the vorticity formulation \eqref{pertubation of NS vor polar coordinates}:
\begin{align}
\left\{
\begin{aligned}
&\partial_tw-\nu(\partial_r^2+\f{1}{r}\partial_r+\f{1}{r^2}\partial_{\theta}^2)w+(A+\f{B}{r^2})\partial_{\theta}w+\f{1}{r}(\partial_r\varphi\partial_{\theta}w-\partial_{\theta}\varphi\partial_rw)=0,\\
&(\partial_r^2+\f{1}{r}\partial_r+\f{1}{r^2}\partial_{\theta}^2)\varphi=w,\quad w|_{r=1,R}=0\quad \text{with }(r,\theta)\in[1,R]\times\mathbb{S}^1 \text{ and } t\geq0.
\end{aligned}
\right.
\end{align}
Recall in Section \ref{Derivation of the perturbative equation} we convert this equation to the following system
\begin{align}
\label{main nonlinear equation}\left\{
\begin{aligned}
&\partial_tw_k+L_kw_k+\f{1}{r}[ikf_1-r^{\f12}\partial_{r}(r^{\f12}f_2)]=0,\\
&w_k(0)=w_k|_{t=0},\quad w_k|_{r=1,R}=0
\end{aligned}
\right.
\end{align}
via introducing 
\begin{equation}
    \begin{split}\label{def wk phik}
        &w_k(t, r):=r^{\f12}e^{ikAt}\hat{w}_k(t, r)= \f{1}{2\pi}\int_0^{2\pi} r^{\f12}e^{ikAt}w(t, r, \theta)e^{-ik \theta} d\theta, \\ & \varphi_k(t, r):= r^{\f12}e^{ikAt}\hat{\varphi}_k(t, r)=\f{1}{2\pi}\int_0^{2\pi} r^{\f12}e^{ikAt}\varphi(t, r, \theta)e^{-ik \theta} d\theta.
    \end{split}
\end{equation}
 Here
$L_k=-\nu(\partial_r^2-\f{k^2-\f14}{r^2})+i\f{kB}{r^2}$ and $f_1, f_2$ stand for nonlinear forms given by
\begin{align}
\label{nonlinear part}f_1=\sum_{l\in\mathbb{Z}}\partial_r(r^{-\f12}\varphi_l)w_{k-l},\quad f_2=\sum_{l\in\mathbb{Z}}ilr^{-\f32}\varphi_lw_{k-l},
\end{align}
where $\varphi_k$  satisfies
\begin{align*}
&(\partial_r^2-\f{k^2-\f14}{r^2})\varphi_k=w_k,\quad \varphi_k|_{r=1,R}=0,\quad r\in[1,R],\quad t\geq0.
\end{align*}

\subsection{Space-time estimates for non-zero frequency}\label{3.2-non-zero}
\begin{comment}
Recall the full nonlinear equation \eqref{main nonlinear equation} as follows
\begin{align*}
\left\{
\begin{aligned}
&\partial_tw_k+L_kw_k+\f{1}{r}[ikf_1-r^{\f12}\partial_{r}(r^{\f12}f_2)]=0,\\
&w_k(0)=w_k|_{t=0},\quad w_k|_{r=1,R}=0.
\end{aligned}
\right.
\end{align*}
Here
$L_k=-\nu(\partial_r^2-\f{k^2-\f14}{r^2})+i\f{kB}{r^2}$ and $f_1, f_2$ are nonlinear terms as
\begin{align*}
f_1=\sum_{l\in\mathbb{Z}}\partial_r(r^{-\f12}\varphi_l)w_{k-l},\quad f_2=\sum_{l\in\mathbb{Z}}ilr^{-\f32}\varphi_lw_{k-l},
\end{align*}
where $\varphi_k$  satisfying
\begin{align*}
&(\partial_r^2-\f{k^2-\f14}{r^2})\varphi_k=w_k,\quad \varphi_k|_{r=1,R}=0,\quad r\in[1,R],\quad t\geq0.
\end{align*}
\end{comment}

We decompose the solution $w_k$ to \eqref{main nonlinear equation} into two parts. Let $w_k=w_k^l+w_k^{n}$, with $w_k^l$ fulfilling the homogeneous linear equation
\begin{align*}
&\partial_tw_k^l+L_kw_k^l=0,\quad w_k^l(0)=w_k(0),
\end{align*}
and $w_k^{n}$ solving the inhomogeneous linear equation with zero initial data
\begin{align}
\label{nonzero mode nonlinear}
&\partial_tw_k^{n}+L_kw_k^{n}+\f{1}{r}[ikf_1-r^{\f12}\partial_{r}(r^{\f12}f_2)]=0,\quad w_k^{n}(0)=0.
\end{align}
Correspondingly, we also decompose $\varphi_k$ as $\varphi_k=\varphi^{l}_k+\varphi^{n}_k$, where $\varphi_k^l$ and $\varphi_k^n$ satisfy
\begin{equation*}
    (\partial_r^2-\f{k^2-\f14}{r^2})\varphi^{l}_k=w_k^l,\quad \varphi^{l}_k|_{r=1,R}=0,\quad \textrm{with} \ r\in[1,R] \ \textrm{and} \ t\geq0,
\end{equation*}
and 
\begin{equation*}
    (\partial_r^2-\f{k^2-\f14}{r^2})\varphi^{n}_k=w_k^n,\quad \varphi^{n}_k|_{r=1,R}=0,\quad \textrm{with} \ r\in[1,R] \ \textrm{and} \ t\geq0.
\end{equation*}
 Now we state the main conclusion of this section, which is also named as the space-time estimate for equation \eqref{main nonlinear equation}.
\begin{proposition}\label{nonzero nonlinear-original}Assume $w_k$ is a solution to \eqref{main nonlinear equation} with $w_k(0)\in L^2$, then there exists a constant $c'>0$ independent of $\nu,B,k,R$ so that it holds
	\begin{align*}
&\|e^{c'(\nu k^2)^{\f13}|B|^{\f23}R^{-2}t}w_k\|_{L^{\infty}L^2}+(\nu k^2)^{\f16}|B|^{\f13}R^{-1}\|e^{c'(\nu k^2)^{\f13}|B|^{\f23}R^{-2}t}w_k\|_{L^2L^2}\\
&+\nu^{\f12}\|e^{c'(\nu k^2)^{\f13}|B|^{\f23}R^{-2}t}\partial_r w_k\|_{L^2L^2}+(\nu k^2)^{\f12}\|e^{c'(\nu k^2)^{\f13}|B|^{\f23}R^{-2}t}\f{w_k}{r}\|_{L^2L^2}\\
&+|B|^{\f12}R^{-2}\left(|k|\|e^{c'(\nu k^2)^{\f13}|B|^{\f23}R^{-2}t}\partial_r\varphi_k\|_{L^2L^2}+k^2\|e^{c'(\nu k^2)^{\f13}|B|^{\f23}R^{-2}t}\f{\varphi_k}{r}\|_{L^2L^2}\right)\\
\lesssim&\|w_k(0)\|_{L^2}+R^{-2}(\log R)^{-1}\|r^2w_k(0)\|_{L^2}+R^{3}\|\f{w_k(0)}{r^3}\|_{L^2}+(\f{\nu}{|kB|})^{\f13}R\|\partial_rw_k(0)\|_{L^2}\\
&+\left(1+(\f{\nu}{|kB|})^{\f13}\log R\right)^{\f12} \cdot \left(R^{-1}\|rw_k(0)\|_{L^2}+R\|\f{w_k(0)}{r}\|_{L^2}(\f{\nu}{|kB|})^{\f23}k^2\right)\\
&+\left(|kB|^{-\f12} (\log R)^{\f12} +\nu^{-\f16} |kB|^{-\f13} \right) \cdot\|e^{c'(\nu k^2)^{\f13}|B|^{\f23}R^{-2}t}kf_1\|_{L^2L^2}\\
&+\nu^{-\f12} \|e^{c'(\nu k^2)^{\f13}|B|^{\f23}R^{-2}t}f_2\|_{L^2L^2}.
	\end{align*}
\end{proposition}
\begin{proof}
Employing Lemma \ref{nonzero frequency linear part} and Lemma \ref{nonzero frequency linear part for varphi}, we can bound the linear part $w_k^l$ and $\varphi_k^l$:
\begin{align*}
&\|e^{c'(\nu k^2)^{\f13}|B|^{\f23}R^{-2}t}w_k^l\|_{L^{\infty}L^2}+(\nu k^2)^{\f16}|B|^{\f13}R^{-1}\|e^{c'(\nu k^2)^{\f13}|B|^{\f23}R^{-2}t}w_k^l\|_{L^2L^2}\\
&+\nu^{\f12}\|e^{c'(\nu k^2)^{\f13}|B|^{\f23}R^{-2}t}\partial_r w_k^l\|_{L^2L^2}+(\nu k^2)^{\f12}\|e^{c'(\nu k^2)^{\f13}|B|^{\f23}R^{-2}t}\f{w_k^l}{r}\|_{L^2L^2}\\
&+|B|^{\f12}R^{-2}\left(|k|\|e^{c'(\nu k^2)^{\f13}|B|^{\f23}R^{-2}t}\partial_r\varphi_k^l\|_{L^2L^2}+k^2\|e^{c'(\nu k^2)^{\f13}|B|^{\f23}R^{-2}t}\f{\varphi_k^l}{r}\|_{L^2L^2}\right)\\
\lesssim&\|w_k(0)\|_{L^2}+R^{-2}(\log R)^{-1}\|r^2w_k(0)\|_{L^2}+R^{3}\|\f{w_k(0)}{r^3}\|_{L^2}+(\f{\nu}{|kB|})^{\f13}R\|\partial_rw_k(0)\|_{L^2}\\
&+\left(1+(\f{\nu}{|kB|})^{\f13}\log R\right)^{\f12} \cdot \left(R^{-1}\|rw_k(0)\|_{L^2}+R\|\f{w_k(0)}{r}\|_{L^2}(\f{\nu}{|kB|})^{\f23}k^2\right).
\end{align*}
 To obtain the desired bounds for the remainder terms $w_k^n$ and $\varphi_k^n$, we utilize Proposition \ref{spacetime estimate proposition-nonlinear} with the choices of $h_1 = -\frac{1}{r}ikf_1$, $h_2 = -r^{\frac{1}{2}}f_2$, and $g(r) = r^{-\frac{1}{2}}$ to deduce
 \begin{align*}
&\|e^{c'(\nu k^2)^{\f13}|B|^{\f23}R^{-2}t}w_k^n\|_{L^{\infty}L^2}+(\nu k^2)^{\f16}|B|^{\f13}R^{-1}\|e^{c'(\nu k^2)^{\f13}|B|^{\f23}R^{-2}t}w_k^n\|_{L^2L^2}\\
&+\nu^{\f12}\|e^{c'(\nu k^2)^{\f13}|B|^{\f23}R^{-2}t}\partial_r w_k^n\|_{L^2L^2}+(\nu k^2)^{\f12}\|e^{c'(\nu k^2)^{\f13}|B|^{\f23}R^{-2}t}\f{w_k^n}{r}\|_{L^2L^2}\\
&+|B|^{\f12}R^{-2}\left(|k|\|e^{c'(\nu k^2)^{\f13}|B|^{\f23}R^{-2}t}\partial_r\varphi_k^n\|_{L^2L^2}+k^2\|e^{c'(\nu k^2)^{\f13}|B|^{\f23}R^{-2}t}\f{\varphi_k^n}{r}\|_{L^2L^2}\right)\\
\lesssim& \left(|kB|^{-\f12} (\log R)^{\f12} +\nu^{-\f16} |kB|^{-\f13} \right) \cdot\|e^{c'(\nu k^2)^{\f13}|B|^{\f23}R^{-2}t}kf_1\|_{L^2L^2}\\
&+\nu^{-\f12} \|e^{c'(\nu k^2)^{\f13}|B|^{\f23}R^{-2}t}f_2\|_{L^2L^2}.
\end{align*}
 Combining these two estimates, we obtain the desired inequality.
\end{proof}
When $\log R\lesssim\nu^{-\f13}|B|^{\f13}$, the $(\f{\nu}{|kB|})^{\f13}\log R$ term on the right for the above estimate can be eliminated, and thus it infers that
\begin{proposition}\label{nonzero nonlinear-1}Under the same conditions as in Proposition \ref{nonzero nonlinear-original}. If $\log R\lesssim\nu^{-\f13}|B|^{\f13}$, then there holds
	\begin{align*}
&\|e^{c'(\nu k^2)^{\f13}|B|^{\f23}R^{-2}t}w_k\|_{L^{\infty}L^2}+(\nu k^2)^{\f16}|B|^{\f13}R^{-1}\|e^{c'(\nu k^2)^{\f13}|B|^{\f23}R^{-2}t}w_k\|_{L^2L^2}\\
&+\nu^{\f12}\|e^{c'(\nu k^2)^{\f13}|B|^{\f23}R^{-2}t}w_k'\|_{L^2L^2}+(\nu k^2)^{\f12}\|e^{c'(\nu k^2)^{\f13}|B|^{\f23}R^{-2}t}\f{w_k}{r}\|_{L^2L^2}\\
&+|B|^{\f12}R^{-2}(|k|\|e^{c'(\nu k^2)^{\f13}|B|^{\f23}R^{-2}t}\partial_r\varphi_k\|_{L^2L^2}+k^2\|e^{c'(\nu k^2)^{\f13}|B|^{\f23}R^{-2}t}\f{\varphi_k}{r}\|_{L^2L^2})\\
\lesssim&\|w_k(0)\|_{L^2}+R^{-2}(\log R)^{-\f32}\|r^2w_k(0)\|_{L^2}+R^{3}\|\f{w_k(0)}{r^3}\|_{L^2}+(\f{\nu}{|kB|})^{\f13}R\|\partial_rw_k(0)\|_{L^2}\\
&+R\|\f{w_k(0)}{r}\|_{L^2}(\f{\nu k^2}{|B|})^{\f23}+\nu^{-\f16} |kB|^{-\f13}\|e^{c'(\nu k^2)^{\f13}|B|^{\f23}R^{-2}t}kf_1\|_{L^2L^2}\\
&+\nu^{-\f12} \|e^{c'(\nu k^2)^{\f13}|B|^{\f23}R^{-2}t}f_2\|_{L^2L^2}.
	\end{align*}
\end{proposition}
Notice that when $\nu k^2 R^{-2} \ge (\nu k^2)^{\f13}|B|^{\f23}R^{-2}$, i.e., $\nu k^2\ge|B|$, the heat dissipation effect becomes more significant compared with the enhanced dissipation. Therefore, in the regime where $\nu k^2\ge |B|$, the following space-time estimates offer a more precise control of $w_k$. 
\begin{proposition}\label{Dominated by the heat equation}Let $w_k$ be the solution to \eqref{main nonlinear equation} with $w_k(0)\in L^2$. Then there exists a constant $c'>0$ independent of $\nu,B,k,R$, such that 
\begin{align*} 
&\|{e^{c'\nu k^2 R^{-2}t}}w_k\|_{L^{\infty}L^2}+\nu\|{e^{c'\nu k^2 R^{-2}t}}\partial_{r}w_k\|_{L^2L^2}+(\nu k^2)^{\f12}\|{e^{c'\nu k^2 R^{-2}t}}\f{w_k}{r}\|_{L^2L^2}\\
&+(\nu k^2)^{\f12}R^{-2}(|k|\|{e^{c'\nu k^2 R^{-2}t}}\varphi'_k\|_{L^2}+k^2\|{e^{c'\nu k^2 R^{-2}t}}\f{\varphi_k}{r}\|_{L^2})\\
\lesssim&\|w_k(0)\|_{L^2}+\nu^{-\f12}(\|{e^{c'\nu k^2 R^{-2}t}}f_1\|_{L^2L^2} +\|{e^{c'\nu k^2 R^{-2}t}}f_2\|_{L^2L^2}).
	\end{align*}
\end{proposition}
\begin{proof}Conducting the integration by parts, we obtain
	\begin{align*}
	&\Re\langle\partial_tw_k-\nu (\partial_r^2-\f{k^2-\f14}{r^2})w_k+\f{ikB}{r^2}w_k+\f{1}{r}[ikf_1-r^{\f12}\partial_{r}(r^{\f12}f_2)],w_k\rangle\\ =&\f12\partial_t\|w_k\|_{L^2}^2+\nu\|\partial_rw_k\|_{L^2}^2+\nu(k^2-\f14)\|\f{w_k}{r}\|_{L^2}^2-\langle f_1,ik\f{w_k}{r}\rangle+\langle f_2,r^{\f12}\partial_r(r^{-\f12}w_k)\rangle=0.
	\end{align*}
It infers
\begin{align*}
\partial_t\|w_k\|_{L^2}^2+\nu\|\partial_rw_k\|_{L^2}^2+\nu k^2\|\f{w_k}{r}\|_{L^2}^2\lesssim\|f_1\|_{L^2}|k|\|\f{w_k}{r}\|_{L^2} +\| f_2\|_{L^2}\|r^{\f12}\partial_r(r^{-\f12}w_k)\|_{L^2}.
	\end{align*}
By applying Cauchy-Schwarz inequality, we then deduce
\begin{align*}
\partial_t\|w_k\|_{L^2}^2+\nu\|\partial_rw_k\|_{L^2}^2+\nu k^2\|\f{w_k}{r}\|_{L^2}^2\lesssim\nu^{-1}(\|f_1\|_{L^2}^2 +\| f_2\|_{L^2}^2).
	\end{align*}
 Noticing that $\|\f{w_k}{r}\|_{L^2}\ge R^{-1}\|w_k\|_{L^2}$, it follows 
 \begin{equation*}
     \partial_t\|w_k\|_{L^2}^2+\nu\|\partial_rw_k\|_{L^2}^2+\nu k^2R^{-2}\|w_k\|_{L^2}^2+\nu k^2\|\f{w_k}{r}\|_{L^2}^2\lesssim\nu^{-1}(\|f_1\|_{L^2}^2 +\| f_2\|_{L^2}^2).
 \end{equation*}
Therefore, we can multiply $e^{2c'\nu k^2R^{-2}t}$ on both sides of above inequality. With $c'$ being a small constant independent of $\nu,B,k,R$, we obtain
	\begin{align*} &\partial_t\|e^{c'\nu k^2R^{-2}t}w_k\|_{L^2}^2+\nu\|e^{c'\nu k^2R^{-2}t}\partial_{r}w_k\|_{L^2}^2+\nu k^2\|e^{c'\nu k^2R^{-2}t}\f{w_k}{r}\|_{L^2}^2\\\lesssim&\nu^{-1}(\|e^{c'\nu k^2R^{-2}t}f_1\|_{L^2}^2 +\|e^{c'\nu k^2R^{-2}t}f_2\|_{L^2}^2).
	\end{align*}
This further implies
\begin{align*} 
&\|e^{c'\nu k^2R^{-2}t}w_k\|_{L^{\infty}L^2}^2+\nu\|e^{c'\nu k^2R^{-2}t}\partial_{r}w_k\|_{L^2L^2}^2+\nu k^2\|e^{c'\nu k^2R^{-2}t}\f{w_k}{r}\|_{L^2L^2}^2\\
\lesssim&\|w_k(0)\|_{L^2}^2+\nu^{-1}(\|e^{c'\nu k^2R^{-2}t}f_1\|_{L^2L^2}^2 +\|e^{c'\nu k^2R^{-2}t}f_2\|_{L^2L^2}^2).
	\end{align*}
In view of Lemma \ref{Appendix A4}, it also holds
\begin{align*}
R^{-2}(|k|\|\varphi'_k\|_{L^2}+k^2\|\f{\varphi_k}{r}\|_{L^2})\lesssim R^{-2}\|rw_k\|_{L^2}\le\|\f{w_k}{r}\|_{L^2}.
\end{align*}
Combining these two estimates above yields
\begin{align*} 
&\|e^{c'\nu k^2R^{-2}t}w_k\|_{L^{\infty}L^2}+\nu\|e^{c'\nu k^2R^{-2}t}\partial_{r}w_k\|_{L^2L^2}+(\nu k^2)^{\f12}\|e^{c'\nu k^2R^{-2}t}\f{w_k}{r}\|_{L^2L^2}\\
&+(\nu k^2)^{\f12}R^{-2}(|k|\|e^{c'\nu k^2R^{-2}t}\varphi'_k\|_{L^2}+k^2\|e^{c'\nu k^2R^{-2}t}\f{\varphi_k}{r}\|_{L^2})\\
\lesssim&\|w_k(0)\|_{L^2}+\nu^{-\f12}(\|e^{c'\nu k^2R^{-2}t}f_1\|_{L^2L^2} +\|e^{c'\nu k^2R^{-2}t}f_2\|_{L^2L^2}),
	\end{align*}
which completes the proof.
\end{proof}

\subsection{Space-time estimates for zero frequency}\label{3.3-zero}
To establish the space-time estimate for zero mode of solutions, we directly utilize the heat dissipative structure of equation and then perform integration by parts.

Recall that the nonlinear perturbation equation reads
\begin{align}\label{full nonlinear eqn}
&\partial_tw-\nu(\partial_r^2+\f{1}{r}\partial_r+\f{1}{r^2}\partial_{\theta}^2)w+(A+\f{B}{r^2})\partial_{\theta}w+\f{1}{r}(\partial_r\varphi\partial_{\theta}w-\partial_{\theta}\varphi\partial_rw)=0,
\end{align}
where $w|_{r=1,R}=0$ and $(\partial_r^2+\f{1}{r}\partial_r+\f{1}{r^2}\partial_{\theta}^2)\varphi=w$. Denote the zero mode of function $f(r, \theta)$ to be 
$$f_{=}:=\f{1}{2\pi}\int_{0}^{2\pi} f(r, \theta) d\theta $$ and set $f_{\neq}:=f-f_{=}$. We then have that the zero frequency part of equation \eqref{full nonlinear eqn} takes the form
\begin{align}
\label{zero part}\partial_tw_{=}-\nu(\partial_r^2+\f{1}{r}\partial_r)w_{=}+\f{1}{r}(\partial_r\varphi\partial_{\theta}w-\partial_{\theta}\varphi\partial_rw)_{=}=0.
\end{align}
Subsequently, we derive the spacetime estimate for $w_{=}$.
\begin{lemma}\label{zero frequency part}Let $w_{=}$ be the solution to \eqref{zero part} and assume $r^{\f12}w_{=}(0)\in L^2$. Denoting $w_0:=r^{\f12} w_{=}$, then the following inequality holds 
\begin{align}
\label{zero frequency part energy}
\|w_0\|_{L^{\infty}L^2}^2+\nu\|r^{\f12}\partial_r(r^{-\f12}w_0)\|_{L^2L^2}^2\leq\nu^{-1}\|\sum_{l\in\mathbb{Z}\backslash\{0\}}\f{l\varphi_lw_{-l}}{r^{\f32}}\|_{L^2L^2}^2+\|w_0(0)\|_{L^2}^2.
\end{align}
Here $\varphi_l$ and $w_{l}$ are defined in \eqref{def wk phik}.
\end{lemma}
\begin{proof}Employing integration by parts for equation \eqref{zero part}, we deduce
\begin{align*}
0=&\Re\langle\partial_tw_{=}-\nu(\partial_r^2+\f{1}{r}\partial_r)w_{=}+\f{1}{r}(\partial_r\varphi\partial_{\theta}w-\partial_{\theta}\varphi\partial_rw)_{=},rw_{=}\rangle\\
=&\f12\partial_t\|r^{\f12}w_{=}(t)\|_{L^2}^2+\nu\|r^{\f12}\partial_rw_{=}\|_{L^2}^2+\Re\langle [\partial_r(\varphi\partial_{\theta}w)-\partial_{\theta}(\varphi\partial_rw)]_{=},w_{=}\rangle.
\end{align*}
Observe that
\begin{align*}
[\partial_r(\varphi\partial_{\theta}w)]_{=}=\partial_r(\varphi\partial_{\theta}w)_{=} \ \textrm{and} \ [\partial_{\theta}(\varphi\partial_rw)]_{=}=\partial_{\theta}(\varphi\partial_rw)_{=}=0.
\end{align*}
Conducting integration by parts again, we then arrive at
\begin{align*}
&\f12\partial_t\|r^{\f12}w_{=}(t)\|_{L^2}^2+\nu\|r^{\f12}\partial_rw_{=}\|_{L^2}^2+\Re\langle \partial_r(\varphi\partial_{\theta}w)_{=},w_{=}\rangle\\
=&\f12\partial_t\|r^{\f12}w_{=}(t)\|_{L^2}^2+\nu\|r^{\f12}\partial_rw_{=}\|_{L^2}^2-\Re\langle (\varphi\partial_{\theta}w)_{=},\partial_rw_{=}\rangle=0.
\end{align*}
This leads to the following inequality
\begin{align*}
\partial_t\|r^{\f12}w_{=}(t)\|_{L^2}^2+2\nu\|r^{\f12}\partial_rw_{=}\|_{L^2}^2\leq2\|r^{-\f12}(\varphi\partial_{\theta}w)_{=}\|_{L^2}\|r^{\f12}\partial_rw_{=}\|_{L^2},
\end{align*}
which renders
\begin{align*}
\partial_t\|r^{\f12}w_{=}(t)\|_{L^2}^2+\nu\|r^{\f12}\partial_rw_{=}\|_{L^2}^2\leq\nu^{-1}\|\f{(\varphi\partial_{\theta}w)_{=}}{r^{\f12}}\|_{L^2}^2.
\end{align*}
We proceed to integrate the above inequality in $t$ variable and get
\begin{align*}
\|r^{\f12}w_{=}(t)\|_{L^2}^2+\nu\int_0^t\|r^{\f12}\partial_rw_{=}(s)\|_{L^2}^2ds\leq\nu^{-1}\int_0^t\|\f{(\varphi\partial_{\theta}w)_{=}(s)}{r^{\f12}}\|_{L^2}^2ds+\|r^{\f12}w_{=}(0)\|_{L^2}^2.
\end{align*}
In view of the definition for $w_k, \varphi_k$ from \eqref{def wk phik}: 
\begin{equation*}
    \begin{split}
        &w_k(t, r):=r^{\f12}e^{ikAt}\hat{w}_k(t, r)= \f{1}{2\pi}\int_0^{2\pi} r^{\f12}e^{ikAt}w(t, r, \theta)e^{-ik \theta} d\theta, \\ & \varphi_k(t, r):= r^{\f12}e^{ikAt}\hat{\varphi}_k(t, r)=\f{1}{2\pi}\int_0^{2\pi} r^{\f12}e^{ikAt}\varphi(t, r, \theta)e^{-ik \theta} d\theta,
    \end{split}
\end{equation*}
we can further express
\begin{align*}
\int_0^t\|\f{(\varphi\partial_{\theta}w)_{=}(s)}{r^{\f12}}\|_{L^2}^2ds=\|\sum_{l\in\mathbb{Z}\backslash\{0\}}r^{-\f12}\hat{\varphi}_l(il)\hat{w}_{-l}\|_{L^2L^2}^2=\|\sum_{l\in\mathbb{Z}\backslash\{0\}}\f{l\varphi_lw_{-l}}{r^{\f32}}\|_{L^2L^2}^2.
\end{align*}
This completes the proof of \eqref{zero frequency part energy}.
\end{proof}

\section{\textbf{Nonlinear stability}}	
This section is devoted to the proof of our main Theorem \ref{main-transition threshold} using a bootstrap argument. Thanks to the transformation conducted in Section \ref{Derivation of the perturbative equation}, we can focus on  analyzing the  system \eqref{main nonlinear equation}, which is presented in the following form
\begin{align*}
\left\{
\begin{aligned}
&\partial_tw_k+L_kw_k+\f{1}{r}[ikf_1-r^{\f12}\partial_{r}(r^{\f12}f_2)]=0,\\
&w_k(0)=w_k|_{t=0},\quad w_k|_{r=1,R}=0.
\end{aligned}
\right.
\end{align*}
Here
$L_k=-\nu(\partial_r^2-\f{k^2-\f14}{r^2})+i\f{kB}{r^2}$ and $f_1, f_2$ have the nonlinear structure of
\begin{align*}
f_1=\sum_{l\in\mathbb{Z}}\partial_r(r^{-\f12}\varphi_l)w_{k-l},\quad f_2=\sum_{l\in\mathbb{Z}}ilr^{-\f32}\varphi_lw_{k-l}.
\end{align*}
Note that $\varphi_k$  satisfies
\begin{align*}
&(\partial_r^2-\f{k^2-\f14}{r^2})\varphi_k=w_k,\quad \varphi_k|_{r=1,R}=0,\quad  \textrm{with} \ r\in[1,R] \ \textrm{and} \ t\geq0.
\end{align*}
For notational simplicity, we set $\mu_k:=\max\left\{(\nu k^2)^{\f13} |B|^{\f23} R^{-2}, \nu k^2 R^{-2} \right\}$. We then define the initial energy $\mathcal{M}_0(0)=\|w_0(0)\|_{L^2}$ and for $k\in\mathbb{Z}\backslash\{ 0 \}$ we let
\begin{align*}
\mathcal{M}_k(0)=&
\left\{
\begin{aligned}
    &R^{-2}(\log R)^{-1}\|r^2w_k(0)\|_{L^2}+R^{3}\|\f{w_k(0)}{r^3}\|_{L^2}+R\|\partial_rw_k(0)\|_{L^2},\quad \textrm{when} \ \nu k^2\leq|B|  ,\\
&\|w_k(0)\|_{L^2}\quad  \textrm{when} \ \nu k^2\geq|B|.
\end{aligned}
\right.    
\end{align*}
As a result of Proposition \ref{nonzero nonlinear-1} and Proposition \ref{Dominated by the heat equation}, we readily deduce
\begin{proposition}\label{nonzero nonlinear-all cases}For $k\in \mathbb{Z}\backslash \{0 \}$, let $w_k$ be the solution to \eqref{main nonlinear equation} with $w_k(0)\in L^2$. Given $\log R\lesssim\nu^{-\f13}|B|^{\f13}$, then there exists a constant $c'>0$ independent of $\nu,B,k,R$, such that it holds
	\begin{align*}
&\|e^{c'\mu_kt}w_k\|_{L^{\infty}L^2}+\mu_k^{\f12}\|e^{c'\mu_kt}w_k\|_{L^2L^2}+\nu^{\f12}\|e^{c'\mu_kt}w_k'\|_{L^2L^2}+(\nu k^2)^{\f12}\|e^{c'\mu_kt}\f{w_k}{r}\|_{L^2L^2}\\
&+\max\{|B|^{\f12}R^{-2},  (\nu k^2)^{\f12} R^{-2} \}\cdot\left(|k|\|e^{c'\mu_kt}\partial_r\varphi_k\|_{L^2L^2}+k^2\|e^{c'\mu_kt}\f{\varphi_k}{r}\|_{L^2L^2}\right)\\
\lesssim&\mathcal{M}_k(0)+\mu_k^{-\f12}R^{-1}\|e^{c'\mu_kt}kf_1\|_{L^2L^2}+\nu^{-\f12} \|e^{c'\mu_kt}f_2\|_{L^2L^2}.
\end{align*}
\end{proposition}

\begin{comment}
\begin{proposition}\label{nonzero nonlinear-nu k2 geq B}Assume $w_k$ is a solution to \eqref{main nonlinear equation} with $w_k(0)\in L^2$, when $\nu k^2\geq |B|$, then there exists a constant $c'>0$ independent of $\nu,B,k,R$ so that 
\begin{align*} 
&\|e^{c'\nu k^2t}w_k\|_{L^{\infty}L^2}+\nu\|e^{c'\nu k^2t}\partial_{r}w_k\|_{L^2L^2}+(\nu k^2)^{\f12}\|e^{c'\nu k^2t}\f{w_k}{r}\|_{L^2L^2}\\
&+(\nu k^2)^{\f12}R^{-2}(|k|\|e^{c'\nu k^2t}\varphi'_k\|_{L^2}+k^2\|e^{c'\nu k^2t}\f{\varphi_k}{r}\|_{L^2})\\
\lesssim&\|w_k(0)\|_{L^2}+\nu^{-\f12}(\|e^{c'\nu k^2t}f_1\|_{L^2L^2} +\|e^{c'\nu k^2t}f_2\|_{L^2L^2}).
	\end{align*}
\end{proposition}
\begin{proof}This is a direct corollary of Proposition \ref{Dominated by the heat equation}.
\end{proof}
\end{comment}
We proceed to construct and control below the energy functionals, with 
\begin{align*}
E_0=&\|w_0\|_{L^{\infty}L^2},\\
E_k=
&\|e^{c'\mu_kt}w_k\|_{L^{\infty}L^2}+\mu_k^{\f12}\|e^{c'\mu_kt}w_k\|_{L^2L^2}+|B|^{\f12}|k|^{\f32}R^{-2}\|e^{c'\mu_kt}\f{\varphi_k}{r^{\f12}}\|_{L^2L^{\infty}} \quad  \textrm{for} \ |k|\ge 1.  
\end{align*}
The derivation of a priori estimates for $E_k$ crucially relies on the space-time estimates obtained in Section \ref{3-space time}. We start from presenting the bound for the energy at zero frequency:

\begin{lemma}\label{E0} Given $\log R\lesssim\nu^{-\f13}|B|^{\f13}$. We have the following energy inequality
\begin{align*}
E_0
\lesssim&\mathcal{M}_0(0)+\nu^{-\f12}|B|^{-\f12}R^2\sum_{l\in\mathbb{Z}\backslash\{0\}}E_lE_{-l}.
\end{align*}
\end{lemma}
\begin{proof}
Directly applying Lemma \ref{zero frequency part}, we obtain
\begin{align*}
E_0\lesssim&\mathcal{M}_0(0)+\nu^{-\f12}\|\sum_{l\in\mathbb{Z}\backslash\{0\}}\f{l\varphi_lw_{-l}}{r^{\f32}}\|_{L^2L^2}\\
\lesssim&\mathcal{M}_0(0)+\nu^{-\f12}\|\sum_{l\in\mathbb{Z}\backslash\{0\}}|l|\|\f{\varphi_l}{r^{\f12}}\|_{L^2L^{\infty}}\|\f{w_{-l}}{r}\|_{L^{\infty}L^2}\\
\lesssim&\mathcal{M}_0(0)+\nu^{-\f12}|B|^{-\f12}R^2\sum_{l\in\mathbb{Z}\backslash\{0\}}E_lE_{-l},
\end{align*}
which proves Lemma \ref{E0}.
\end{proof}
Next, we turn to control  $E_k$ for $k\neq 0$. We first show that the nonlinear terms $f_1$ and $f_2$ obey the below $L^2$ estimates.

\begin{lemma}\label{f1 f2}For $k\in \mathbb{Z}\backslash\{0 \}$, we have
\begin{align*}
\|f_1\|_{L^2}
\lesssim& \sum_{l\in\mathbb{Z}\backslash\{0, k\}}\|\partial_r(r^{-\f12}\varphi_l)\|_{L^{\infty}}\|w_{k-l}\|_{L^2}+R\|w_0\|_{L^2}\|w_{k}\|_{L^2},\\
\|f_2\|_{L^2}\leq& \sum_{l\in\mathbb{Z}\backslash\{0\}}|l|\|\f{\varphi_l}{r^{\f12}}\|_{L^{\infty}}\|\f{w_{k-l}}{r}\|_{L^2}.
\end{align*}
\end{lemma}
\begin{proof}
With the expressions of nonlinear terms
\begin{align*}
f_1=\sum_{l\in\mathbb{Z}}\partial_r(r^{-\f12}\varphi_l)w_{k-l},\quad f_2=\sum_{l\in\mathbb{Z}}ilr^{-\f32}\varphi_lw_{k-l},
\end{align*}
we get
\begin{align*}
\|f_1\|_{L^2}=&\|\sum_{l\in\mathbb{Z}}\partial_r(r^{-\f12}\varphi_l)w_{k-l}\|_{L^2}\\
\leq& \sum_{l\in\mathbb{Z}\backslash\{0, k\}}\|\partial_r(r^{-\f12}\varphi_l)\|_{L^{\infty}}\|w_{k-l}\|_{L^2}+\|\partial_r(r^{-\f12}\varphi_0)\|_{L^{\infty}}\|w_{k}\|_{L^2}+\|\partial_r(r^{-\f12}\varphi_k)\|_{L^{\infty}}\|w_{0}\|_{L^2} \\
\lesssim& (\f{R}{R-1})^{\f12}\cdot(1+\log R)\cdot\Big(\sum_{l\in\mathbb{Z}\backslash\{0, k\}}\|rw_l\|_{L^{2}}\|w_{k-l}\|_{L^2}+\|rw_0\|_{L^2}\|w_{k}\|_{L^2}+\|rw_k\|_{L^2}\|w_{k}\|_{L^2} \Big)\\
\lesssim& R (\f{R}{R-1})^{\f12}\cdot(1+\log R)\cdot\Big(\sum_{l\in\mathbb{Z}\backslash\{0, k\}}\|w_l\|_{L^{2}}\|w_{k-l}\|_{L^2}+\|w_0\|_{L^2}\|w_{k}\|_{L^2} \Big),
\end{align*}
where in the third line we employ Lemma \ref{Appendix A4} and Lemma \ref{Appendix A1-6}.

And the $L^2$ norm of $f_2$ can be directly controlled by
\begin{align*}
\|f_2\|_{L^2}=\|\sum_{l\in\mathbb{Z}\backslash\{0\}}\f{l\varphi_lw_{k-l}}{r^{\f32}}\|_{L^2}\leq \sum_{l\in\mathbb{Z}\backslash\{0\}}|l|\|\f{\varphi_l}{r^{\f12}}\|_{L^{\infty}}\|\f{w_{k-l}}{r}\|_{L^2}.
\end{align*}
\end{proof}
The above proposition leads to a control of $E_k$ with $k\in \mathbb{Z}\backslash \{0 \}$.
\begin{lemma}\label{E_k k non-zero}
    Given $\log R\lesssim\nu^{-\f13}|B|^{\f13}$. For $k\in \mathbb{Z}\backslash \{0 \}$, we have
    \begin{equation*}
       E_k\lesssim
\mathcal{M}_k(0)
+\nu^{-\f12}|B|^{-\f12}R^2 (\f{R}{R-1})^{\f12}(1+\log R) \sum_{l\in\mathbb{Z}}E_lE_{k-l}.
    \end{equation*}
\end{lemma}
\begin{proof}
Observe the following basic inequality
\begin{align*}
|k|^{\alpha}\leq|l|^{\alpha}+|k-l|^{\alpha} \quad \textrm{for any}\  k, l\in\mathbb{Z} \ \textrm{and any} \ \alpha\in (0, 1].
\end{align*}
Thanks to the definition $$\mu_k:=\max\left\{(\nu k^2)^{\f13} |B|^{\f23} R^{-2}, \nu k^2 R^{-2} \right\},$$
this yields
\begin{equation}\label{basic 1}
    \mu_k \le \mu_l+\mu_{k-l}  \quad \textrm{for any}\  k, l\in\mathbb{Z}
\end{equation}
by considering the scenarios  $\nu k^2<|B|$ and $\nu k^2\ge|B|$ separately. In view of Proposition \ref{nonzero nonlinear-all cases}, Lemma \ref{Appendix A1-1} and Lemma \ref{Appendix A4}, we now have
\begin{equation}\label{Ek f1 f2}
    E_k\lesssim \mathcal{M}_k(0)+\mu_k^{-\f12}R^{-1}\|e^{c'\mu_kt}kf_1\|_{L^2L^2}+\nu^{-\f12} \|e^{c'\mu_kt}f_2\|_{L^2L^2}.
\end{equation}
Utilizing Lemma \ref{f1 f2}, it holds
\begin{equation}\label{f2}
    \begin{split}
        &\mu_k^{-\f12}R^{-1}\|e^{c'(\nu k^2)^{\f13}|B|^{\f23}R^{-2}t}kf_1\|_{L^2L^2} \\
\lesssim&\mu_k^{-\f12}|k|(\f{R}{R-1})^{\f12}(1+\log R)\cdot\\ &\Big(\|w_0\|_{L^{\infty}L^2}\|e^{c'\mu_k t}w_{k}\|_{L^2L^2}+\sum_{l\in\mathbb{Z}\backslash\{0, k\}}\|e^{c\mu_l t}w_l\|_{L^{\infty}L^{2}}\|e^{c'\mu_{k-l}t}w_{k-l}\|_{L^2L^2}\Big)\\
\leq&\mu_k^{-\f12}|k|(\f{R}{R-1})^{\f12}(1+\log R)\cdot\Big(E_0\mu_k^{-\f12} E_k+\sum_{l\in\mathbb{Z}\backslash\{0, k\}}|l|^{-\f12}E_l\mu_{k-l}^{-\f12} E_{k-l}\Big)\\
\leq&\mu_k^{-\f12}|k|R(\f{R}{R-1})^{\f12}(1+\log R)\cdot \\&\Big( (\nu k^2)^{-\f16}|B|^{-\f13}RE_0E_k+\sum_{l\in\mathbb{Z}\backslash\{0, k\}} |l|^{-\f12}(\nu (k-l)^2)^{-\f16}|B|^{-\f13}RE_lE_{k-l}\Big)\\
\leq& \mu_k^{-\f12}|k|(\nu k^2)^{-\f16}|B|^{-\f13}R(\f{R}{R-1})^{\f12}(1+\log R)\Big(E_0E_k+\sum_{l\in\mathbb{Z}\backslash\{0, k\}} E_lE_{k-l}\Big),
    \end{split}
\end{equation}
where in the last line we use the fact 
\begin{equation*}
    |k|\le 2|l||k-l| \quad \textrm{for any}\  k\in\mathbb{Z}, l\in\mathbb{Z}\backslash\{0, k\}.
\end{equation*}
Note that
\begin{equation*}
  \mu_k^{-\f12}|k|(\nu k^2)^{-\f16}|B|^{-\f13}R\le   \nu^{-\f12} |B|^{-\f12}R^2.
\end{equation*}
Hence we obtain
\begin{equation*}
    \mu_k^{-\f12}R^{-1}\|e^{c'(\nu k^2)^{\f13}|B|^{\f23}R^{-2}t}kf_1\|_{L^2L^2}\lesssim \nu^{-\f12} |B|^{-\f12}R^2  (\f{R}{R-1})^{\f12}(1+\log R) \sum_{l\in\mathbb{Z}}E_lE_{k-l}.
\end{equation*}
To estimate the term involving $f_2$ on the right of \eqref{Ek f1 f2}, we apply Lemma \ref{f1 f2} again to deduce
\begin{equation}\label{f2}
\begin{split}
\nu^{-\f12} \|e^{c'\mu_k t}f_2\|_{L^2L^2}
\leq&\nu^{-\f12} \sum_{l\in\mathbb{Z}\backslash\{0\}}|l|\|e^{c'\mu_l t}\f{\varphi_l}{r^{\f12}}\|_{L^2L^{\infty}}\|e^{c' \mu_{k-l} t}\f{w_{k-l}}{r}\|_{L^{\infty}L^2} \\
 \leq&\nu^{-\f12} |B|^{-\f12}R^2\sum_{l\in\mathbb{Z}\backslash\{0\}}E_lE_{k-l}.
\end{split}
\end{equation}
Together with \eqref{Ek f1 f2}, the above two bounds give the desired control for $E_k$ when $|k|\ge 1$.
\end{proof}

Gathering all above estimates, we now conclude our main theorem:

\begin{theorem}\label{transition threshold} Given $\log R\lesssim\nu^{-\f13}|B|^{\f13}$. There exist constants $\nu_0$ and $c_0, C,c'>0$ independent of $\nu,B,R$,  such that if 
\begin{align*}
  \mathcal{M}(0)=\sum_{k\in\mathbb{Z}}\mathcal{M}_k(0)\le c_0\nu^\f12 B^{\f12}R^{-2}(\f{R}{R-1})^{-\f12}(1+\log R)^{-1},\quad 0<\nu\leq \nu_0,
\end{align*}
then the solution $w$ to the system \eqref{scaling nonlinear} exists globally in time and satisfies the decaying and stability estimates
\begin{align*}
\sum_{k\in\mathbb{Z}\backslash\{0\}}\|w_k(t)\|_{L^2}&\le Ce^{-c'\nu^{\f13}B^{\f23}R^{-2}t}\mathcal{M}(0) \quad \textrm{and} \quad
\|w_{0}(t)\|_{L^2}\le C\mathcal{M}(0).
\end{align*}
\end{theorem}
\begin{proof}
Employing Lemma \ref{E0} and Lemma \ref{E_k k non-zero}, we have
\begin{align*}
\sum_{k\in\mathbb{Z}}E_k\le&\sum_{k\in\mathbb{Z}}
C\mathcal{M}_k(0)
+C\nu^{-\f12}|B|^{-\f12}R^2(\f{R}{R-1})^{\f12}(1+\log R)\sum_{k\in\mathbb{Z}}\sum_{l\in\mathbb{Z}}E_lE_{k-l}.
\end{align*}
We perform a bootstrap argument to prove our main theorem. Assume that for a small constant $c_0>0$ there holds
\begin{align}\label{bootstrap assumption}
\sum_{k\in\mathbb{Z}}E_k\leq 2C\mathcal{M}(0)\le 2C c_0 \nu^{\f12} |B|^{\f12}R^{-2}(\f{R}{R-1})^{-\f12}(1+\log R)^{-1},
\end{align}
we then obtain
\begin{align*}
\sum_{k\in\mathbb{Z}}E_k
\leq & C\sum_{k\in\mathbb{Z}}
\mathcal{M}_k(0)
+C\nu^{-\f12}|B|^{-\f12}R^2(\f{R}{R-1})^{\f12}(1+\log R)(\sum_{k\in\mathbb{Z}}E_k)^2\\
\leq&C\sum_{k\in\mathbb{Z}}
\mathcal{M}_k(0)
+2Cc_0\sum_{k\in\mathbb{Z}}E_k.
\end{align*}
Choosing $c_0\leq\f13 C^{-1}$, we arrive at
\begin{align*}
&\sum_{k\in\mathbb{Z}}E_k\leq 1.5 C\sum_{k\in\mathbb{Z}}
\mathcal{M}_k(0),
\end{align*}
which improves the bootstrap assumption \eqref{bootstrap assumption}. By a standard continuity argument, we therefore prove that the system \eqref{scaling nonlinear} admits a global-in-time solution, and its $k$-th mode $w_k$ obeys the desired bounds
\begin{align*}
\sum_{k\in\mathbb{Z}\backslash\{0\}}\|w_k(t)\|_{L^2}&\le Ce^{-c'\nu^{\f13}B^{\f23}R^{-2}t}\mathcal{M}(0),\\
\|w_{0}(t)\|_{L^2}&\le C\mathcal{M}(0).
\end{align*}
This completes the proof of our main theorem.
\end{proof}

\appendix
	\section{Basic Estimates}
\setcounter{equation}{0}
\renewcommand\theequation{A.\arabic{equation}}
 In the appendix, we provide the proof of several elementary but fundamental calculus lemmas which are continually utilized throughout this paper. The following two Sobolev-type inequalities are crucial in the derivation of resolvent estimates in Section \ref{Section Resolvent estimate}.
\begin{lemma}\label{Appendix A1-1}For any $w\in L^2(1,R)$ with $R>1$, it holds
\begin{equation*}
    \|w\|_{L^{\infty}}^2\leq 2\|\f{w}{r^{\f12}}\|_{L^2}\|r^{\f12}w'\|_{L^2}+\f{R}{R-1}\|\f{w}{r^{\f12}}\|_{L^2}^2.
\end{equation*}
Moreover, if $w|_{r=1 ,R}=0$, then
\begin{align*}
&\|w\|_{L^{\infty}}^2\leq 2\|\f{w}{r^{\f12}}\|_{L^2}\|r^{\f12}w'\|_{L^2}.
\end{align*}
\end{lemma}
\begin{proof}
Choose $r_0\in [1, R]$ such that
\begin{equation*}
   (R-1) |\f{w(r_0)}{r_0^{\f12}}|^2\le \|\f{w}{r^{\f12}}\|_{L^2}^2.
\end{equation*}
It then follows
\begin{equation*}
    |w(r_0)|^2\le  \frac{R}{R-1}\|\f{w}{r^{\f12}}\|_{L^2}^2.
    \end{equation*}
And we have
\begin{equation*}
    |w(r)|^2-|w(r_0)|^2=\int_{r_0}^r\partial_s(|w(s)|^2)ds=\int_{r_0}^r (w'(s)\overline{w}(s)+w(s)\overline{w}'(s))ds\leq 2\|\f{w}{r^{\f12}}\|_{L^2}\|r^{\f12}w'\|_{L^2},
\end{equation*}
which gives
\begin{equation*}
    \|w\|_{L^{\infty}}^2\leq 2\|\f{w}{r^{\f12}}\|_{L^2}\|r^{\f12}w'\|_{L^2}+\f{R}{R-1}\|\f{w}{r^{\f12}}\|_{L^2}^2.
\end{equation*}
When $w|_{r=1,R}=0$, it can be inferred that
\begin{align*} |w(r)|^2=&\int_1^r\partial_s(|w(s)|^2)ds=\int_1^r (w'(s)\overline{w}(s)+w(s)\overline{w}'(s))ds\leq 2\|\f{w}{r^{\f12}}\|_{L^2}\|r^{\f12}w'\|_{L^2}.
\end{align*}
\end{proof}

For the weighted quantity $r^{\f12} w$, we have
\begin{lemma}\label{Appendix A1}Let $w\in L^2(1,R), w|_{r=1,R}=0$, where $R>1$. Then the following inequality holds
\begin{align*}
&\|r^{\f12}w\|_{L^{\infty}}^2\leq 4\|r^{\f12}w\|_{L^2}\|r^{\f12}w'\|_{L^2}.
\end{align*}
\end{lemma}
\begin{proof}Due to the vanishment of $w$ on the boundary, it can be deduced 
\begin{align}
\nonumber r|w(r)|^2=&\int_1^r\partial_s(s|w(s)|^2)ds=\int_1^r|w(s)|^2ds+\int_1^rs\partial_s(|w(s)|^2)ds\\
\label{A1}=&\int_1^r|w(s)|^2ds+\int_1^rs(w'(s)\overline{w}(s)+w(s)\overline{w}'(s))ds.
\end{align}
Meanwhile, through integration by part we can write
\begin{comment}
we have
\begin{align}
\nonumber r|w(r)|^2=&\int_1^r\partial_s(s|w(s)|^2)ds=\int_1^r|w(s)|^2ds+\int_1^rs\partial_s(|w(s)|^2)ds\\
\label{A2}=&\int_1^r|w(s)|^2ds+\int_1^rs(w'(s)\overline{w}(s)+w(s)\overline{w}'(s))ds,
\end{align}
which along with \eqref{A1} gives
\begin{align*}
&\int_1^r|w(s)|^2ds+\int_1^rs(w'(s)\overline{w}(s)+w(s)\overline{w}'(s))ds=\int_R^r|w(s)|^2ds+\int_R^rs(w'(s)\overline{w}(s)+w(s)\overline{w}'(s))ds.
\end{align*}
Thus, we get
\begin{align*}
\int_1^r|w(s)|^2ds-\int_R^r|w(s)|^2ds=-\int_1^rs(w'(s)\overline{w}(s)+w(s)\overline{w}'(s))ds+\int_R^rs(w'(s)\overline{w}(s)+w(s)\overline{w}'(s))ds,
\end{align*}
which gives
\begin{align*}
&\int_1^r|w(s)|^2ds+\int_r^R|w(s)|^2ds=\int_1^R|w(s)|^2ds\\
=&-\int_1^rs(w'(s)\overline{w}(s)+w(s)\overline{w}'(s))ds-\int_r^Rs(w'(s)\overline{w}(s)+w(s)\overline{w}'(s))ds\\
=&-\int_1^Rs(w'(s)\overline{w}(s)+w(s)\overline{w}'(s))ds.
\end{align*}
Thus, we obtain
\end{comment}
\begin{align}
\label{A3}\int_1^R|w(s)|^2ds=&-\int_1^Rs(w'(s)\overline{w}(s)+w(s)\overline{w}'(s))ds\leq 2\|r^{\f12}w\|_{L^2}\|r^{\f12}w'\|_{L^2}.
\end{align}
It then follows from \eqref{A1} that
\begin{align*}
\sup_{r\in[1,R]}r|w(r)|^2
\leq&\sup_{r\in[1,R]}\int_1^r|w(s)|^2ds+\sup_{r\in[1,R]}|\int_1^rs(w'(s)\overline{w}(s)+w(s)\overline{w}'(s))ds|\\
\leq&2\|r^{\f12}w\|_{L^2}\|r^{\f12}w'\|_{L^2}+2\|r^{\f12}w\|_{L^2}\|r^{\f12}w'\|_{L^2}=4\|r^{\f12}w\|_{L^2}\|r^{\f12}w'\|_{L^2}.
\end{align*}
\end{proof}

To control the stream function $\varphi_k$ with the vorticity $w_k$, we develop a series of (weighted) elliptic estimates as below.

\begin{lemma}\label{Appendix A4}
Given $|k|\geq 1$. Let $w=(\partial_r^2-\f{k^2-\f{1}{4}}{r^2} )\varphi$ with $\varphi|_{r=1, R}=0$. Then the following inequality holds 
\begin{align*}
&\|\varphi'\|_{L^2}^2+|k|^2\|\f{\varphi}{r}\|_{L^2}^2\lesssim |\langle w, \varphi\rangle | \lesssim |k|^{-2} \|rw\|_{L^2}^2, \\
&\|\varphi'\|_{L^2}^2+|k|^2\|\f{\varphi}{r}\|_{L^2}^2\lesssim |k|^{-1}\|r^{\f12}w\|_{L^1}^2, \\
&\|r^{\f12} \varphi'\|_{L^\infty}+|k|\|r^{-\f12}\varphi\|_{L^\infty}\lesssim (\f{R}{R-1})^{\f12}|k|^{-\f12}\|rw\|_{L^2}.
\end{align*}
\end{lemma}
\begin{proof}
Noting that 
\begin{equation*}
    -\langle w, \varphi\rangle=\|\varphi'\|_{L^2}^2+(k^2-\f{1}{4})\|\f{\varphi}{r}\|_{L^2}^2,
\end{equation*}
we immediately obtain the first inequality. It then follows from Cauchy-Schwarz inequality that
\begin{equation*}
    |k|^2\|\f{\varphi}{r}\|_{L^2}^2\lesssim \|rw\|_{L^2} \|\f{\varphi}{r}\|_{L^2},
\end{equation*}
which implies 
$$|k|^2\|\f{\varphi}{r}\|_{L^2}\lesssim \|rw\|_{L^2},$$
and
\begin{equation*}
    |\langle w, \varphi\rangle |\le \|rw\|_{L^2} \|\f{\varphi}{r}\|_{L^2} \lesssim |k|^{-2}\|rw\|_{L^2}.
\end{equation*}
On the other hand, in view of Cauchy-Schwarz inequality again and using Lemma \ref{Appendix A1-1}, it can be estimated
\begin{align*}
    \|\varphi'\|_{L^2}^2+(k^2-\f{1}{4})\|\f{\varphi}{r}\|_{L^2}^2=&-\langle w, \varphi\rangle\le \|r^{\f12}w\|_{L^1}\|\f{\varphi}{r^{\f12}}\|_{L^\infty} \\ \lesssim& \|r^{\f12}w\|_{L^1}\|\f{\varphi}{r}\|_{L^2}^{\f12}(\|\f{\varphi}{r}\|_{L^2}+\|\varphi'\|_{L^2}^2) \\
    \lesssim& |k|^{-\f12}\|r^{\f12}w\|_{L^1}  \left( \|\varphi'\|_{L^2}^2+(k^2-\f{1}{4})\|\f{\varphi}{r}\|_{L^2}^2 \right)^{\f12},
\end{align*}
which yields the second inequality.

For the last inequality, we observe that
\begin{equation*}
    \|rw\|_{L^2}^2=\|r\varphi''\|_{L^2}^2+2(k^2-\f14)\|\varphi'\|_{L^2}^2+(k^2-\f14)^2\|\f{\varphi}{r}\|_{L^2}^2 .
\end{equation*}
By applying Lemma \ref{Appendix A1-1} and inequalities above, we then prove
\begin{align*}
    \|r^{\f12}\varphi'\|_{L^\infty}^2\lesssim& \|\varphi'\|_{L^2} \|r^{\f12}\partial_r (r^{\f12}\varphi')\|_{L^2}+\f{R}{R-1} \|\varphi'\|_{L^2}^2 \\
    \lesssim& \|\varphi'\|_{L^2}(\|r\varphi''\|_{L^2}+\|\varphi'\|_{L^2})+\f{R}{R-1} \|\varphi'\|_{L^2}^2 
    \lesssim\f{R}{R-1}|k|^{-1}\|rw\|_{L^2}^2,
\end{align*}
and
\begin{equation*}
\begin{split}
    \|r^{-\f12}\varphi\|_{L^\infty}^2\lesssim& \|\f{\varphi}{r}\|_{L^2} \|r^{\f12}\partial_r (r^{-\f12}\varphi) \|_{L^2} 
    \lesssim |k|^{-3} \|rw\|_{L^2}^2.
\end{split}
\end{equation*}
\end{proof}

\begin{comment}

\begin{lemma}\label{Poincare inequality 2}For any $w\in L^2(1,R), \int_{1}^R w dr=0$ with $R>1$, the following Poincar$\acute{e}$ inequality type holds
\begin{align*}
&\|w\|_{L^2}\le 2\|rw'\|_{L^2}.
\end{align*}
\end{lemma}
\begin{proof}
   Denote $F(r)=\int_1^r w(s) ds$, then one obtains $F(0)=F(R)=0$ and $F'=w$. Through integration by parts, we estimate
   \begin{align*}
       \|w\|_{L^2}^2=\int_1^R \overline{w} dF=-\int_1^R F\overline{w}'dr \le \|\f{F}{r}\|_{L^2}\|rw'\|_{L^2}\le 2\|w\|_{L^2}\|rw'\|_{L^2},
   \end{align*}
   where in the last line we utilize Lemma \ref{Poincare inequality} for $F$.
\end{proof}
\end{comment}
To establish the elliptic estimate for the zero mode of stream function, we require the below Poincar$\acute{e}$-type inequalities.
\begin{lemma}\label{Poincare inequality}For any $w\in L^2(1,R), w|_{r=1,R}=0$ with $R>1$, the following Poincar$\acute{e}$ type inequalities hold
\begin{align*}
\|w\|_{L^2}\le 2 \|rw'\|_{L^2} \quad \textrm{and} \quad \|\f{w}{r^{\f12}}\|_{L^2}\le 2\log R\|r^{\f12}w'\|_{L^2}.
\end{align*}
\end{lemma}
\begin{proof}
    Via integration by parts, we immediately obatin
    \begin{align*}
       \|w\|_{L^2}^2= \int_{1}^{R} |w|^2 dr=&-\int_1^R r(w'\overline{w}+w\overline{w}') dr \le 2\|rw'\|_{L^2}\|w\|_{L^2}.
    \end{align*}
    and 
    \begin{align*}
        \|\f{w}{r^{\f12}}\|_{L^2}^2=\int_{1}^{R} |w|^2 d(\log r)=-\int_1^R \log r(w'\overline{w}+w\overline{w}') dr\le \log R\|r^{\f12}w'\|_{L^2}\|\f{w}{r^{\f12}}\|_{L^2}.
    \end{align*}
\end{proof}
Consequently, we obtain the following elliptic estimate which involves the outer radius $R$.
\begin{lemma}\label{Appendix A1-6}
Let $w=(\partial_r^2+\f{1}{r}\partial_r)\varphi$ with $\varphi|_{r=1, R}=w|_{r=1, R}=0$. The following inequality holds 
\begin{align*}
&\|\varphi'\|_{L^{\infty}} \lesssim (\f{R}{R-1})^{\f12}(1+\log R)\|r^{\f32}w\|_{L^2}.
\end{align*}
\end{lemma}
\begin{proof}
 Notice that
 \begin{align*}
    -\langle w, r\varphi \rangle=-\langle  \varphi''+\f{1}{r}\varphi', r\varphi \rangle = \|r^{\f12}\varphi'\|_{L^2}.   
 \end{align*}
 In view of Lemma \ref{Poincare inequality}, it can be inferred
 \begin{equation*}
     \|r^{\f12}\varphi'\|_{L^2}^2\le \|r^{\f32}w\|_{L^2}\|\f{\varphi}{r^{\f12}}\|_{L^2} \lesssim \log R \|r^{\f32}w\|_{L^2}\|r^{\f12} \varphi'\|_{L^2},
 \end{equation*}
 which yields
 \begin{equation*}
     \|r^{\f12}\varphi'\|_{L^2}+\|\f{\varphi}{r^{\f12}}\|_{L^2}\lesssim \log R\|r^{\f32}w\|_{L^2},
 \end{equation*}
 and
 \begin{equation*}
     \|r^{\f32}\varphi''\|_{L^2}\le \|r^{\f32}w\|_{L^2}+\|r^{\f12}\varphi'\|_{L^2}\lesssim (1+\log R)\|r^{\f32}w\|_{L^2}.
 \end{equation*}
 Applying Lemma \ref{Appendix A1-1}, we derive the $L^\infty$ bound for $\varphi'$ as below
 \begin{align*}
     \|\varphi'\|_{L^\infty}^2\lesssim& \|r^{-\f12}\varphi'\|_{L^2}\|r^{\f12}\varphi''\|_{L^2}+\f{R}{R-1}\|r^{-\f12}\varphi'\|_{L^2}^2\\\lesssim& \left(\f{R}{R-1}+\log R(1+\log R)\right)\|r^{\f32}w\|_{L^2}^2\lesssim \f{R}{R-1}(1+\log R)^2\|r^{\f32}w\|_{L^2}^2.
 \end{align*}
 This completes the proof of this lemma.
\end{proof}
Following  is a collection of calculus inequalities which provide quantitative information of the interval $[1, R]\cap \{r> 0: \  |1-\lambda r^2|\le \delta \}$ where $\delta>0$ is sufficiently small. These estimates play a significant role in the proof of resolvent estimates in Section \ref{Section Resolvent estimate}. 

\begin{lemma}\label{Appendix A1-2}For any $\lambda\in[\f{1}{R^2},1]$ and $0<\delta\ll1$ being a small constant, it holds
\begin{align*}
\int_{r_{-}}^{r_{+}}\f{1}{r}dr\lesssim\delta,\quad \quad\quad\int_{r_{-}-\f{\delta }{\sqrt{\lambda}}}^{r_{+}+\f{\delta }{\sqrt{\lambda}}}\f{1}{r}dr\lesssim\delta,
\end{align*}
where $r_{-}\in(\sqrt{\f{1-\delta}{\lambda}}-\f{\delta }{\sqrt{\lambda}},\sqrt{\f{1-\delta}{\lambda}})$ and $r_{+}\in(\sqrt{\f{1+\delta}{\lambda}},\sqrt{\f{1+\delta}{\lambda}}+\f{\delta }{\sqrt{\lambda}})$.
\end{lemma}
\begin{proof}
A direct computation yields
\begin{align*}
&\int_{r_{-}}^{r_{+}}\f{1}{r}dr=\ln\f{r_{+}}{r_{-}}\leq\ln\f{\sqrt{1+\delta}+\delta}{\sqrt{1-\delta}-\delta}\\
=&\f12\Big[\ln(1+\delta+2\delta\sqrt{1+\delta}+\delta^2)-\ln(1-\delta-2\delta\sqrt{1-\delta}+\delta^2)\Big]\lesssim\delta
\end{align*}
and
\begin{align*}
&\int_{r_{-}-\f{\delta }{\sqrt{\lambda}}}^{r_{+}+\f{\delta }{\sqrt{\lambda}}}\f{1}{r}dr=\ln\f{r_{+}+\f{\delta }{\sqrt{\lambda}}}{r_{-}-\f{\delta }{\sqrt{\lambda}}}\leq\ln\f{\sqrt{1+\delta}+2\delta}{\sqrt{1-\delta}-2\delta}\\
=&\f12\Big[\ln(1+\delta+4\delta\sqrt{1+\delta}+4\delta^2)-\ln(1-\delta-4\delta\sqrt{1-\delta}+4\delta^2)\Big]\lesssim\delta.
\end{align*}
\end{proof}

The $L^2$ norm of $r^{-\f12}(1-\lambda r^2)^{-1}$ in the region $\{r\in [1, R]: \  |1-\lambda r^2|> \delta \}$ can be bounded in terms of $R$ and $\delta$.
\begin{lemma}\label{Appendix A6}
For any $\delta >0, R>1 $ and $\lambda\in \mathbb{R}$. Let $E=\{r> 0: \  |1-\lambda r^2|\le \delta \}$ and $E^c=(0, \infty)\backslash E$. Then it holds
\begin{equation}\label{Appendix A6 ineq}
    \int_{[1, R] \cap E^c} \f{dr}{r(1-\lambda r^2)^2} \lesssim \log R+ \delta^{-1}.
\end{equation}
\end{lemma}
\begin{proof}
If $\lambda\le 0$, then one has
\begin{equation*}
     \int_{[1, R] \cap E^c} \f{dr}{r(1-\lambda r^2)^2} \le \int_1^R \f{dr}{r}=\log R.
\end{equation*}
Now we assume $\lambda>0$. Denote $r_{\pm}=\sqrt{\f{1\pm\delta}{\lambda}}$. We consider the following five scenarios:
\begin{enumerate}
    \item If $1\le r_-<r_+\le R$, then a direct calculation implies
 \begin{equation}\label{Appendix A6 ineq part1}
     \begin{split}
         \int_1^{r_-} \f{dr}{r(1-\lambda r^2)^2}=\f{1}{2}& \int_{\delta}^{1-\lambda} \f{ds}{(1-s)s^2} 
         =\f12\left(\log|\f{s}{s-1}|-\f{1}{s} \right)\Big|^{1-\lambda}_{\delta}\\
         \lesssim& \log(\lambda^{-1}-1)+\delta^{-1} 
         \lesssim \log R+ \delta^{-1}.
     \end{split}
 \end{equation}
    Here in the first line we use the change of variables $s=1-\lambda r^2$ and the last line follows from the condition $r_{-}\in[1, R]$.
    
    Similarly, taking into account the condition $r_+\in[1, R]$, we deduce
    \begin{equation}\label{Appendix A6 ineq part2}
         \int_{r_+}^{R} \f{dr}{r(1-\lambda r^2)^2}  \lesssim \log R+ \delta^{-1}.
    \end{equation}
    Combining these two inequalities above, we obtain the desired result.
    \item If $1\le r_-\le R<r_+$,  then \eqref{Appendix A6 ineq} directly follows from \eqref{Appendix A6 ineq part1} and the fact that $[1 ,R]\cap E^c=[1, r_-]$.
    \item If $r_-<1\le r_+\le R$, then \eqref{Appendix A6 ineq} can be readily obtained thanks to \eqref{Appendix A6 ineq part2} and the fact that $[1 ,R]\cap E^c=[r_+, R]$.
    \item If $r_+<1$ or $r_->R$, then  $[1, R]\cap E^c=[1, R]$. Assume first $r_->R$, it is hence clear that $1-\lambda R^2\ge \delta$. A suitable change of variables implies
     \begin{equation*}
     \begin{split}
         \int_1^{R} \f{dr}{r(1-\lambda r^2)^2}=\f{1}{2}& \int_{1-\lambda R^2}^{1-\lambda} \f{ds}{(1-s)s^2} 
         =\f12\left(\log(\f{s}{s-1})-\f{1}{s} \right)\Big|^{1-\lambda}_{1-\lambda R^2}\\
         \lesssim& \log(R^2\f{1-\lambda}{1-\lambda R^2})+\f{1}{1-\lambda R^2}
         \lesssim \log R+ \delta^{-1}.
     \end{split}
 \end{equation*}
    The case when $r_+<1$ can be treated in a similar manner.
    \item If $r_-<1<R<r_+$, then immediately we see \eqref{Appendix A6 ineq} holds since $[1, R]\cap E^c=\varnothing$.
\end{enumerate}
\end{proof}
Finally, in order to prove the resolvent estimates for $\varphi'$ in Lemma \ref{resolvent estimate-4}, it is necessary to construct a piecewise $C^1$ cut-off function satisfying the following properties.
\begin{lemma}\label{Appendix A5}
For any $0<\delta\ll 1$ and $\lambda\in \mathbb{R}$. Let $E=\{r> 0: \  |1-\lambda r^2|\le \delta \}$ and $E^c=(0, \infty)\backslash E$. Then there exists a piecewise $C^1$ function $\chi: (0, \infty) \to \mathbb{R}$ satisfying the following properties:
\begin{enumerate}
    \item $\chi(r)=\f{1}{\f{1}{r^2}-\lambda}$ on $E^c$.
    \item there exists a constant $C>0$ independent of $\lambda$ and $\delta$, such that it holds \[\|\f{\chi}{r^2}\|_{L^\infty(0, \infty)}+\delta^{\f12}\|r^{\f12}(\f{\chi}{r^2})'\|_{L^2(0, \infty)}\le C\delta^{-1}.\]
    \item $\|\f{1}{r}\|_{L^1(E)}\lesssim \delta, \ \|(1-\lambda r^2)\f{\chi}{r^{\f{5}{2}}}\|_{L^2(E)}\lesssim \delta^{\f12}$.
\end{enumerate}
\end{lemma}
\begin{proof}
If $\lambda\le 0$, then we have $E=\varnothing$, since $|1-\lambda r^2|\ge1$. Setting $\chi(r)=\f{1}{\f{1}{r^2}-\lambda}$, one can directly check that $(1)$ and $(2)$ are satisfied. Now we assume $\lambda>0$. Denote $r_{\pm}=\sqrt{\f{1\pm\delta}{\lambda}}$. We can see that $E=(r_-, r_+)$ and $m(E)\approx \delta \lambda^{-\f12}$. 

Choose the piecewise $C^1$ function $\chi: (0, \infty) \to \mathbb{R}$ as below:
\begin{align*}
\chi(r)=\left\{
\begin{aligned}
&\f{1}{\f{1}{r^2}-\lambda},\quad r\in E^c,\\
&-\f{r_+^2+r_-^2}{\delta (r_+-r_-)}(r-r_-)+\f{r_-^2}{\delta},\quad r\in E.
\end{aligned}
\right.
\end{align*}
Observing that $\|\chi\|_{L^\infty(E)}\lesssim \delta^{-1}\lambda^{-1}$, $\|\chi'\|_{L^\infty(E)}\lesssim \delta^{-2}\lambda^{-\f12}$, and $r\approx \lambda^{-\f12}$ on $E$, we then verify 
\begin{equation*}
    \|\f{1}{r}\|_{L^1(E)}\lesssim \delta, \ \|(1-\lambda r^2)\f{\chi}{r^{\f{5}{2}}}\|_{L^2(E)}\lesssim \delta^{\f12},
\end{equation*}
and
\begin{equation*}
    \|\f{\chi}{r^2}\|_{L^\infty(E)}+\delta^{\f12}\|r^{\f12}(\f{\chi}{r^2})'\|_{L^2(E)}\lesssim \delta^{-1}.
\end{equation*}
It remains to establish
\begin{equation*}
    \|\f{\chi}{r^2}\|_{L^\infty(E^c)}+\delta^{\f12}\|r^{\f12}(\f{\chi}{r^2})'\|_{L^2(E^c)}\lesssim \delta^{-1}.
\end{equation*}
The first term on the left can be easily estimated due to the fact that $\f{\chi}{r^2}=\f{1}{1-\lambda r^2}$ on $E^c=\{r> 0: \  |1-\lambda r^2|> \delta \}$.

For the second term, a direct computation yields
\begin{align*}
\|r^{\f12}(\f{\chi}{r^2})'\|_{L^2(E^c)}^2= \int_{E^c} \f{4\lambda^2 r^3}{(1-\lambda r^2)^4} dr 
=\int_{\delta}^1 \f{2(1-y)}{y^4}dy+\int_{-\infty}^{-\delta} \f{2(1-y)}{y^4}dy  \lesssim \delta^{-3},
\end{align*}
which gives the desired inequality.
\end{proof}

\section*{Acknowledgments}
XA is supported by MOE Tier 1 grants A-0004287-00-00, A-0008492-00-00 and MOE Tier 2 grant A-8000977-00-00.  TL is supported by MOE Tier 1 grant A-0004287-00-00.

\end{CJK*}

\bigskip


\begin{thebibliography}{99}

\bibitem{AHL-1}X. An, T. He and T. Li, {\it Enhanced dissipation and nonlinear asymptotic stability of the Taylor-Couette flow for the 2D Navier-Stokes equations}, arXiv:2112.15357v1.

\bibitem{BGM-bams}J. Bedrossian, P. Germain and N. Masmoudi, {\it Stability of the Couette flow at high Reynolds number in
2D and 3D}, Bull. Amer. Math. Soc. (N.S.), 56 (2019), 373-414.

\bibitem{BH}J. Bedrossian and S. He, {\it Inviscid damping and enhanced dissipation of the boundary
layer for 2D Navier-Stokes linearized around Couette flow in a channel}, Comm. Math. Phys.
379 (2020), no. 1, 177–226.

\bibitem{BMV} J. Bedrossian, N. Masmoudi and V. Vicol,  {\it Enhanced dissipation and inviscid damping in the inviscid limit of the Navier-Stokes equations near the two dimensional Couette flow},  Arch. Ration. Mech. Anal., 219(2016), 1087-1159.

\bibitem{BWV}J. Bedrossian, F. Wang and V. Vicol, {\it The Sobolev stability threshold for 2D shear flows near Couette}, J. Nonlinear Sci.,  28 (2018),  2051-2075.

\bibitem{Cha} S. J. Chapman, {\it Subcritical transition in channel flows},
            J.  Fluid Mech., 451(2002), 35-97.


\bibitem{CLWZ-2D-C}Q. Chen, T. Li, D. Wei and Z. Zhang, {\it  Transition threshold for the 2-D Couette flow in a finite channel}, Arch. Ration. Mech. Anal., 238 (2020), 125-183.


\bibitem{CI}P. Chossat, G. Iooss, {\it The Couette-Taylor Problem}, 
Applied Mathematical Sciences 102, Springer-Verlag, New York, 1994.

\bibitem{DHB}F. Daviaud, J. Hagseth and P. Berg$\acute{e}$, {\it  Subcritical transition to turbulence in plane Couette flow}, Phys. Rev. Lett., 69(1992), 2511-2514.

\bibitem{DW}P. Drazin and W. Reid, {\it  Hydrodynamic Stability}, Cambridge Monographs Mech. Appl. Math., Cambridge Univ. Press, New York, 1981.



\bibitem{F}G. Falkovich, {\it Fluid Mechanics}, Cambridge University Press, 2018.



\bibitem{Ga0}T. Gallay, {\it  Stability and interaction of vortices in two-dimensional viscous flows}, Discr. Cont. Dyn. Systems Ser. S 5 (2012), 1091-1131.

\bibitem{Ga}T. Gallay, {\it  Enhanced dissipation and axisymmetrization of two-dimensional viscous vortices}, Archive for Rational Mechanics and Analysis, volume 230 (2018), 939-975.

\bibitem{GV}T. Gallay and V. $\check{S}$ver$\acute{a}$k, {\it  Arnold's variational principle and its application to the stability of planar vortices}, arXiv:2110.13739.

\bibitem{GW}T. Gallay and C. E. Wayne, {\it  Global stability of vortex soluotions of the two-dimensional Navier-Stokes equation}, Comm. Math. Phys., 255 (2005), 97-129.

\bibitem{HS}B. Helffer and J. Sj\"{o}strand, Improving semigroup bounds with resolvent estimates, Integral Equations Operator Theory 93
(2021), no. 3, Paper No. 36, 41.

\bibitem{K}Robert H. Kraichnan, {\it Inertial Ranges in Two-Dimensional Turbulence}, The Physics
of Fluids, 10 (1967), 1417-1423.

\bibitem{Kato}T. Kato, \textit{Perturbation Theory for Linear Operators}, Grundlehren der Mathematischen Wissenschaften 132. Berlin: Springer, 1966.





\bibitem{Kel}L. Kelvin, {\it Stability of fluid motion-rectilinear motion of viscous fluid between two parallel plates}, Phil. Mag., 24(1887), 188-196.



\bibitem{LWZ-O} T. Li, D. Wei and Z. Zhang, {\it Pseudospectral and spectral bounds for the Oseen vortices operator}, Annales Scientifiques of Ecole Normale Sup\'erieure, 53 (2020),  993-1035.




\bibitem{LHR}A. Lundbladh, D. Henningson and S. Reddy, {\it Threshold amplitudes for transition in channel flows}, in Transition, Turbulence and Combustion, Springer-Verlag, New York, 1994, pp. 309-318.

\bibitem{L}T. Lundgren, {\it Strained spiral vortex model for turbulent fine structure}, Physics of Fluids, 25 (1982), 2193-2203.


\bibitem{MHA}S. Maretzke, B. Hof and M. Avila, {\it Transient growth in linearly stable Taylor-Couette flows}, Journal of Fluid Mechanics , 742(2014), 254-290.

\bibitem{MZ-1}N. Masmoudi and W. Zhao, {\it Stability threshold of two-dimensional Couette flow in
Sobolev spaces}, Ann. Inst. H. Poincar$\acute{e}$ C Anal. Non Lin$\acute{e}$aire 39 (2022), no. 2, 245–325.


\bibitem{MZ-2}N. Masmoudi and W. Zhao, {\it Enhanced dissipation for the 2D Couette flow in critical space}, 
Communications in Partial Differential Equations, 45 (2020), no. 12, 1682-1701.

\bibitem{MVGL}R.O. M$\acute{o}$nico, R. Verzicco, S. Grossmann and D. Lohse, {\it Turbulence decay towards the linearly stable regime of Taylor-Couette flow}, Journal of Fluid Mechanics , 748 (2014), R3.

\bibitem{Pa}A. Pazy, \textit{Semigroups of linear operators and applications to partial differential equations}, Applied Mathematical Sciences 44, Springer-Verlag, New York, 1983.

\bibitem{Po}B. Pope,{\it Turbulent Flows}, Cambridge University Press, 2000.



\bibitem{RSBH}S. Reddy, P. Schmid, J. Baggett and D. Henningson, {\it On stability of streamwise streaks and transition thresholds in plane channel flows}, J. Fluid Mech., 365(1998), 269-303.

\bibitem{Rey}O. Reynolds, {\it An experimental investigation of the circumstances which determine whether the motion of water shall be direct or sinuous, and of the law of resistance in parallel channels}, Proc. R. Soc. Lond., 35(1883), 84.

\bibitem{RY}P. B. Rhines and W. R. Young, {\it How rapidly is a passive scalar mixed within closed streamlines?}, J. Fluid Mech., 133 (1983).

\bibitem{SH}P. Schmid and D. Henningson, {\it Stability and Transition in Shear Flows}, Applied Mathematical Sciences 142, Springer-Verlag, New York, 2001.

\bibitem{Se}J. Serrin, {\it On the stability of viscous fluid motions}, Arch. Rational Mech. Anal. 3 (1959), 1–13.


\bibitem{TA}N. Tillmark and P. H. Alfredsson, {\it Experiments on transition in plane Couette flow}, J. Fluid Mech., 235(1992), 89-102.



\bibitem{TTRD}L. Trefethen, A. Trefethen, S. Reddy and T. Driscoll, {\it Hydrodynamic stability without eigenvalues}, Science, 261(1993), 578-584.


\bibitem{Wei}D. Wei, \textit{Diffusion and mixing in fluid flow via the resolvent estimate}, preprint.

\bibitem{WZ}D. Wei and Z. Zhang, {\it Enhanced dissipation for the Kolmogorov flow via the hypocoercivity method}, Sci. China Math. 62 (2019), no. 6, 1219–1232.

\bibitem{WZZ}D. Wei, Z. Zhang, and W. Zhao, {\it Linear inviscid damping and enhanced
dissipation for the Kolmogorov flow}, Adv. Math. 362 (2020), 106963, 103.

\bibitem{Ya}A. Yaglom, {\it Hydrodynamic instability and transition to turbulence}, Fluid Mech. Appl. 100, Springer-Verlag, New York, 2012.




\end{thebibliography}
\end{document}